# Bridging Algebra and Nature: Toward a Deformable 3D Hyper-complex framework for Modeling Dynamic Systems

Abdon Atangana

Institute for Groundwater Studies, Faculty of Natural and Agricultural Science, University of the Free State, 9300, Bloemfontein, South Africa. Email: atanganaa@ufs.ac.za

**Abstract**

In this paper, we present a new hyper-complex number system, Trinition, that has an unusual structure of commutativity, non-commutativity, non-associativity, and deformability. Challenging traditional algebra, it is generative and degenerative, so it is able to form a dynamic construction that can self-regulate. We are able to extend basic mathematical primitives, such as a norm and distance metric, to a reversible and deformable definition, ultimately creating an entirely new geometry. This then leads to new 3D geometrical figures, such as a deformable sphere and isosurface representations, yielding a new family of 3D fractal system described using recursive Trinition transformations. We then prove a new collection of inequalities which also generalize classical results to this hyper-complex scenario. Along with geometric configurations, we derive new velocity, acceleration, energy functions and respective Laplace and Fourier transforms for a wider range of physics and engineering applications. Using the process of extending classical mechanics into Trinition space, we find the new Hamiltonian and Lagrangian formalism associated with Trinition mechanics, which produces a generalized equation of motion that incorporates hyper-complex interactions. To relate these new evolution stages, we now introduce Archiometry, a new form of hyper-complex trigonometry in which the angles, which we call Archiotwists, define the rotational transformations in Trinition space. The culmination of these results outlines a novel mathematical background that integrates algebra, geometry, physics, and analysis all within the same frame, which may pave the way toward future quantum mechanics developments, complex systems modeling, and room temperature higher-dimension physics applications. We have reformulated spacetime equation, which led to a formula that could open doors for more investigation and our understanding of nature. We have presented a new norm that deforms the geometry by injecting deformation parameters directly into the inner product, resulting in a deformation Laplace-Beltrami operator. This operator in return, yield generalized partial differential equations from which anisotropic and non-Euclidian entities may be appropriately extracted. We prove the existence, uniqueness, and stability of the solutions in the corresponding inequalities with deformation effects. In particular, numerical simulations in MATLAB clearly show how a simple modification of the deformation parameters can lead to significant changes in wave propagation, illustrating the power of this approach to model complex systems where classical convex Euclidean metrics fail to reproduce the system dynamics. Our results pave the way for future studies in anisotropic media, higher-order PDE model and new probabilistic structures.

Contents

**List of tables**

**List of figures**

# Introduction

The development of mathematics has historically been motivated by the desire to reflect the complexities of the natural world [1]. From real numbers describing scalar quantities to complex numbers representing planar rotations, every expansion has brought new worlds of meaning [2,3]. But in this encounter with the modern world of quantum entanglement [4], curved space-time [5], and high-dimensional data, we find ourselves facing an important limitation: the current mathematical framework is unable to unify algebraic structure, geometric intuition, and the fundamental three-dimensionality of physical reality [6].

The dictatorship of commutativity and associativity: Classical hypercomplex systems like quaternions [7] and octonions brought arithmetic beyond complex numbers but added rigid constraints. Quaternions, for example, give up commutativity for the goal of modeling 3D rotations [8], and octonions give up commutativity and associativity altogether [9]. However, these trade-offs, while revolutionary, make them unsuitable for systems with sequence-dependent or contextual interactions situations (biochemical pathways, where the outcomes of reactions depend on the ordering of molecular collisions, or quantum systems, where the results of a measurement depend on the order of past and future measurements) where interactivity is actually time-ordered or interacting in a contextually defined way. Forcing unnatural symmetries on this procedure by insisting that all traditional algebraic properties be preserved limits the capacity of these systems to model adaptive and often chaotic interactions found in nature.

Deformation, the road to flexible mathematics: Deformation theory [10], the study of how algebraic structures deform under parameterized changes, provides an encouraging escape from that rigidity. These parameters were adjusted until they became more closely aligned with what was observed, and in this sense, they adjust the rules of a mathematical context to fine-tune how signs of phenomena can be interpreted, like a lens that is adjusted to focus more accurately. This has proven fruitful in quantum groups where the q-deformations [12] independently deform the commutativity relations so as to encode particle statistics. But these deformations are still abstract and are divorced from the actual geometric structure of physical space. A single system that can deform both algebraic and geometric features would help fill this gap, allowing mathematicians and scientists to tune their models to match experimental data or emerging theoretical ideas.

The 3D Imperative [13] and beyond: While our universe is inherently 3d, most hypercomplex systems overshoot into 4d (quaternions) or undershoot into 2d (complex numbers) [11]. This mismatch pushes researchers to compress their 3D phenomena into non-matching alien dimensional spaces with distortions and artifacts resulting. However, a native 3D hyper-complex system one which treats all axes as interdependent partners rather than isolated dimensions would allow the mathematical tools to be matched to the physical world's intrinsic geometric axioms. Such a system would not only make calculations simpler but also would uncovers hidden patterns that are hidden by dimensional mismatch, from protein folding pathways, to the subtle warping of space-time.

Forthcoming directions, dynamic and circumstantial mathematics: Now imagine a number system where numbers embed natively in 3D [13], whose algebraic rules deform via parameters to emulates dynamic-tightening, and where multiplication couples geometry to algebra so that distances and angles warp in tandem. Such data-driven strategies can elucidate the emergent phenomena in complex systems from addressing non-commutative observables in quantum mechanics, optimization of high-dimensional structures in the data science to replicate emergent characteristics in materials science. Instead, by celebrating all that, every hypercomplex relationship would accept the reality of asymmetry, and, because it would emerge from the world, the new paradigm would operate by not only describing reality but adapting with it, allowing us an innovative possibility for decoding the greatest mysteries the universe has in store. This work suggests such a system, trained to go beyond the bounds of conventional hypercomplex numbers without sacrificing its analytic strength. By introducing a new approach that combines algebraic creativity, intuitive geometry, and deformation methods, we pioneer the development of a mathematics that grows in parallel to the systems it attempts to describe, and introduces a new level of precision, flexibility, and insight into complex systems.

# Motivation and introduction of new hypercomplex number

The obstinate faith in commutativity and associativity whose basis is rooted in the relatively unambiguous forms of classical mathematics prevents us from seeing that reality is a market place of non-commuting transactions between quantum fields. Whose arcane choreography is dominated by creation and annihilation operators, and whose ordered movements decide the birth of particles, let alone the entanglement of space-time itself. Gauge theories [15] made from non-Abelian Lie-algebras [14,15] that weave the cords of the standard model, in which the system itself trades singlets for triplets. In addition, the notion of commutativity is an illusion and colloquies of energy transmutations bind quarks through gluons, which dance in space-time. The black holes [16], those anti-entropies [17], tell us that holographic dualities exist, and that our monolithic concepts of associative locality tear under the tension of quantum error-correcting codes. The octonions of non-associativity haunt M-theory's dream-space and the geometry of forces begins to dissolve, while biology speaks the subversive language of quantum tunneling and we learn that we are indulgers in the illegal economy of subatomic crime. The neural nets [18] operate on deeply-pathological classical code, and do not perceive the deep quantum art of entanglement that is the soul; even time itself, that macroscopic fantasy cast by the macroscope, unravels and with it our sense of quantum retro-causation, so that observables can only commute in a narrow band of causal symmetry. While AI, hobbled by associative Boolean logic, stumbles toward AGI, while failing to realize the braided syntax of topological qubits and the contextuality of quantum cognition. Moreover, the very cosmos as it was born, in Planckian infancy, a foam of non-commuting coordinates, deconstructing the false smoothness of Riemannian geometry [19] all of this points to a singular truth. The nature's laws do not exist as axioms in a human-based algebra, but rather fugues in a dissonant calculus, and require us to abandon the false comfort of commutativity and the rigid scaffolding of associativity. Then admit a mathematics of relations, adjunctions, and heterarchies, such that logic can become topological, operations are modal, and the very act of measurement is an act of participatory treachery to deeper, unformalizable symmetries.

***Definition:*** Trinition [13], a hypercomplex number with two imaginaries. Let $T$ be the Trinition, $Z \in T$, $Z = a + bi + cj$ where $a, b, c$ are real numbers.

The basis elements of this system are $\{1, i, j\}$, where $i^2 = j^2 = -1$. Here a derived term $k$, a dependent operator that is defined explicitly by the product $ij$. We introduce a parameter $\alpha \in [0,1]$*(note that alpha can be any real finite number, alpha can be piecewise and for generalization, alpha can be a function)*, this parameter will smoothly transition the system between commutative, when the fractional parameter is zero and to a quaternion-like non-commutative when the parameter is 1. In this system, we define the following multiplication rule:

a) $i^2 = j^2 = -1$
b) $i.j = 1 - \alpha + \alpha k$
c) $j.i = 1 - \alpha - \alpha k$ this helps us to introduce non-commutativity when $0 \neq \alpha$

The non-independence of $k$, it is a linear factor that combine the basis element of the Trinition, derived from the product, noting that when the fractional parameter is non null, $k$ can be expressed as follows: (1)

$$k = \frac{ij - (1-\alpha)}{\alpha}$$

When $\alpha = 0$, $ij = ji = 1$, we obtain a fully commutative structure, complex number-like structure.

When $\alpha = 1$, $ij = k$, $ji = -k$, we obtain a fully non-commutative structure, indeed quaternions-like structure.

When $0 < \alpha < 1$ for example $ij$, we obtain a mixed of commutative and non-commutative behaviour. The table below present a resume of the multiplication under some prescribed conditions.

**Table 1**: Multiplication within Trinition

|   | 1 | $i$ | $j$ | $k$ |
|---|---|---|---|---|
| 1 | 1 | $i$ | $j$ | $k$ |
| $i$ | $i$ | $-1$ | $1 - \alpha + \alpha k$ | $-j$ (if $k = i.j$) |
| $j$ | $j$ | $1 - \alpha - \alpha k$ | $-1$ | $i$ (if $k = i.j$) |
| $k$ | $k$ | $-j$ | $i$ | $-1$ (if $k^2 = -1$) |

We recall that, this system is non-associative for each value of the fractional number

$$(ij)i = (1 - \alpha + \alpha k)i \neq (1 - \alpha - \alpha k)i$$

## Properties and interpretations

From the above table, the norm of a Trinition number is given by:

$$\|Z\|_\alpha = \sqrt{x^2 + y^2 + z^2 - 2(1-\alpha)yz} \qquad (2)$$

Thus, the norm for a Trinition number $Z = x + yi + zj$ encodes a geometry parametrized by the fractional order $\alpha$, this will lead of to some important interpretations:

1) The role of the fractional order $\alpha$: Tuning coupling between components.

   - When $\alpha = 1$, the norm simplified to the well-known Euclidean norm $\|Z\| = \sqrt{x^2 + y^2 + z^2}$. This corresponds to a decoupled state where the y-axis and z-axis are orthogonal, akin to quaternions. Indeed, the system behave like a classical 3D Euclidean space.
   - When $\alpha = 0$ the norm becomes $\sqrt{x^2 + y^2 + z^2 - 2yz}$. Here, the term $-2yz$ introduces an anti-correlation between $y$ and $z$. Geometrically, one can notice that, this resembles a space with a non-Euclidian metric, where $y$ and $z$ are coupled.
   - Finally, when $0 < \alpha < 1$ for example, the coupling between $y$ and $z$ interpolates between full anti-correlation when the fractional order is 0 and the orthogonality when it is 1. The coefficient $-2(1-\alpha)$ serves as a tunable interaction strength between $y$ and $z$

For the geometric interpretation, the norm defines a metric tensor with off-diagonal terms:
$$\|Z\|_\alpha^2 = x^2 + y^2 + z^2 - 2(1-\alpha)yz$$

The component $x^2$ is independent, while $y^2 + z^2 - 2(1-\alpha)yz$ is coupled. Here also, the metric tensor $G$ in the $(y, z)-$subspace is:

$$G = \begin{bmatrix} 1 & -(1-\alpha) \\ -(1-\alpha) & 1 \end{bmatrix}$$

The metric obtain here describes a warped plane where $y$ and $z$ are not orthogonal. The warping is controlled by the fractional order. We will formally define the Trinition space $T_\alpha$, a 3 D parameterized by the fractional parameter $\alpha$, blending algebraic and geometric structures. We present below its axiomatic framework, properties and implications. We will start with inner product

### The inner product

$$\langle Z_1, Z_2 \rangle_\alpha = x_1 x_2 + y_1 y_2 + z_1 z_2 - (1-\alpha)(y_1 z_2 + y_2 z_1) \tag{3}$$

Symmetry $\langle Z_1, Z_2 \rangle_\alpha = \langle Z_2, Z_1 \rangle_\alpha$

Linearity: $\langle a Z_1 + b Z_2, Z_3 \rangle_\alpha = a \langle Z_1, Z_3 \rangle_\alpha + b \langle Z_2, Z_3 \rangle_\alpha$

Positivity definiteness $\langle Z, Z \rangle_\alpha \geq 0$, with equality if and only if $Z = 0$ for $\alpha > 0$

We define a distance:

$$d_\alpha(Z_1, Z_2) = \|Z_1 - Z_2\|_\alpha = \sqrt{(\Delta x)^2 + (\Delta y)^2 + (\Delta z)^2 - 2(1-\alpha)\Delta y \Delta z} \tag{4}$$

What is then our geometry structure? Let us start with the metric tensor,

$$G = \begin{bmatrix} 1 & 0 & 0 \\ 0 & 1 & -(1-\alpha) \\ 0 & -(1-\alpha) & 1 \end{bmatrix}$$

For curvature, the space is flat for $\alpha = 1$ this correspond to the Euclidean, since

$$d_1(Z_1, Z_2) = \|Z_1 - Z_2\|_1 = \sqrt{(\Delta x)^2 + (\Delta y)^2 + (\Delta z)^2}. \tag{5}$$

But when $0 < \alpha < 1$, for example, we have curved leading a different geometry. This is indeed a tunable interaction because the distances in the y-z plane are compressed or stretched accordingly to the parameter. For example if we consider the following Trinition points $A(0,1,0)$ and $B(0,0,1)$ the distance between these two points varies accordingly to the value of alpha. The simulation is presented in Figure 1 as a function of alpha. (6)

$$d_\alpha(A - B) = \sqrt{0^2 + 1^2 + 1^2 - 2(1-\alpha)(1)(1)} = \sqrt{2\alpha}$$

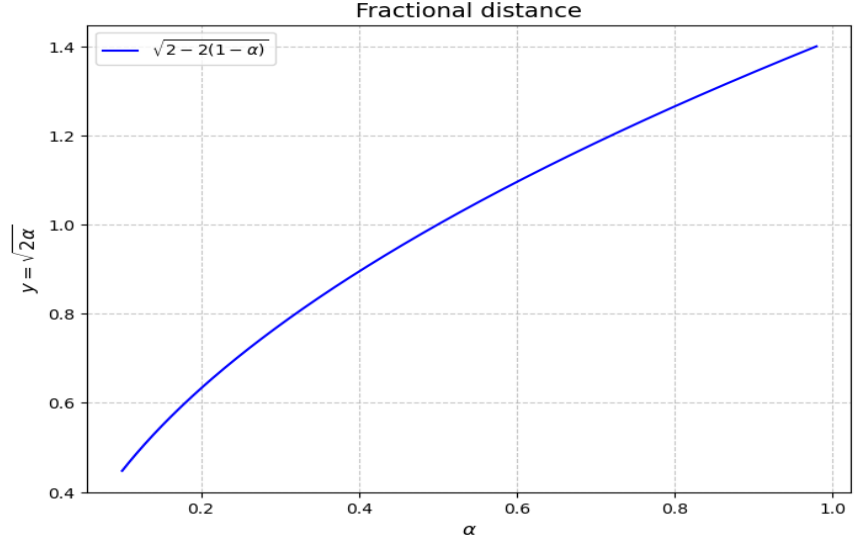

**Figure 1**: An $\alpha-$distance

$\alpha = 0$, this degenerates into a semi-metric space. (7)

$$d_0(Z_1, Z_2) = \sqrt{(\Delta x)^2 + (\Delta y)^2 + (\Delta z)^2 - 2\Delta y \Delta z} = \sqrt{(\Delta x)^2 + (\Delta y - \Delta z)^2}$$

We not only get a metric via the Trinition distance, but also a whole family of geometries with dimensions interacting in a controlled way via $\alpha$. It has ramifications in mathematics (deformation theory, hyper-complex algebras), physics (toy models of curved space-time), machine learning (adaptive embeddings). Though this does not go to displacing Euclidean space and quaternions, it adds a very useful toolset for building a system where correlation (or curvature) can be adjusted as a parameter, which can be viewed as a contribution to both the theoretical and applied mathematics. We shall now verify the triangular inequality within the framework of Trinition.

**Theorem:** Triangle inequality for Trinition space

For vectors $u = (u_1, u_2, u_3)$ and $v = (v_1, v_2, v_3)$ in the Trinition space $T_\alpha$, then

$$\|u + v\|_\alpha \leq \|u\|_\alpha + \|v\|_\alpha \qquad (8)$$

**Proof**:

$$\|u + v\|_\alpha^2 = (u_1 + v_1)^2 + (u_2 + v_2)^2 + (u_3 + v_3)^2 - 2(1-\alpha)(u_2 + v_2)(u_3 + v_3) \quad (9)$$

We can expand the above expression and grouping terms, we obtain: (10)

$$\|u+v\|_\alpha^2 = \{u_1^2 + u_2^2 + u_3^2 - 2(1-\alpha)u_2 u_3\} + \{v_1^2 + v_2^2 + v_3^2 - 2(1-\alpha)v_2 v_3\}$$
$$+ 2\{u_1 v_1 + u_2 v_2 + u_3 v_3 - (1-\alpha)(u_2 v_3 + u_3 v_2)\}$$

Therefore, we have: (11)

$$\|u+v\|_\alpha^2 = \|u\|_\alpha^2 + \|v\|_\alpha^2 + 2\langle u, v\rangle_\alpha$$

We shall now apply the Cauchy-Schwarz inequality on the inner product. We now hope to establish that (12)

$$\langle u, v\rangle_\alpha \le \|u\|_\alpha \|v\|_\alpha$$

We will present proof for the Cauchy-Schwarz. First we will make use of the defined G metric tensor. We notice that, the inner product (13)

$$\langle u, v\rangle_\alpha = u^T \begin{bmatrix} 1 & 0 & 0 \\ 0 & 1 & -(1-\alpha) \\ 0 & -(1-\alpha) & 1 \end{bmatrix} v$$

Indeed the fractional parameter $\alpha \in [0, \alpha]$, $G$ is positive definite (in fact all eigenvalues are positive, $\lambda = 1, \lambda = 1 \mp (1-\alpha) > 0$. At $\alpha = 0$, $G$ is positive semi-definite. We shall now apply the Cauchy-Schwarz inequality [20]. We let, $t \in \mathbb{R}$, the squared norm of $u + tv$ is non-negative (14)

$$\|u+v\|_\alpha^2 = \langle u+tv, u+tv\rangle_\alpha \ge 0$$

Expanding the above: (15)

$$\langle u+tv, u+tv\rangle_\alpha = \langle u, u\rangle_\alpha + 2t\langle v, u\rangle_\alpha + t^2\langle v, v\rangle_\alpha = t^2\|v\|_\alpha^2 + 2t\langle v, u\rangle_\alpha + \|u\|_\alpha^2 \ge 0$$

We have then obtained a quadratic inequality in $t$. We recall for the quadratic $At^2 + Bt + C \ge 0$ to hold for all $t$, its discriminant must be negative or null. Now applying this on our equation, we obtain: (16)

$$4\langle v, u\rangle_\alpha^2 - 4\|v\|_\alpha^2 \|u\|_\alpha^2 \le 0$$

Therefore, divide by 4 leads to (17)

$$\langle v, u \rangle_\alpha^2 \leq \|v\|_\alpha^2 \|u\|_\alpha^2$$

Therefore, we obtain (18)

$$\|u + v\|_\alpha \leq \|u\|_\alpha + \|v\|_\alpha$$

We have an edge case, this when $\alpha = 0$, then

$$G = \begin{bmatrix} 1 & 0 & 0 \\ 0 & 1 & -1 \\ 0 & -1 & 1 \end{bmatrix}$$

With eigenvalues 1, 2, 0. The null eigenvector (0,1,1) satisfies $\|u\|_0 = 0$, but for non-degenerate vectors $(u, v \neq 0)$, the inequality still holds. Indeed, the triangular and the Cauchy-Schwarz inequality hold for $\alpha \in [0,1]$, this indeed insure that, the defined Trinition space is a valid inner product space with several application in physics, optimization and machine learning. A direct comparison with the Euclidean space leads to the following table

Table 2: Direct comparison with Euclidean space

| Property | Euclidean Space $\alpha = 1$ | Trinition Space ($\alpha \in [0,1)$) |
|---|---|---|
| Orthogonality | Basis vectors orthogonal | Basis vectors coupled $(y - z)$ |
| Algebra | commutative | Non-commutative interactions |
| Curvature | Flat $(k = 0)$ | Curved $(k = 4(1 - \alpha)^2)$ |

Trinition space is a high-dimensional manifold that generalizes Euclidean geometry under a very narrow regime of $\alpha \in [0,1]$, leading to a tunable model for anisotropic systems. By following this interval (without being overly restrictive), all mathematical properties (positive definiteness, curvature, etc.) are indeed rigorously valid, allowing to be used not only in the pure galaxy modelling, but also in physics, machine learning.

## Multiplicativity with the framework of Trinition

We next investigate the multiplicativity under the framework of the Trinition. For all $Y, Z \in T_\alpha$ do we have in general $\|YZ\|_\alpha = \|Y\|_\alpha \|Z\|_\alpha, \forall \alpha$ finite? .

We recall that,

$$\|Z\|_\alpha = \sqrt{x^2 + y^2 + z^2 - 2(1-\alpha)yz}, \quad \|Y\|_\alpha = \sqrt{u^2 + v^2 + w^2 - 2(1-\alpha)wv}$$

$\alpha = 1$: We have a quaternion-like algebra, since the multiplication rule $ij = k$, where $k$ is a new basis element similar to the quaternions and the norm is $\|Z\|_\alpha = \sqrt{x^2 + y^2 + z^2}$. We will require that, $jk, ki$ to mirror quaternions: $i^2 = j^2 = k^2 = -1$, $ij = k, jk = i, ki = j$. Therefore, having the closure, the norm is multiplicative like in the case of quaternions and then, $\|YZ\|_1 = \|Y\|_1\|Z\|_1$, But when the equality does not hold since in the case $\alpha \neq 1$, we have that

$$YZ = xu - yv - zw + (1-\alpha)(yw + zv) + (xv + yu)i + (xw + zu)j + \alpha(yw - zv)k \quad (19)$$

Which is a pure quaternion number providing that $\alpha \neq 0$. Thus by using the quaternion norm, we obtain

$$\|ZY\|_\alpha^Q = \left\{(xu - yv - zw + (1-\alpha)(yw + zv))^2 + (xv + yu)^2 + \left(\alpha(yw + zv)\right)^2\right\}^{\frac{1}{2}} \quad (20)$$

While

$$\|Z\|_\alpha \|Y\|_\alpha = \sqrt{x^2 + y^2 + z^2 - 2(1-\alpha)yz}\sqrt{u^2 + v^2 + w^2 - 2(1-\alpha)wv} \quad (21)$$

### Trinition Gram-Schmidt [21]

1. Let $(V_1, V_2, V_3) \in \Omega^3$ an initial set of vectors, we want an orthonormal set $\{u_1, u_2, u_3\}$ under the Trinition norm. We suggest this procedure;

    I. Normalize $V_1$

    $$u_1 = \frac{V_1}{\|V_1\|_\alpha}, \quad \|V_1\|_\alpha^2 = \langle V_1, V_1\rangle_\alpha \quad (21)$$

    II. To make $V_2$ orthorgonal to $u_1$

    $$\widetilde{V_2} = V_2 - \frac{\langle V_2, u_1\rangle_\alpha}{\langle u_1, u_1\rangle}u_1 \quad (22)$$

Then the normalize is,

$$u_2 = \frac{\widetilde{V_2}}{\|\widetilde{V_2}\|_\alpha} \quad (23)$$

    III. Orthogonalize $V_3$

$$\widetilde{V_3} = V_3 - \frac{\langle V_3, u_1 \rangle_\alpha}{\langle u_1, u_1 \rangle_\alpha} u_1 - \frac{\langle V_3, u_2 \rangle_\alpha}{\langle u_2, u_2 \rangle_\alpha} u_2 \tag{24}$$

Then

$$u_3 = \frac{\widetilde{V_3}}{\|\widetilde{V_3}\|_\alpha} \tag{25}$$

That is the traditional Gram-Schmidt, except every 'dot product' and norm is replaced by the defined Trinition product and norm. Zero divisors or partial commutation can make the overall procedure troublesome or yield strange edge cases. In this case, one will use the more symmetric method, simulator to Householder transformations from classical linear algebra, which is usually numerically stable systematically when it comes to cross terms.

We define the Trinition Householder reflection. Here, we want an operator $H_\alpha^{T\alpha}$ that reflects a vector about a plane orthogonal in Trinition framework to some vector $V$. We recall that in classical geometry, that reflection is;

$$H(\bar{V}) = \bar{V} - 2 \frac{\langle \bar{V}, V \rangle}{\langle V, V \rangle} V \tag{26}$$

Therefore, in the Trinition geometry, we define;

$$H_\alpha^{T\alpha}(\bar{V}) = \bar{V} - 2 \frac{\langle \bar{V}, V \rangle_\alpha}{\langle V, V \rangle_\alpha} V \tag{27}$$

Indeed, we will need $\langle V, V \rangle_\alpha > 0$ so that the fraction is well defined.

# Geometrical figures in Trinition

Below, we present a step-by-step outline of how to derive the 3D geometric equations for five well known shapes including sphere, cylinder, cone, ellipsoid and hyperboloid within the framework of Trinition.

We recall that;

$$\|(x, y, z)\|_\alpha^2 = x^2 + y^2 + z^2 - 2(1 - \alpha)yz \tag{28}$$

## Trinition sphere

Therefore, the equation of a sphere within the framework Trinition is $\|(x, y, z)\|_\alpha^2 = R^2$. Note that when $\alpha = 1$, we will recover the classical sphere $x^2 + y^2 + z^2 = R^2$. When $\alpha = 0$, we obtain $x^2 + y^2 + z^2 - xyz = R^2$ which is a more degenerate shape when $\alpha \in (0,1)$, we have a deformed sphere. We present below in Figure 2 an $\alpha-$sphere, degenerated and deformed

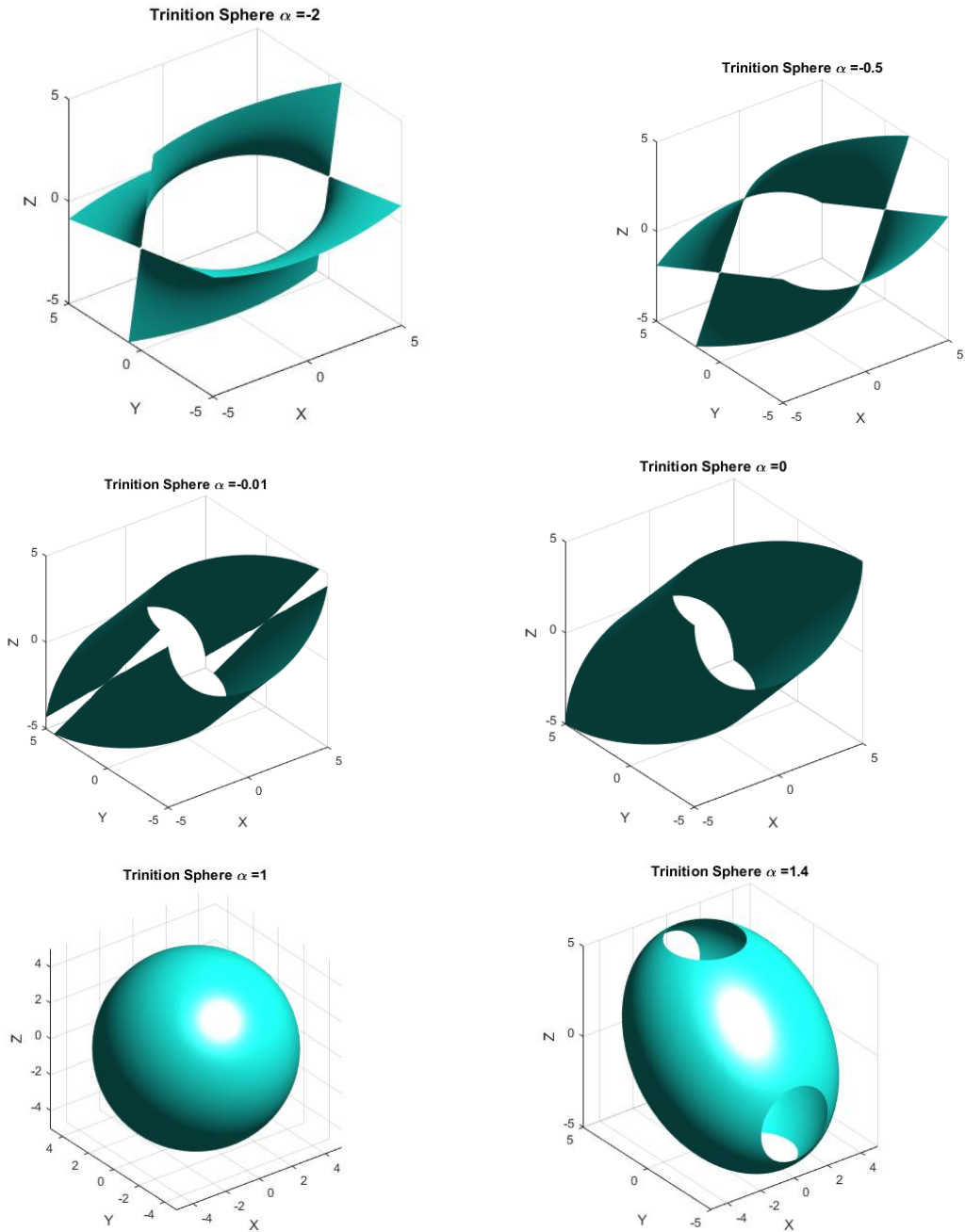

**Figure 2**: An $\alpha$ −sphere deformed and degenerated as a function of alpha for a fixed radius

Trinition cylinder aligned with x-axis

$$\sqrt{y^2 + z^2 - 2(1-\alpha)yz} = R \qquad (29)$$

Which is equivalent to;

$$y^2 + z^2 - 2(1-\alpha)yz = \mathbb{R}^2, x \in \mathbb{R} \qquad (30)$$

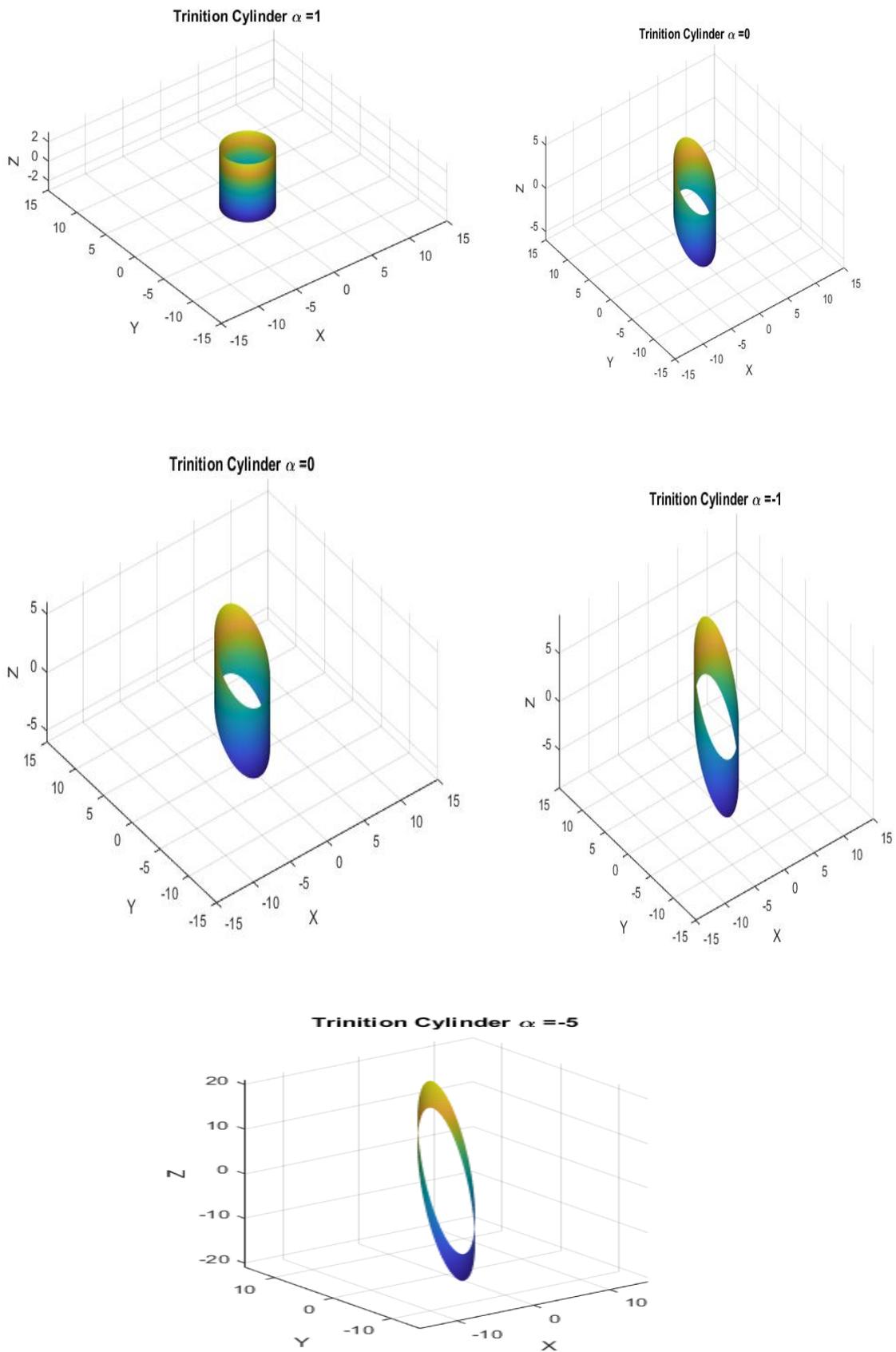

**Figure 3**: An-α −cylinder generated, and deformed for different alphas

## Trinition cone (vertex at the origin axis along the x-axis).

Here, the distance from x-axis in the Trinition sense is;

$$\sqrt{y^2 + z^2 - 2(1-\alpha)yz} \qquad (31)$$

A right circular cone condition says;

$$x^2 = \bar{k}^2[y^2 + z^2 - 2(1-\alpha)yz] \qquad (32)$$

Therefore, the Trinition cone equation is given as;

$$x^2 = \bar{k}^2[y^2 + z^2 - 2(1-\alpha)yz] = 0 \qquad (33)$$

## Trinition ellipsoid

To have the cross-section in $(y, z)$ to follow the Trinition circle structure, we can define;

$$\frac{x^2}{a} + \frac{y^2 + z^2 - 2(1-\alpha)yz}{b^2} = 1 \qquad (34)$$

With $a \neq b$. We present in Figure 4 below a deformed, and degenerated ellipsoid within the framework of Trinition:

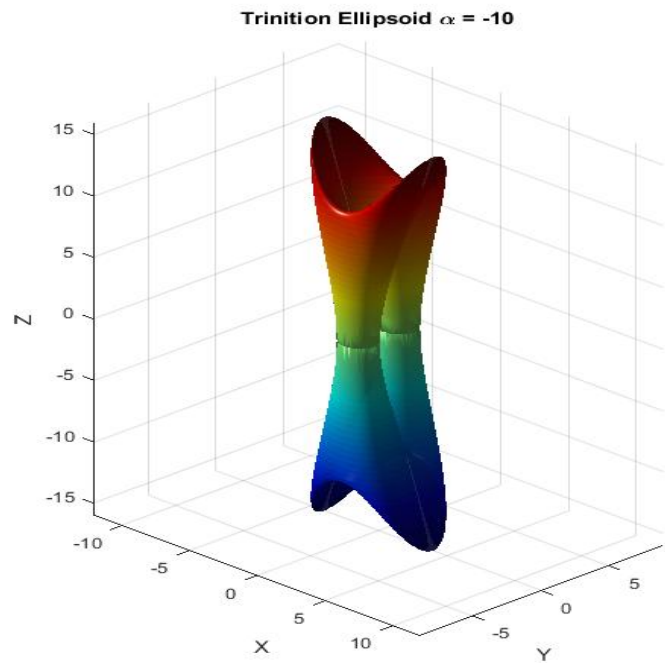

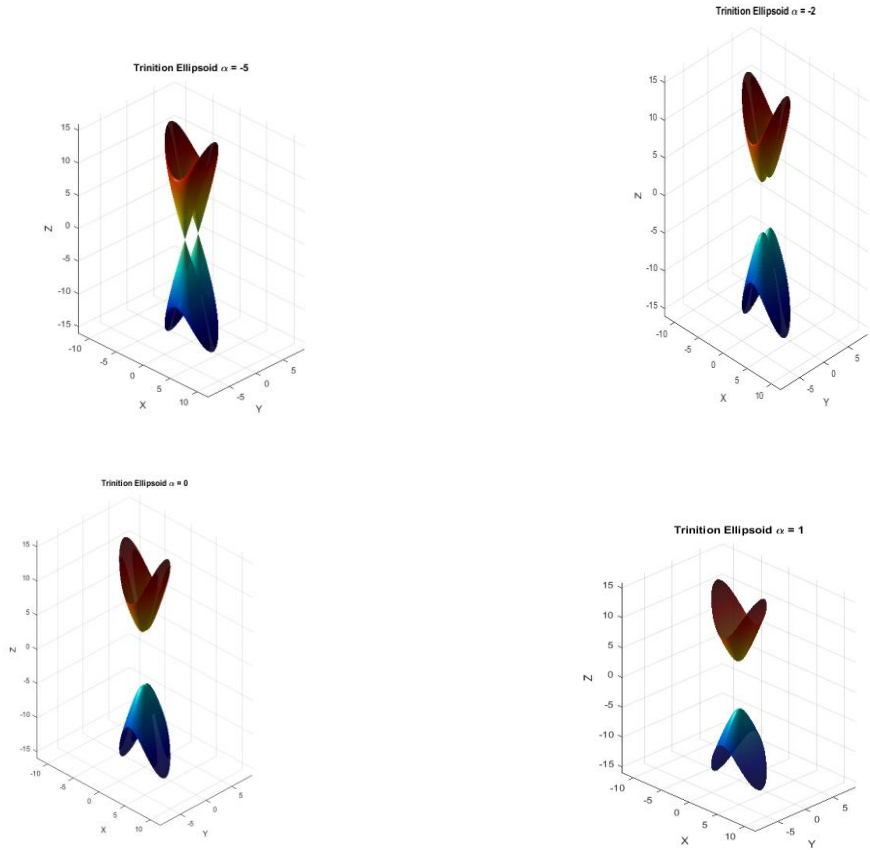

**Figure 4**: An $\alpha-$Ellipsoid (3D and 2D) degenerated and deformed for different values of alpha.

### Trinition hyperboloid

If we want a shape that looks like a hyperboloid by with $(y,z)$ we will define;

$$x^2 + (y^2 + z^2 - 2(1-\alpha)yz) - \lambda w^2 = 1 \qquad (35)$$

This will depend on how many coordinates we use. A simpler version could be;

$$x^2 - y^2 + z^2 - 2(1-\alpha)yz = R^2 \qquad (36)$$

Besides these known shapes, we will introduce new shapes.

### Locus with Trinition distances

We define the following shape.

All points $(x, y, z)$ whose sum of Trinition distances two lines is constant. Trinition distance to a line say the x-axis is;

$$\sqrt{y^2 + z^2 - 2(1-\alpha)yz} \qquad (37)$$

Therefore, one could pick two lines $l_1$ and $l_2$ in 3D and define;

$$\sqrt{\Omega_1(y,z)} + \sqrt{\Omega_2(y,z)} = constant \tag{38}$$

Where $\Omega_i(y,z) = y^2 + z^2 - 2(1-\alpha)yz$.

The resulting shape is clearly new in standard geometry. We can define a shape by mixing Trinition distance in one direction and classical distance in another.

*example,*

$$|x| + \sqrt{x^2 + z^2 - 2(1-\alpha)xy} = constant \tag{39}$$

The above is suggested because of the plane is $x = 0$, the Euclidean distance is $|x|$. If the axis is the y-axis, the Trinition distance is $(x, z)$, plane might be $\sqrt{x^2 + z^2 - 2(1-\alpha)xz}$.

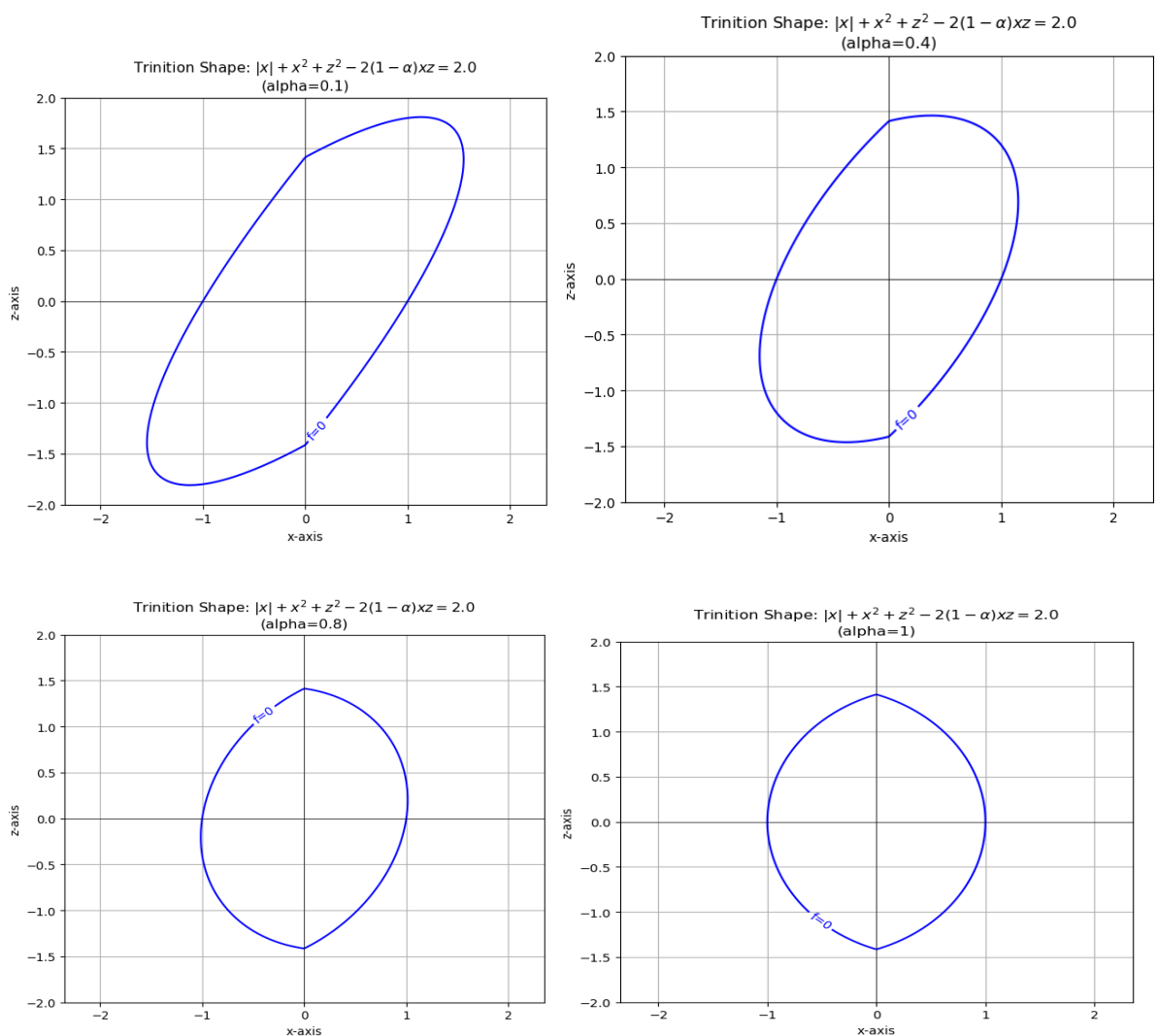

**Figure 5**: Deformable non-conventional circle in in x-z plane

## Trinition reverse of sweeps

In Trinition geometry, the distance to the z-axis is $\sqrt{x^2 + y^2 - 2(1-\alpha)xy}$. If one revolves a 2D curve around the z-axis using the Trinition radius,

$$r_\alpha = \sqrt{x^2 + y^2 - 2(1-\alpha)xy} \tag{40}$$

We can get a new shape that will differ from the standard revolution.

For example, if we revolve the line $r_\alpha = a + bz$, that defines as;

$$\sqrt{x^2 + y^2 - 2(1-\alpha)xy} = a + bz \tag{41}$$

We will obtain a shape that do not exist in standard Euclidean geometry but it is constant with Trinition.

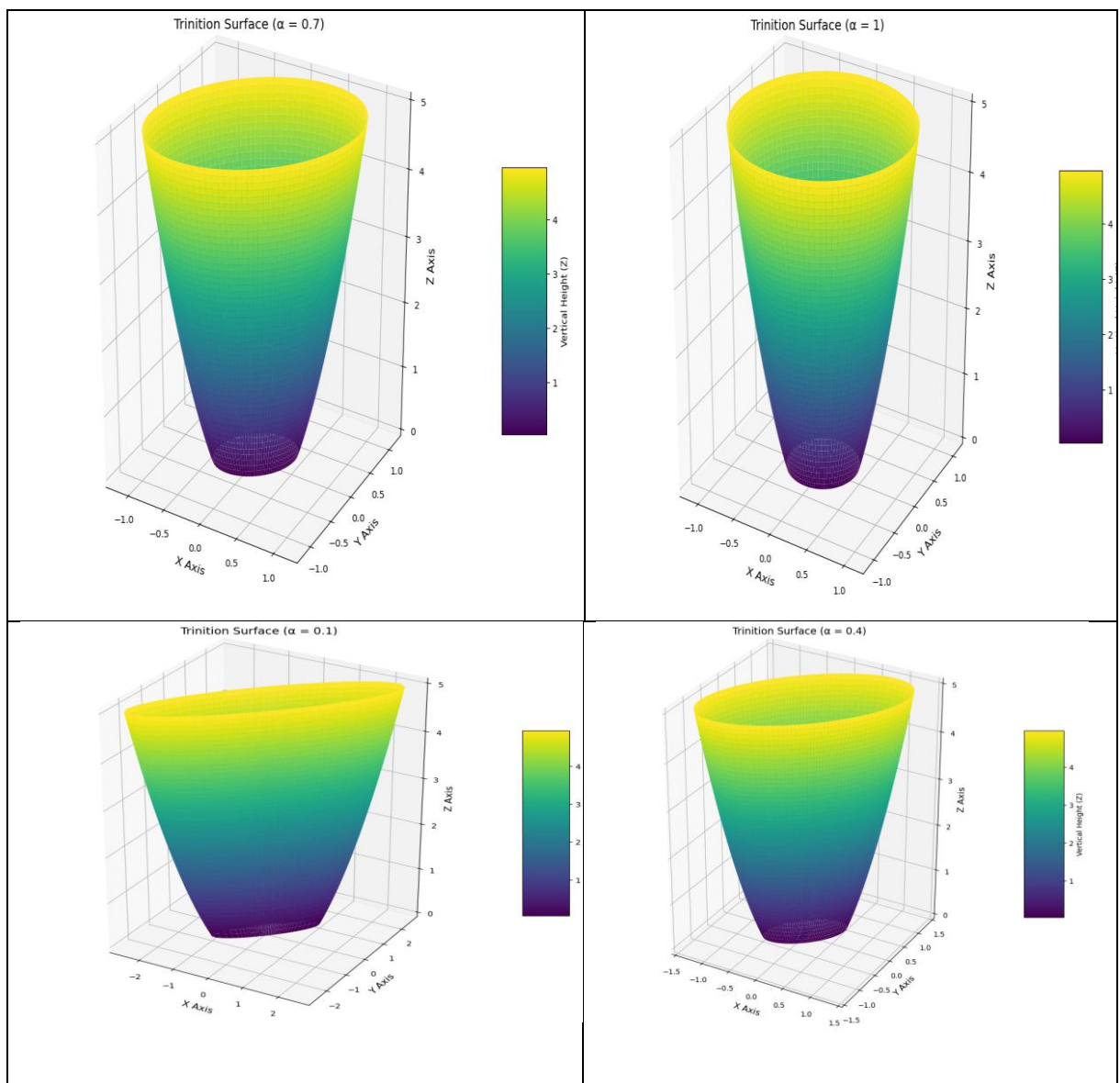

**Figure 6**: Trinition reverse of sweeps for different alphas

## Sphere with deformable radius

We can have a sphere with deformable radius.

$$x^2 + y^2 + z^2 - 2(1-\alpha)xy = a^2 - b^2 + c^2 - 2(1-\alpha)bc \tag{42}$$

$$\|(x,y,z)\|_\alpha^2 = R^2 \tag{43}$$

Circle, cylinder and other geometric figures with radii can be defined similarly.

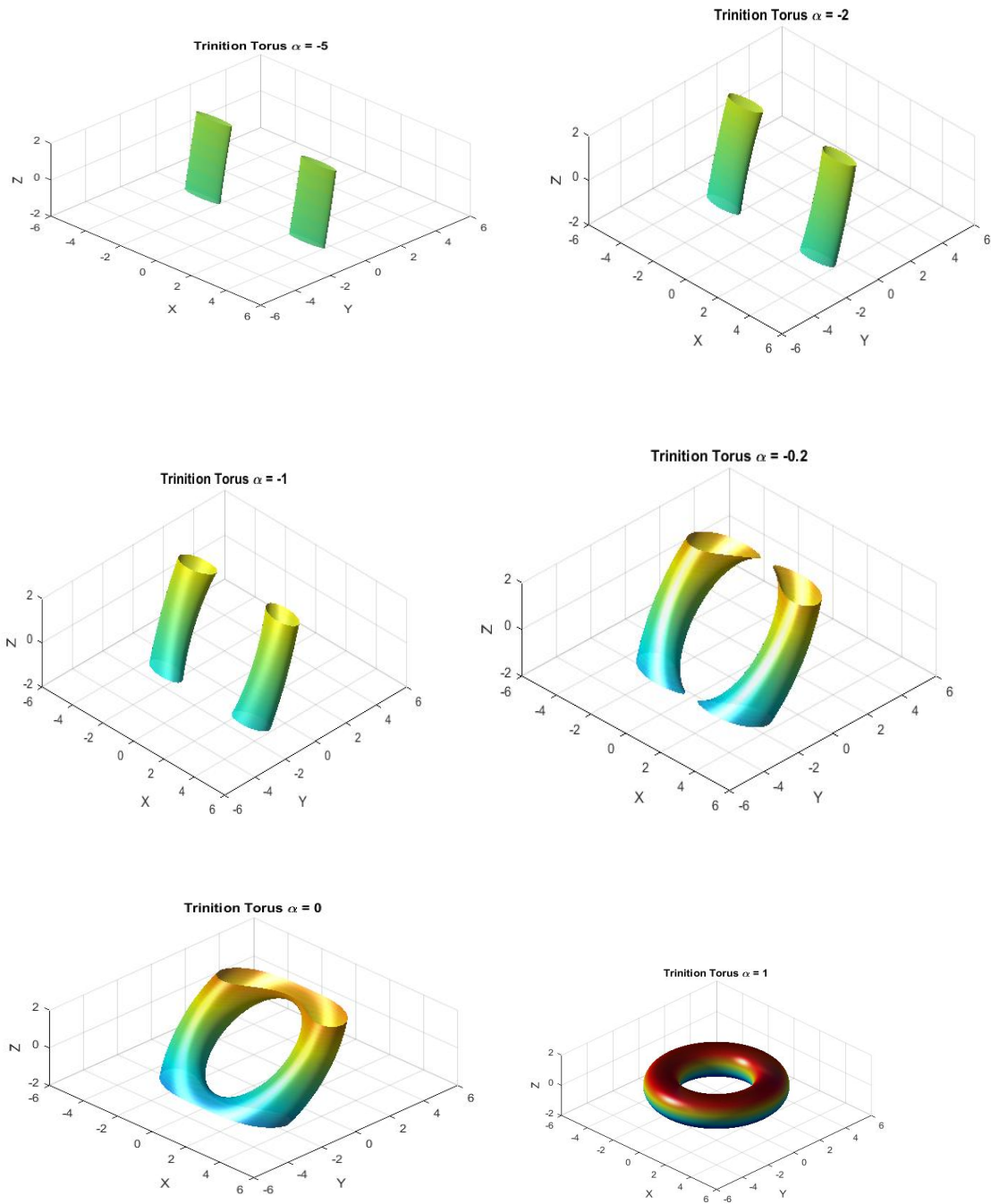

**Figure 7**: An $\alpha-$Trinition Torus deformed and generated for different values of alphas.

# Novel inequality in Trinition space

We use the unique structure of Trinition space to obtain new inequalities that will perhaps generalize or deviate from classical results.

## Parallelogram law deviation.

We recall that in Euclidean space, the parallelogram law holds as follows;

$$\|u + V\|^2 + \|u - V\|^2 = 2\|u\|^2 + 2\|V\|^2 \tag{44}$$

But, in our Trinition space, this law will be altered because of the y-z coupling.

$$\|u + V\|_\alpha^2 + \|u - V\|_\alpha^2 = 2\|u\|_\alpha^2 + 2\|V\|_\alpha^2 + 4(1 - \alpha)(u_2 V_3 + u_3 V_2) \tag{45}$$

The above leads to the following inequality;

$$\|u + V\|_\alpha^2 + \|u - V\|_\alpha^2 = 2\|u\|_\alpha^2 - 2\|V\|_\alpha^2 \leq 4(1 - \alpha)\|u\|_\alpha \|V\|_\alpha \tag{46}$$

Indeed, this qualifies the non-Euclidean of the Trinition space via $\alpha$.

## Eigenvalue-Bounded inner product

Theorem: If $\lambda_{min} = \alpha$ and $\lambda_{max} = 2 - \alpha$ are the extreme eigenvalues of the metric tensor $G$, then $\forall U, V \in T_\alpha$.

$$\lambda_{min} \|U\|_2 \|V\|_2 \leq \langle U, V \rangle_\alpha \leq \lambda_{max} \|U\|_2 \|V\|_2 \tag{47}$$

Where $\|\ \|_2$ is the Euclidean norm.

To prove this, we recall that in matrix form we defined the inner product $\langle x, y \rangle_\alpha = x^T G y$. Since $G$ is real and symmetric, and positive definite, its eigenvalues lines in $[\lambda_{min}, \lambda_{max}]$.

By standard linear algebra, $\forall x, y$:

$$\lambda_{min} x^T y \leq x^T G y \leq \lambda_{max} x^T y \tag{48}$$

Replace, $x = \frac{x}{\|U\|_2}$ and $y = \frac{U}{\|V\|_2}$, then, we obtain the requested inequality. The above result implies that $\langle ., . \rangle_\alpha$ is sandwhiched between $\alpha$ and $2 - \alpha$ times the well-known Euclidean product. This bounds Trinition geometry by classical Euclidean geometry.

## The cross-term AM-GM bound.

We notice that the cross term $-2(1 - \alpha)yz$ in the Trinition norm is bounded by the AM-GM, that is to say;

$$|yz| \leq \frac{y^2 + z^2}{2} \iff 2|yz| \leq y^2 + z^2, \tag{49}$$

$$-(1 - \alpha)(y^2 + z^2) \leq -2(1 - \alpha)|yz|$$

Noting that AM (arithmetic mean) and GM (geometric mean).

$$\|W\|_\alpha = x^2 + y^2 + z^2 - 2(1 - \alpha)yz \tag{50}$$

The negative cross term cannot reduce the sum by more than $(1 - \alpha)\frac{y^2+z^2}{2}$. Thus, we will get a minimal value $\alpha(y^2 + z^2)$ in 2D slices or $\alpha w_2^2 + \alpha w_3^2$ in 3D, plus the term $w_1^2$ such that;

$$\|W\|_\alpha^2 \geq w_1^2 + \alpha w_2^2 + \alpha w_3^2 \tag{51}$$

This inequality informs us that the norm in the Trinition is never drastically smaller than a scaled Euclidean norm, given a lower bound in terms of the classical geometry.

## Hölder-type inequality

$\forall q, p > 1$ with $\frac{1}{p} + \frac{1}{q} = 1$, we have the Trinition Hölder

$$\langle U, V \rangle_\alpha \leq \|U\|_\alpha^p \|V\|_\alpha^q \tag{52}$$

Where, we define,

$$\|w\|_\alpha^p = \left(w_1^p + w_2^p + w_3^p - 2(1 - \alpha)|w_2 w_3|^{\frac{p}{2}}\right)^{\frac{1}{p}} \tag{53}$$

Proof: In the Trinition, the inner product,

$$\langle u, V \rangle_\alpha = u_1 V_1 + u_2 V_2 + u_3 V_3 - (1 - \alpha)(u_2 V_3 + u_3 V_2) \tag{54}$$

Where $u = (u_1, u_2, u_3), V = (V_1, V_2, V_3)$

When $\alpha = 1$, this reduces to classical dot product $u.V$.

But when $\alpha \in [0,1)$ we have an additional cross-term $-(1 - \alpha)(u_2 V_3 + u_3 V_2)$.

The defined $L^P$ −type norm $\|w\|_\alpha^p$.

$$\|w\|_\alpha^p = \left(w_1^p + w_2^p + w_3^p - 2(1 - \alpha)|w_2 w_3|^{\frac{p}{2}}\right)^{\frac{1}{p}} \tag{55}$$

Now if $\alpha = 1$, then $\|w\|_1^p = \left(w_1^p + w_2^p + w_3^p\right)^{\frac{1}{p}}$.

Otherwise the cross-term $-2(1 - \alpha)|w_2 w_3|^{\frac{p}{2}}$ modifies the measure of $(w_2, w_3)$. Now, we have to sow that $|\langle u, V \rangle_\alpha| \leq \|u\|_\alpha^p \|V\|_\alpha^p$ find of all when $p = 2$, we have the Cauchy-Schwarz in Trinition geometry.

$$\langle u, V \rangle_\alpha = \left( u_1 V_1 + u_2 V_2 + u_3 V_3 - (1-\alpha)(u_2 V_3 + u_3 V_2) \right) \tag{56}$$

Taking the absolute value, we obtain,

$$|\langle u, V \rangle_\alpha| \leq |u_1 V_1| + |u_2 V_2| + |u_3 V_3| + (1-\alpha)|u_2 V_3 + u_3 V_2| \tag{57}$$

We have to show that the right-hand side is less than $\|u\|_\alpha^p \|V\|_\alpha^p$.

We can use the Hölder or Young's inequality,

$$|u_1 V_1| \leq \frac{|u_1|^p}{p} + \frac{|V_1|^q}{q}, |u_2 V_2| \leq \frac{|u_2|^p}{p} + \frac{|V_2|^q}{q} \text{ and } |u_3 V_3| \leq \frac{|u_3|^p}{p} + \frac{|V_3|^q}{q}$$

Therefore,

$$|u_1 V_1| + |u_2 V_2| + |u_3 V_3| \leq \sum_{j=1}^{3} \frac{|u_j|^p}{p} + \sum_{j=1}^{3} \frac{|V_j|^q}{q} \tag{58}$$

For the cross-term, we first use the triangular inequality and then the Hölder to obtain;

$$|u_2 V_3 + u_3 V_2| < |u_2 V_3| + |u_3 V_2| \tag{59}$$

$$|u_2 V_3 + u_3 V_2| \leq \frac{|u_2|^p + |u_3|^p}{p} + \frac{|V_2|^q}{q} + \frac{|V_3|^q}{q} \tag{60}$$

By multiplying by $2(1-\alpha)$

$$2(1-\alpha)|u_2 V_3 + u_3 V_2| \leq \frac{2(1-\alpha)}{p}(|u_2|^p + |u_3|^p) + \frac{2(1-\alpha)}{q}(|V_2|^q + |V_3|^q) \tag{61}$$

Thus,

$$\frac{1}{p}\sum_{j=1}^{3}|u_j|^p + \frac{2(1-\alpha)}{q}(|u_2|^q + |u_3|^q) = \frac{1}{p}|u_1|^p + (1 + 2(1-\alpha))|u_2|^p \tag{62}$$

$$+ (1 + 2(1-\alpha)|u_3|^p)$$

But we must note that $(1 + 2(1-\alpha)) = 3 - 2\alpha$.

However, in $\|u\|_\alpha^p$, the cross-term with negative sign $-2(1-\alpha)|u_2 u_3|^{\frac{p}{2}}$ is not accounted for explicitly in the sum. Here, the Riesz-Thorin [22] or Minkowski argument will be used, bounding $|u_2 u_3|^{\frac{p}{2}}$ by the partial AM-GM, we finally get a factor that merges into the $\|u\|_\alpha^p \|V\|_\alpha^p$ product.

In other words, these bounding lines produce a function of,

$$\sum_{j=1}^{3}|u_j|^p - 2(1-\alpha)|u_2 u_3|^{\frac{p}{2}}$$

And

$$\sum_{j=1}^{3}|V_j|^q - 2(1-\alpha)|V_2 V_3|^{\frac{q}{2}}$$

$$\langle u,V\rangle_\alpha \leq \left(\sum_{j=1}^{3}|u_j|^p - 2(1-\alpha)|u_2 u_3|^{\frac{p}{2}}\right)^{\frac{1}{p}} \left(\sum_{j=1}^{3}|V_j|^q - 2(1-\alpha)|V_2 V_3|^{\frac{q}{2}}\right)^{\frac{1}{q}} \quad (63)$$

$$= \|u\|_\alpha^p \|V\|_\alpha^p$$

Note that by defining $\Phi_\alpha(w) = w_1^p + w_2^p + w_3^p - 2(1-\alpha)|w_2 w_3|^{\frac{p}{2}}$ and show that we have convexity when $p > 1$. Then we define the dual functional $\Phi_\alpha^\forall$ that will act like,

$$\Phi_\alpha(w) + \Phi_\alpha^*(V) \geq w_1 V_1 + w_2 V_2 + w_3 V_3 - (1-\alpha)(w_2 V_3 + w_3 V_2) \quad (64)$$

When $p = q = 2$, we get the Trinition Cauchy-Schwarz inequality

$$\langle u,V\rangle_\alpha \leq \|u\|_\alpha^2 \|V\|_\alpha^2 \quad (65)$$

In general, the above inequality is the mixed Minkowski or Riesz-Thorin approach.

The question one will ask is: How useful this equation is? The answer is:

It will be crucial for the analysis of PDEs or functional spaces under Trinition geometry. The inequality will ensure that standard expansions, partial integrals, or discretised summations stay controlled by Trinition norm. This will help us to adapt classical theorems for example, boundness of operators, existence and uniqueness arguments to the Trinition framework.

*Dynamic $\alpha$ − contraction inequality*

Let us assume that $\alpha$ is a function of $t$, but $\alpha = \alpha(t)$, then

$$\frac{d\|w(t)\|_{\alpha(t)}}{dt} \leq \|\dot{w}(t)\|_{\alpha(t)} + |\dot{\alpha}(t)|\frac{|w_2(t)w_3(t)|}{\|w(t)\|_{\alpha(t)}} \quad (66)$$

Note that the proof is via the chain rule.

*Reverse triangle inequality*

$$\forall \alpha, \ \big|\|u\|_\alpha - \|V\|_\alpha\big| \geq \|u - V\|_\alpha - \frac{2(1-\alpha)|u_2 V_3 + u_3 V_2|}{\|u\|_\alpha + \|V\|_\alpha} \quad (67)$$

## Spectral Gap inequality

For $w = (0, y, z) \in \mathbb{R}^3$

$$\|w\|_\alpha = y^2 + z^2 - 2(1-\alpha)yz \geq \alpha(y^2 + z^2) \tag{68}$$

Hence,

$$\|w\|_\alpha \geq \sqrt{\alpha}\sqrt{y^2 + z^2} = \sqrt{\alpha}\|w\|_2 \tag{69}$$

Proof:

$$-2(1-\alpha)yz > -2(1-\alpha)\frac{y^2 + z^2}{2} = -(1-\alpha)(y^2 + z^2) \tag{70}$$

Then $y^2 + z^2 - 2(1-\alpha)yz \geq \alpha(y^2 + z^2)$

This leads to;

$$\|w\|_\alpha \geq \sqrt{\alpha}\|w\|_2 \tag{71}$$

## Trinition Minkowski inequality

Let $\alpha \in [0,1]$ and vectors $(u, V) \in \mathbb{R}^3$.

$$\|w\|_\alpha^2 = w_1^2 + w_2^2 + w_3^2 - 2(1-\alpha)w_2 w_3 \tag{72}$$

There exists a constant $f(\alpha) > 0$ such that;

$$\|u + V\|_\alpha \leq \|u\|_\alpha + \|V\|_\alpha + f(\alpha)|u_2 V_3 + u_3 V_2| \tag{73}$$

The proof is derived by expanding;

$$\|u + V\|_\alpha^2 = (u_1 + V_1)^2 + (u_2 + V_2)^2 + (u_3 + V_3)^2 \tag{74}$$
$$- 2(1-\alpha)(u_2 + V_2)(u_3 + V_3)$$

Then the cross-terms can be bounded by a function of $\|u\|_\alpha^2$ and $\|V\|_\alpha^2$.

$\alpha \to 0$, we recover the standard version. But, for $\alpha < 1$, the Trinition Minkowski is slightly weaker than the classical one.

## The Trinition weighted Poincare inequality.

Let $\Omega \in \mathbb{R}^3$ be a domain and let $\varphi: \Omega \to \mathbb{R}^3$ be a continuous function with zero boundary values in the Trinition case.

Then we can find a constant depending on $\alpha$ $C_\alpha(\Omega)$ such that,

$$\int_\Omega \|\varphi(X)\|_\alpha^2 \leq C_\alpha(\Omega) \int_\Omega \|\nabla\varphi(X)\|_\alpha^2 dx \tag{74}$$

The proof is obtained from the classical Poincare $\|\varphi\|_{L^2(\Omega)}$.

Then we replace each local Euclidean norm by the Trinition norm. We note in addition that $\lambda_{min}(\alpha)$ ensures $\|\varphi\|_\alpha \approx \sqrt{\lambda_{min}}\|\varphi\|_2$. Here the standard Poincare constant could possibly be extended by a bigger constant $C_\alpha(\Omega)$. This provides us with a stability or coercivity for PDE solution in the Trinition geometry. This will ensure that the energy is always controlled by the gradient in the same geometry.

*The Trinition Grönwall-Freedholm inequality.*

Let us consider time-dependent function $w(t) \in \mathbb{R}^3$ that evolves under a Trinition norm $\|\ \|_{\alpha(t)}$. Let us assume that;

$$\frac{d\|w(t)\|_{\alpha(t)}}{dt} \leq A(t)\|w(t)\|_{\alpha(t)} + B(t) \tag{75}$$

For some positive functions $A(t)$ and $B(t)$. Then,

$$\|w(t)\|_{\alpha(t)} \leq \left(\|w(0)\|_{\alpha(0)} + \int_0^t B(\tau)exp\left(\int_0^\tau A(l)dl\right)d\tau\right) exp\left(\int_0^t A(\tau)d\tau\right) \tag{76}$$

These new inequalities have expanded the toolset for analysing geometry, PDEs and systems under Trinition norms. They bridge classical results with extra cross-terms and partial noncommutativity introduced by $\alpha$.

# Resonant Dancers

In this section, we will introduce new sets and mappings that generalize the Julia and Mandelbrot sets and give rise to a sort of generalized fractals. This, however, is found within the creative multidimensional construct of Trinition, and can perhaps be thought of as the Atangana Resonant Dancer; an entity that defies the formal description of fractals yet possesses the complex, organic nature of some naturally arising fractal system. The growing shape has an incredible potential for motion, rotations, degenerations, and formations of new, unusual shapes and gradually embraces asymmetry the further it develops. It mirrors the showering, transitory formations of the distant nebulae, whose gaseous elements are spawned and dispersed by the forces of gravity as well as the gusts of light from incandescent stars. The billowing, restless flicker of turbulent currents upon an ocean surface; the branching, sometimes-capricious systems of pathways in the networks of life. By that token, the Atangana Resonant Dancer becomes a perfect metaphor for these ghostly structures not just as static attractors but as irrational. Poly-dimensional strings dancing and weaving their way in and out of the gap between chaos and order, forever challenging our idea of symmetry and structure and subliminally suggesting the processes that enable both cosmic realities as well as the very process of life. In this section, we will present some important recursive formula that yield extraordinary fractal shape. We shall first develop some well-known concepts like those that Julia set and Mandelbrot sets in the Trinition space.

## Julia set.

For $z = x + yi + zj$ and a constant

$$C = a + bi + cj$$

$$Z_{n+1} = Z_n^2 + c$$

$$Z^2 = x^2 - y^2 - z^2 + 2xyi + 2xzj + zy(ij + ji)$$

Therefor $Z_{n+1} = Z_n^2 + C$   this leads to

$$\begin{cases} x_{n+1} = x_n^2 - y_n^2 - z_n^2 + z_n^2 + 2y_n z_n(1 - \alpha) + c \\ y_{n+1} = 2x_n y_n + b \\ z_{n+1} = 2x_n z_n - c \end{cases} \tag{77}$$

We will then introduce a new set that will produce outstanding 3D fractal.

## Atangana Trifold vortex set

The Atangana Trifold vortex set will combine recursive conjugation, exponential torsion and parametric self-interference.

$$Z_{n+1} = TrinPower(Z_n, 3) + \alpha \exp(-\beta z_n) + \delta Fold(z_{n-1}) \tag{78}$$

Where

$$TrinPower(Z, m) = Z^m \tag{79}$$

The first component is the non-linear torsion, the second part is a damped conjugate and finally, there is recursive feedback. Trin_power is the 3-Dimentional exponentiation is based on Trinition algebra for which a detailed analysis will be presented later. The following Figure will show the simulation when the component of the fold is set to zero

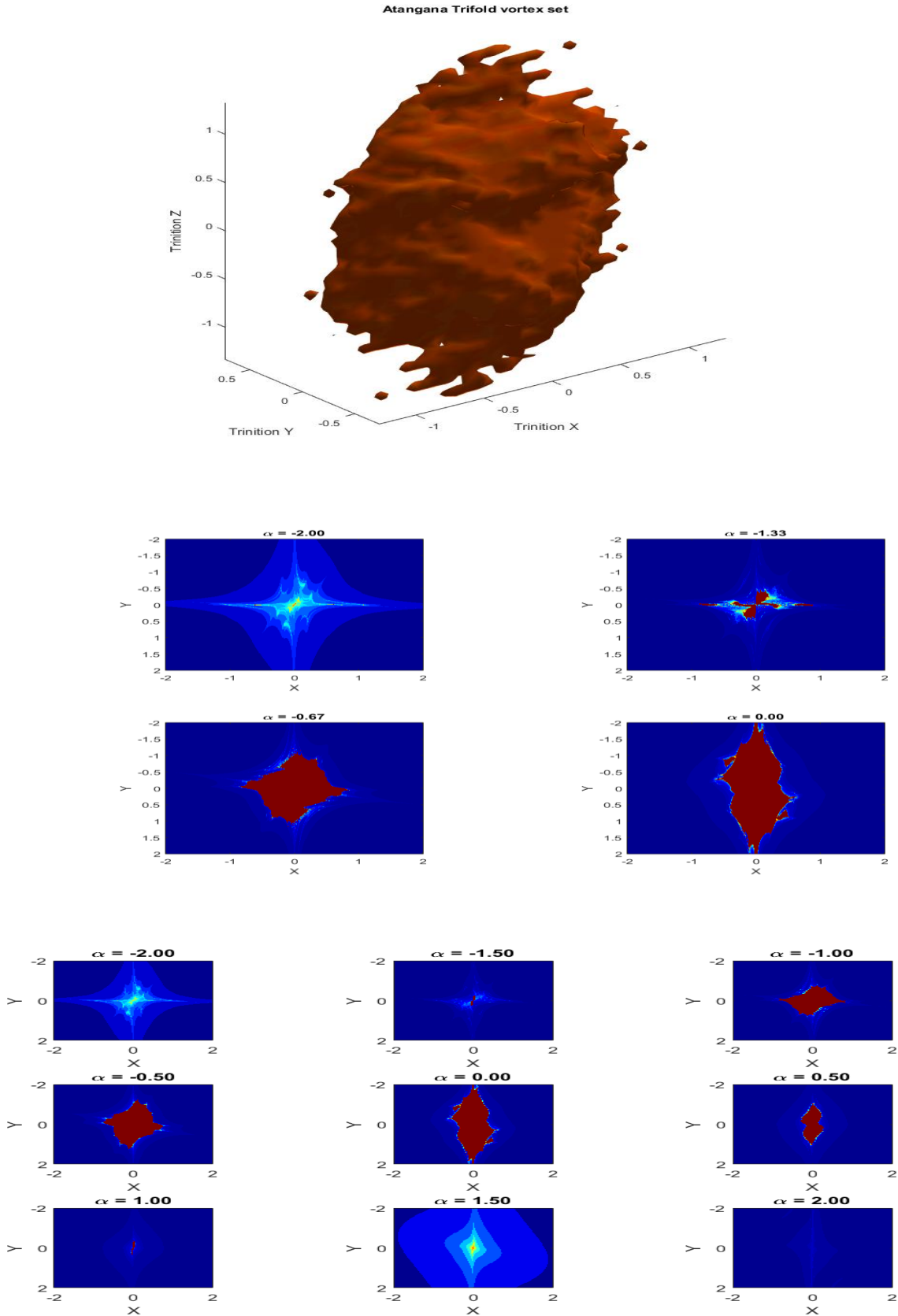

**Figure 8**: Trifold vortex set mapping in 3D and two dimension for different values of alpha

Here we introduce the fold defined as:

$$Fold(z_n, m) = \frac{z_{n-1}}{1 + \|z_{n-1}\|_\alpha^2} \tag{80}$$

The following simulations are obtained

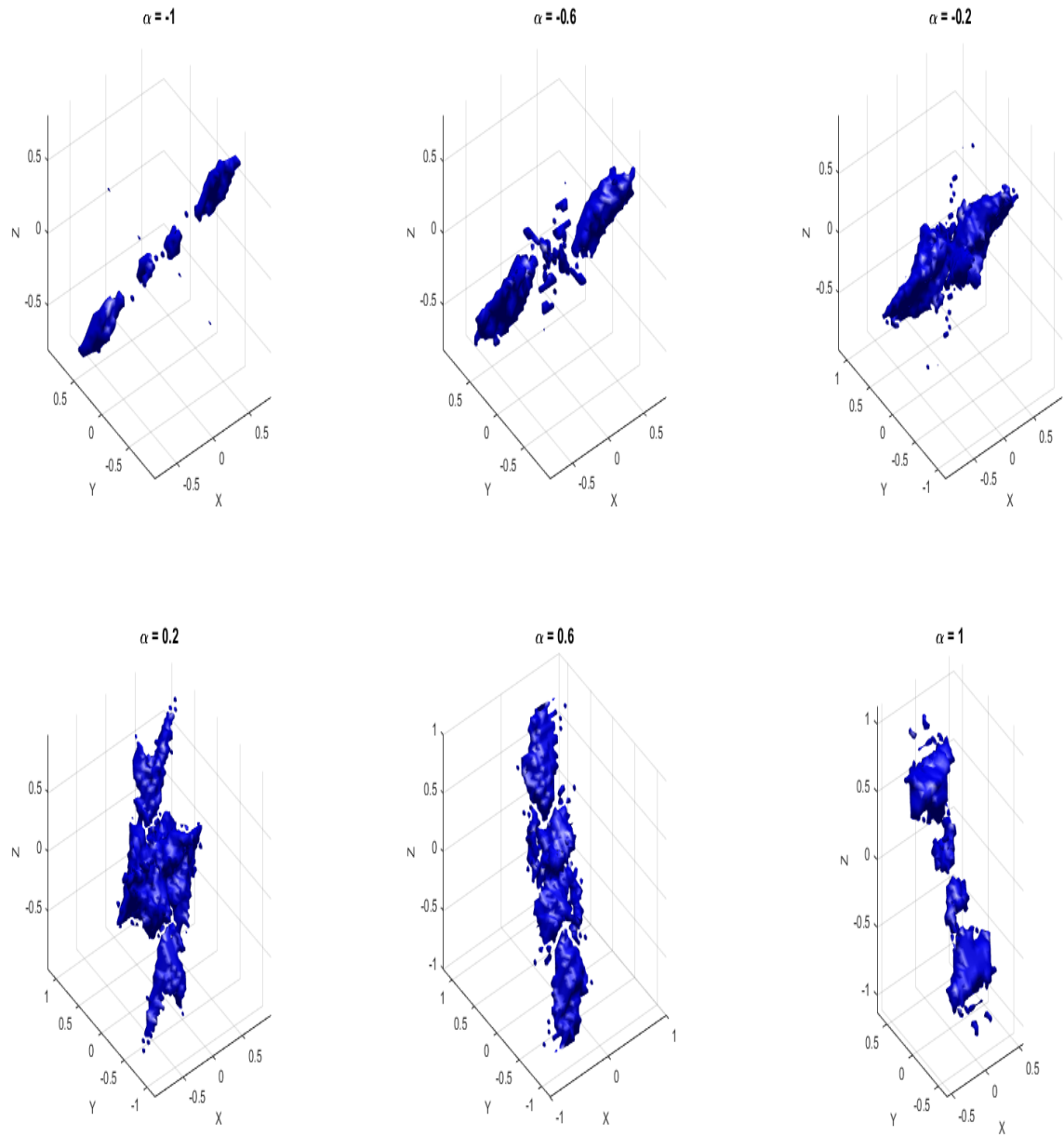

**Figure 9:** 3D Trifold vortex mapping with a fold term for different values of alpha

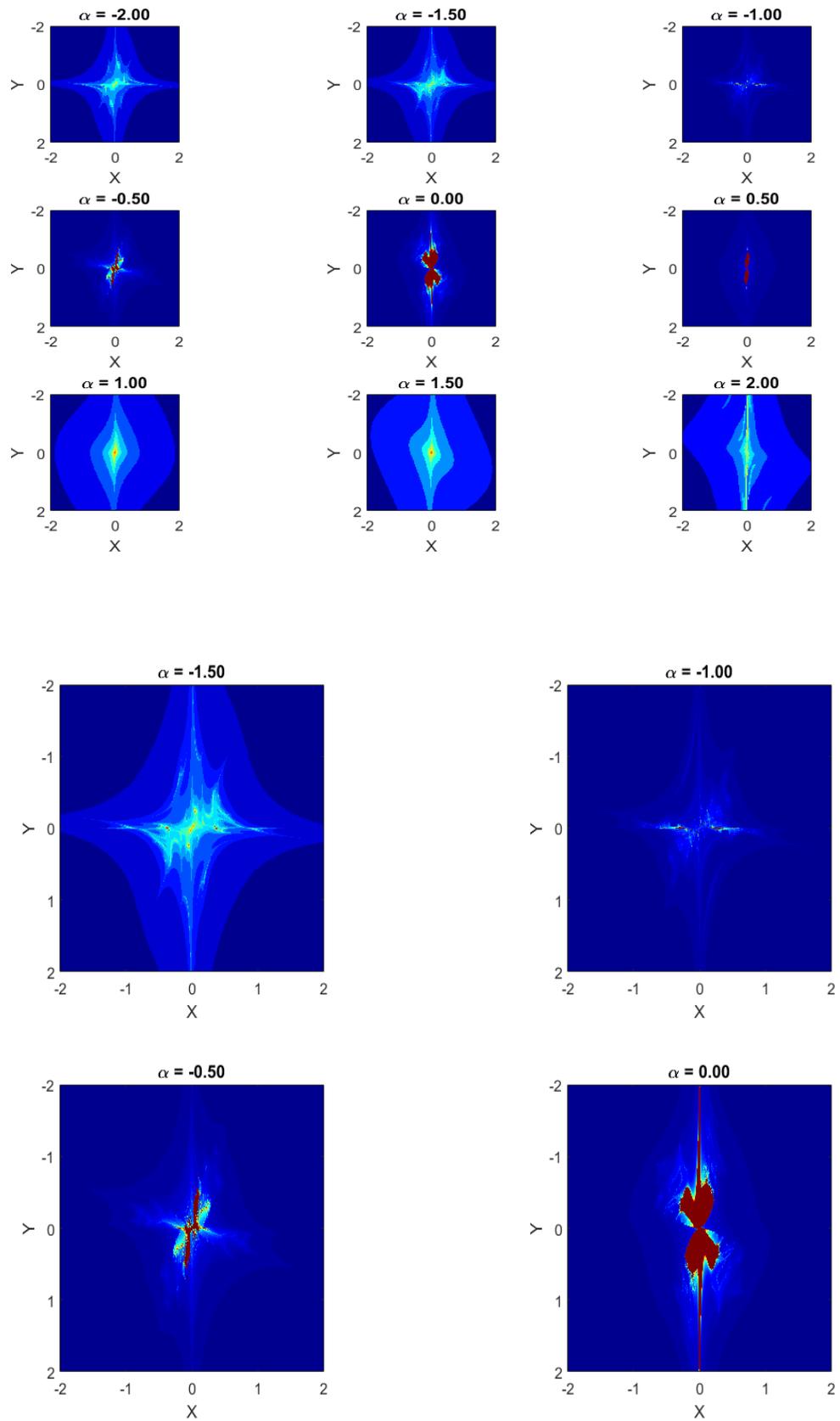

**Figure 10**: 2D Trifold vortex mapping with fold term with different values of alpha

## Atangana polynomial map

We present first a different Atangana resonant dancer that provides structures that we called **Angelus** with reflected component with a formal quadratic map.

$$z_{n+1} = z_n^2 + \alpha \bar{z}_n + c \tag{81}$$

The generated Atangana resonance dancer is derived from the Trinition Fractal. The Trinition Fractal builds on the Mandelbrot-like iterative process but spans three-dimensional space in the Trinition Number System. This Resonance dancer is a 3D extension of the Mandelbrot set, using Trinition algebra. The parameter $\alpha$ adjusts how strongly conjugation is mixed into the iteration, which alters the shape of the fractal.

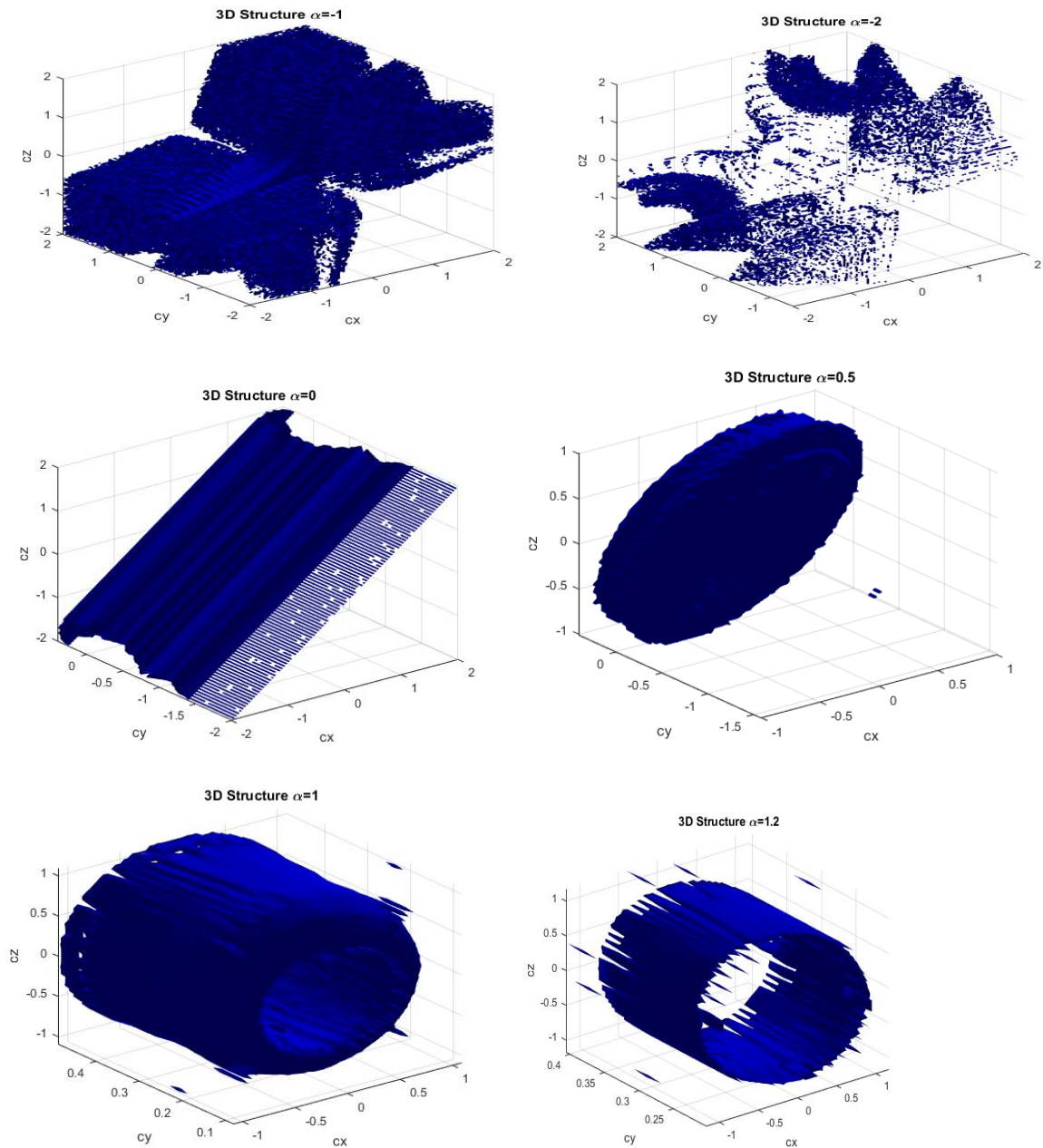

**Figure 11**: 3D Isosurface for different values of alpha (resonance dancer)

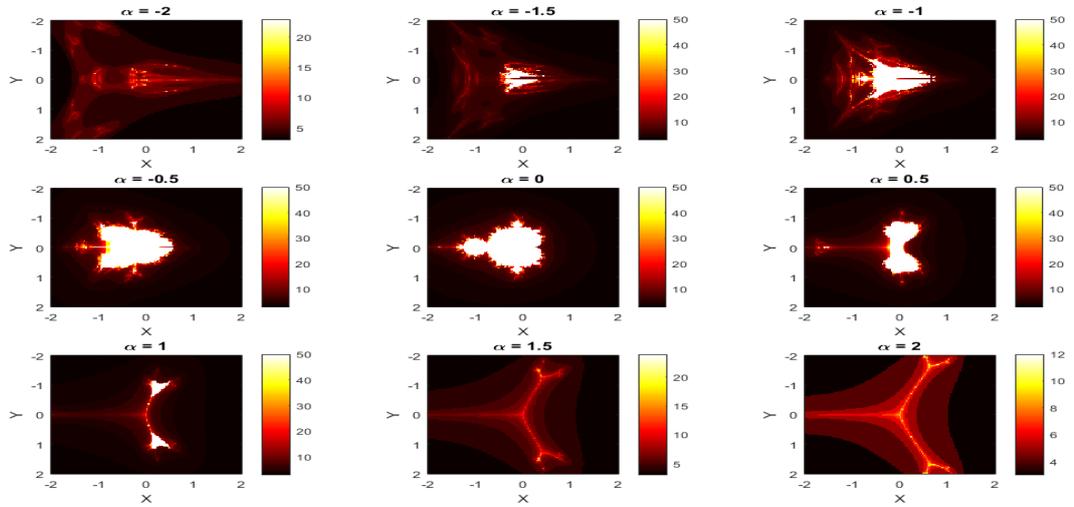

**Figure 12**: Angelus resonance dance with $\alpha \in [-2,2]$

Now by fixing the fractional number that appear in the Trinition equal to 1 and alter the coefficient in front of the conjugate, we obtain the following resonance dance

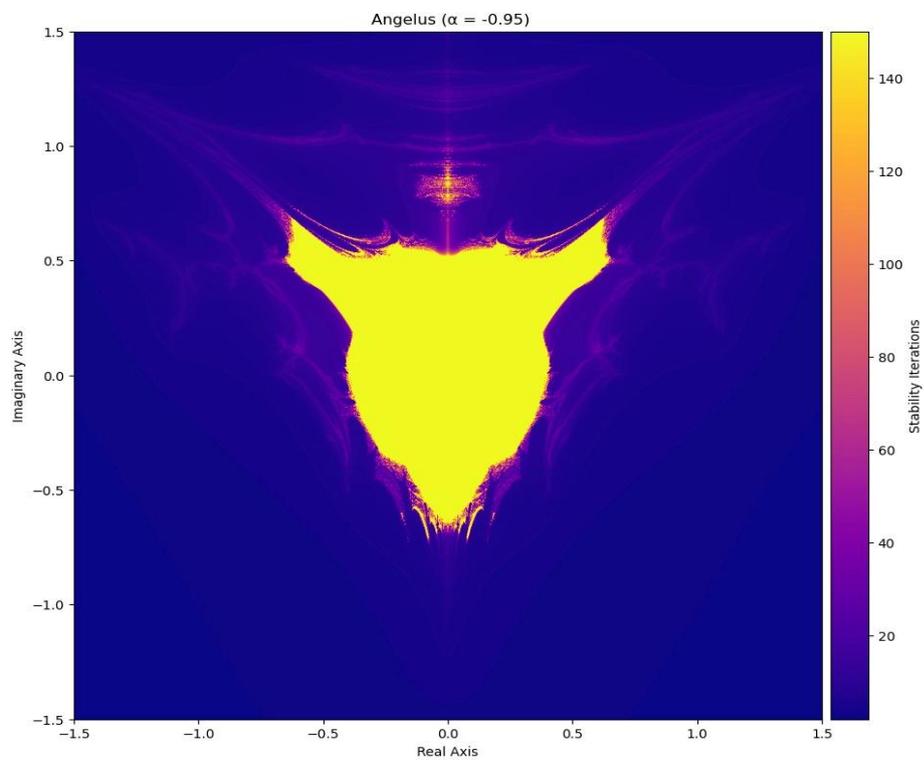

**Figure 13**: Angelus dancing at alpha -0.95

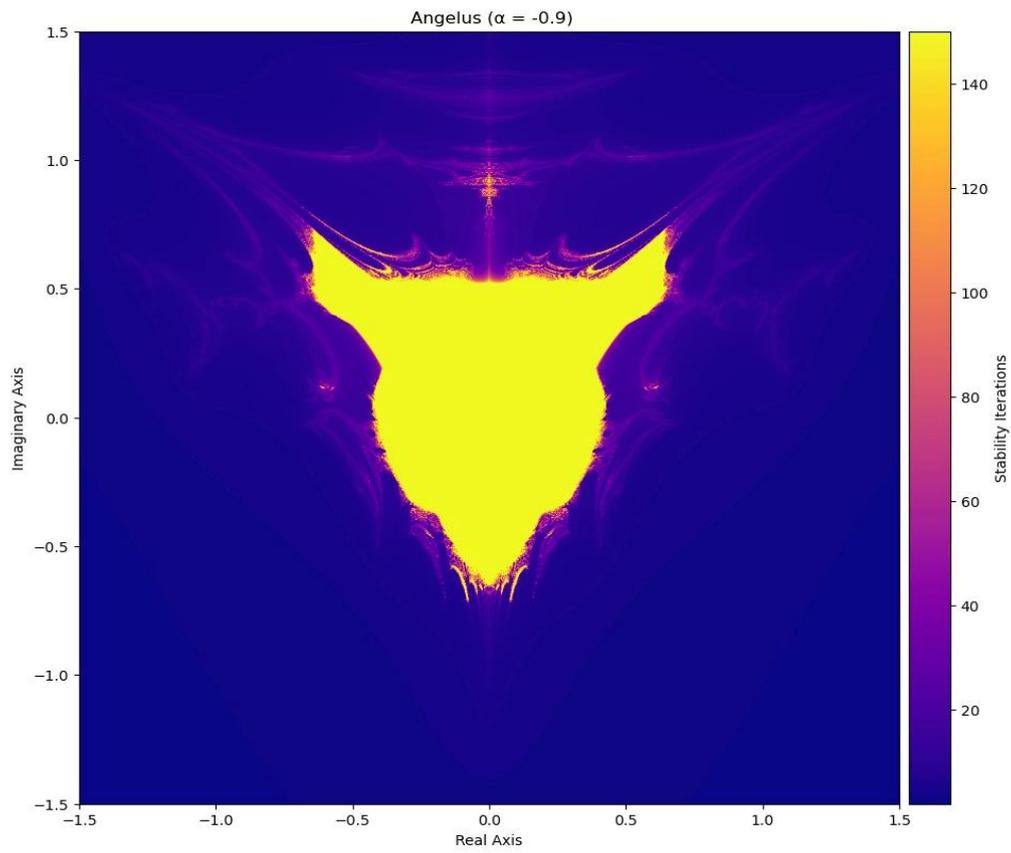

Figure 14: Angelus dancing at alpha -0.9

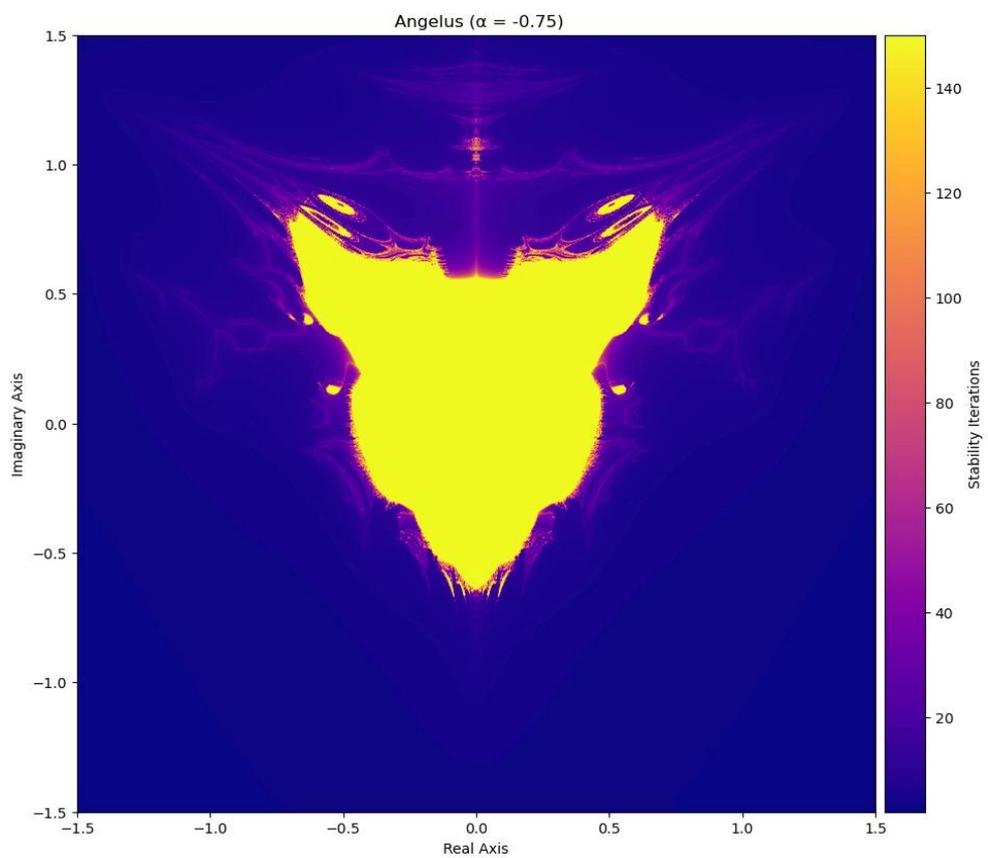

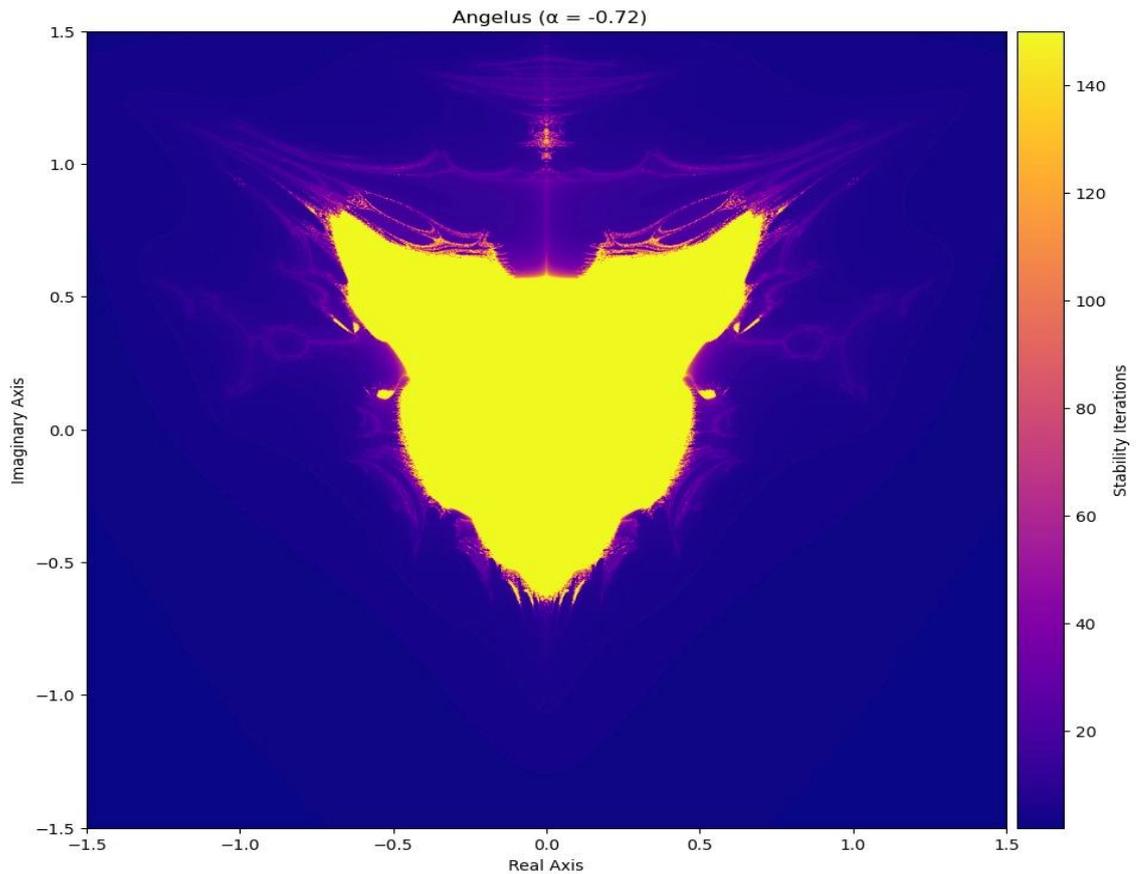

We also suggest the following mapping

$$Z_{n+1} = \sin_\alpha(z_n) + c, \quad \sin_\alpha z = \frac{1}{2u}(\exp(uz) - \exp(-uz)) \tag{82}$$

$\{c \setminus \|z_n\|_\alpha \text{ remains bounded for}\}$

$z_{n+1} = \sin_\alpha(z_n) + c$

Noting that $u = \frac{yi+zj}{r}$

## Atangana conju-quadratic mapping.

Here is now their novel twist on the transition formula, a new variant that replaces a linear conjugate term with its square.

$$z_{n+1} = z_n^2 + (\bar{z}_n)^2 + c \tag{83}$$

By squaring the conjugate, the interplay between standard quadratic growth and the conjugate dynamics is increased; this will yield fractal structure that is unfamiliar to human.

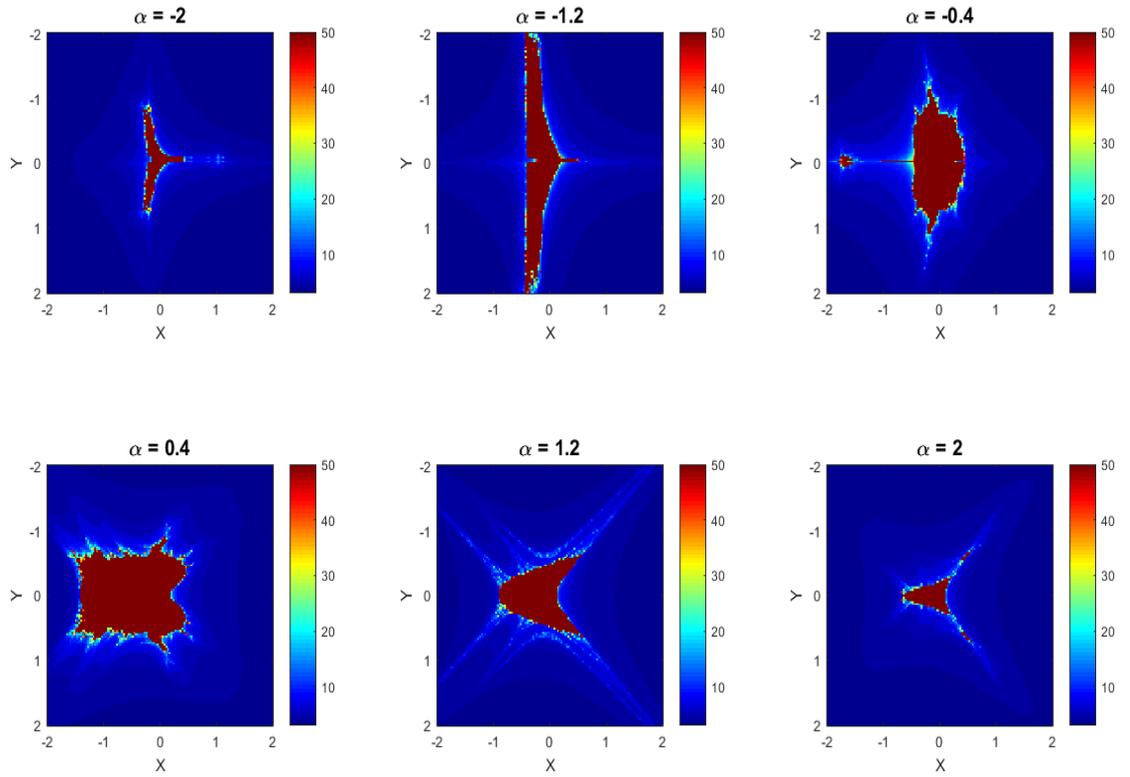

**Figure 15**: Resonance dance for different values of alpha

### Atangana $\alpha$-Morph sets

Here, we will introduce a mapping that will dynamically change $\alpha$ during iteration to introduce fractals that morph between algebraic regimes.

$$\begin{cases} z_{n+1} = z_n^2 + C \\ \sigma_n = \sigma_0 + \beta Re(z_n) \end{cases} \tag{84}$$

Here we will present the three dimension Atangana $\alpha$-Morph sets for different value of alpha. For simplicity also we will replace beta with alpha such that when alpha changes, beta also changes. The numerical simulations are presented below in figure 16 and 17

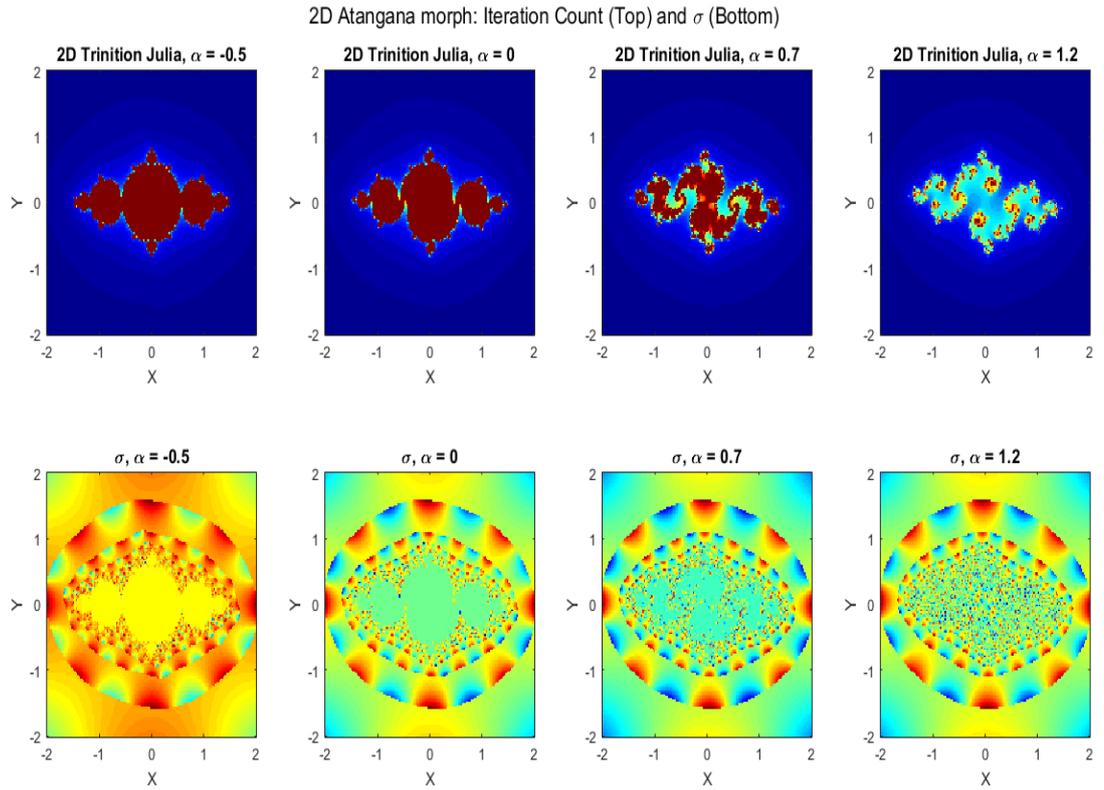

**Figure 16**: *Two dimensional α-Morph structure for different alphas*

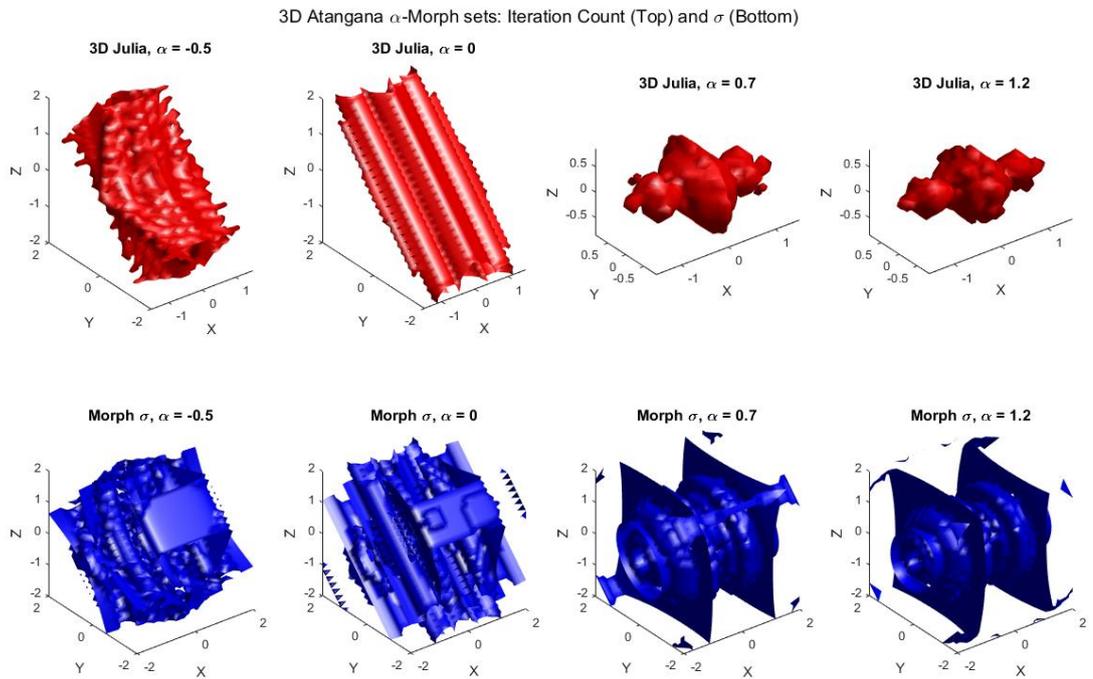

**Figure 17**: 3D dimensional α-Morph structure for different alphas

$$\begin{cases} z_{n+1} = z_n^2 + \beta \bar{z}_n + C \\ \sigma_n = \delta_0 + \beta Re(z_n) \end{cases} \tag{85}$$

$$\begin{cases} z_{n+1} = z_n^2 + \beta (\bar{z}_n)^2 + C \\ \sigma_n = \delta_0 + \beta Re(z_n) \end{cases} \tag{86}$$

$\sigma_0$ is the initial deformation

$\beta$ coupling strength between $\beta$ and $z$

Atangana vortex lattice mapping.

Were the phase modulation is given by $exp(u\gamma\pi)$ that creates a spiral torsion.

### The asymmetric rotative resonance dancer.

To develop this, we shall start with the Trinition multiplication under the Trin_power. We note that the function takes a 3D vector $z(x, y, z)$ and computes a square using.

$$\begin{cases} x_{new} = x^2 + y^2 - (1-\gamma)z^2 \\ y_{new} = 2xy - \gamma(xz - xy) \\ z_{new} = (1+\gamma)xz + \gamma y^2 \end{cases} \tag{87}$$

After the above computation of the new vector, we rotate its $(x, y)$ components by angle to introduce asymmetry and a twist effect

$$\begin{bmatrix} x' \\ y' \\ z' \end{bmatrix} = \begin{bmatrix} \cos(\theta) & -\sin(\theta) & 0 \\ \sin(\theta) & \cos(\theta) & 0 \\ 0 & 0 & 1 \end{bmatrix} \begin{bmatrix} x \\ y \\ z \end{bmatrix} \tag{88}$$

The angle will alters with each iteration as

$$\theta = 0.5 . iteration$$

The iteration for each grid point $c = [x, y, z]$, we start with structure that changes from a crystalline lattice to a chaotic quantum forms the parameter $\gamma$ ranges from $0$ $to$ $1$. This new mapping will combine iteration Trinition with rotation of $Z_{n+1}^2$ plus $coz(\alpha\pi)c$

$$Z_{n+1} = 0.7\{rotate(Z_n^2) + coz(\alpha\pi)c\} \tag{89}$$

The numerical simulation is presented in figure below.

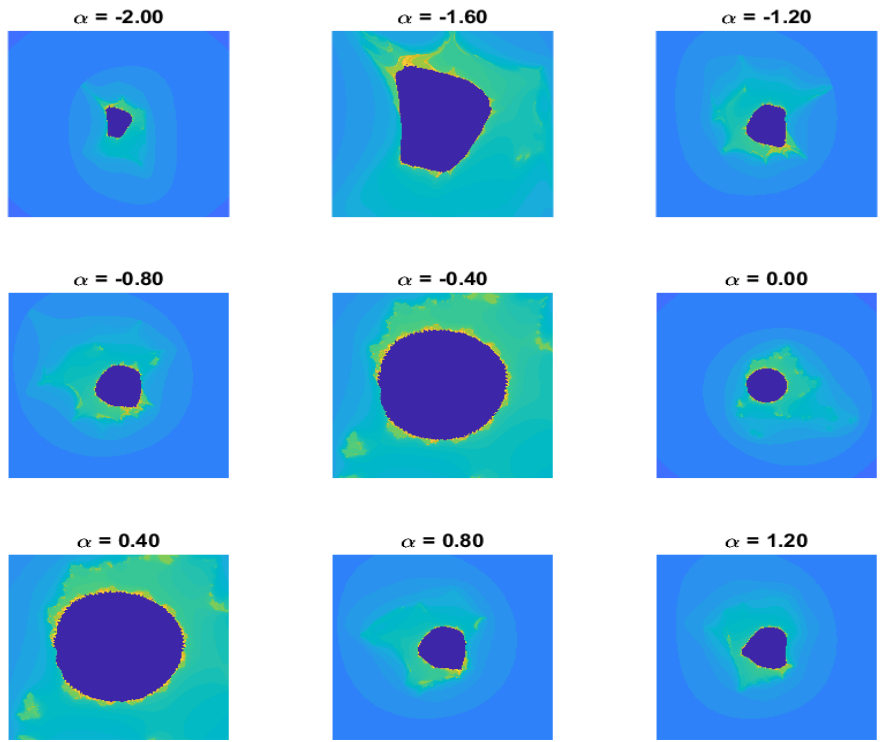

**Figure 18**: Asymmetric resonance dancing

We present some additional resonance dancer and also some new geometrical figures

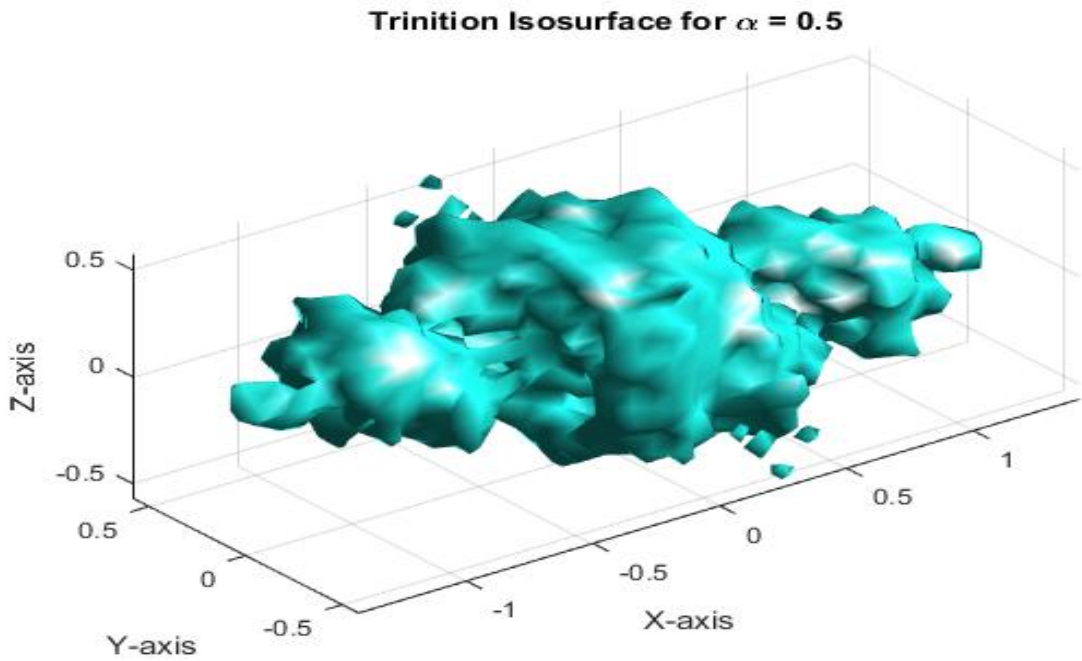

**Figure 19**: Trinition isosurface

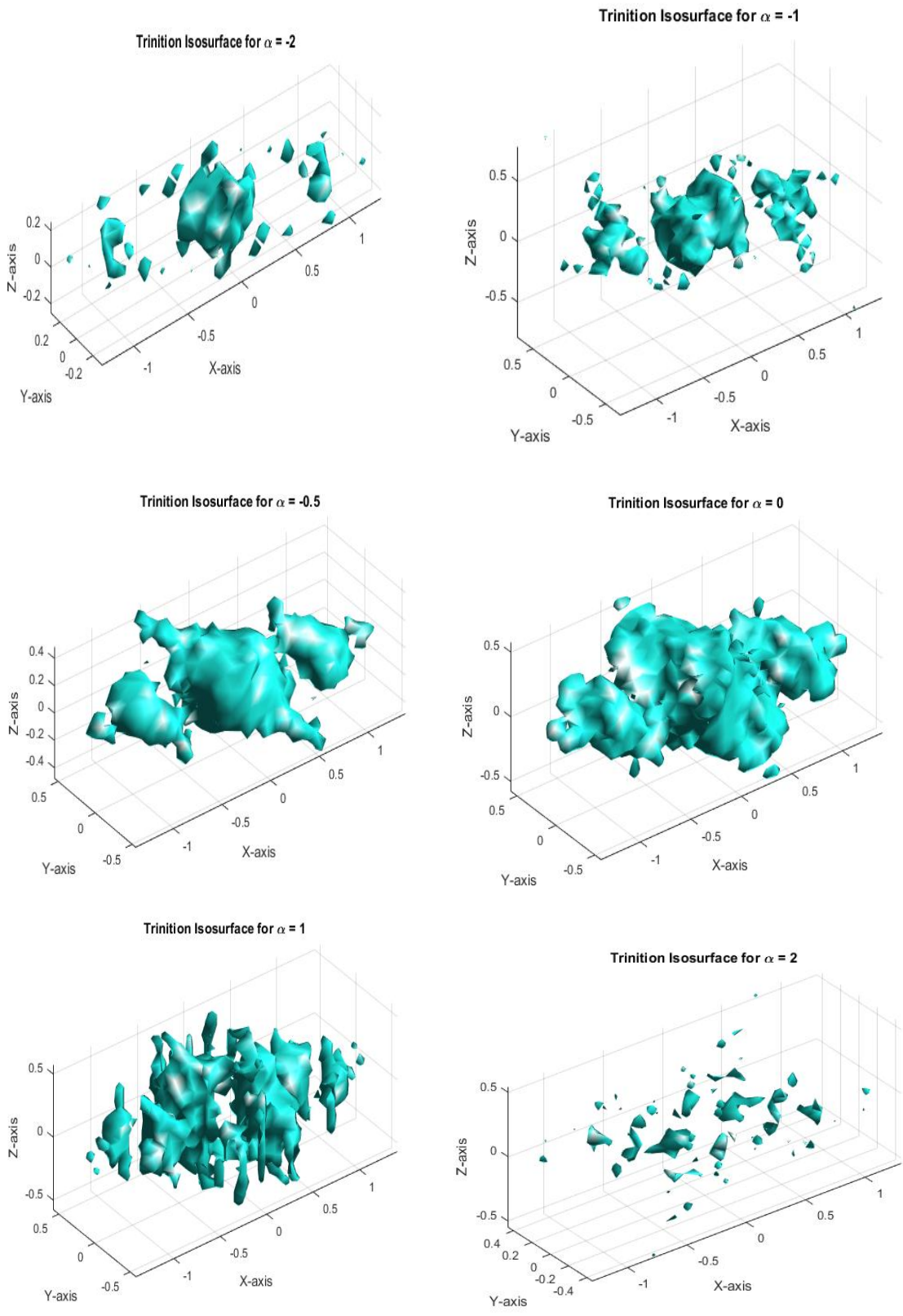

**Figure 20**: Deformable, degenerated and generated 3D isosurface

We present below a new 3D geometrical figure that degenerate, generate, disappears, reappears, deformed accordingly to the value of alpha

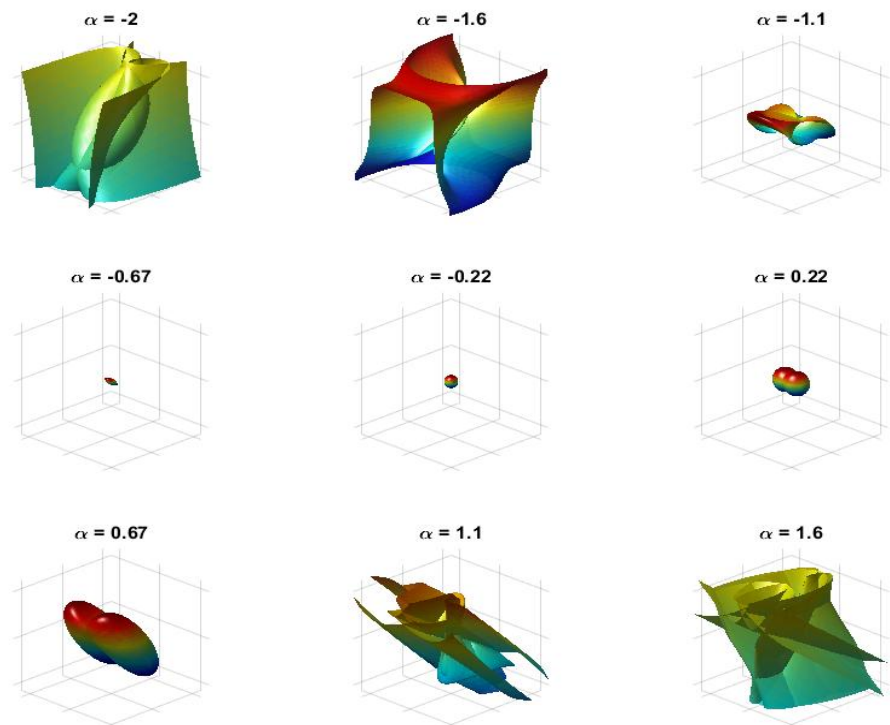

**Figure 21**: Alpha deformable 3D isosurface

We present below a 3D spiral obtained via the Trinition for different values of alpha.

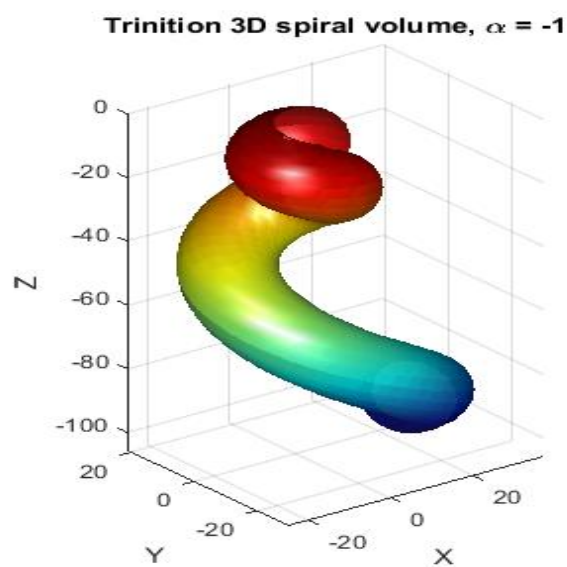

**Figure 22**: A deformable Trinition spiral volume

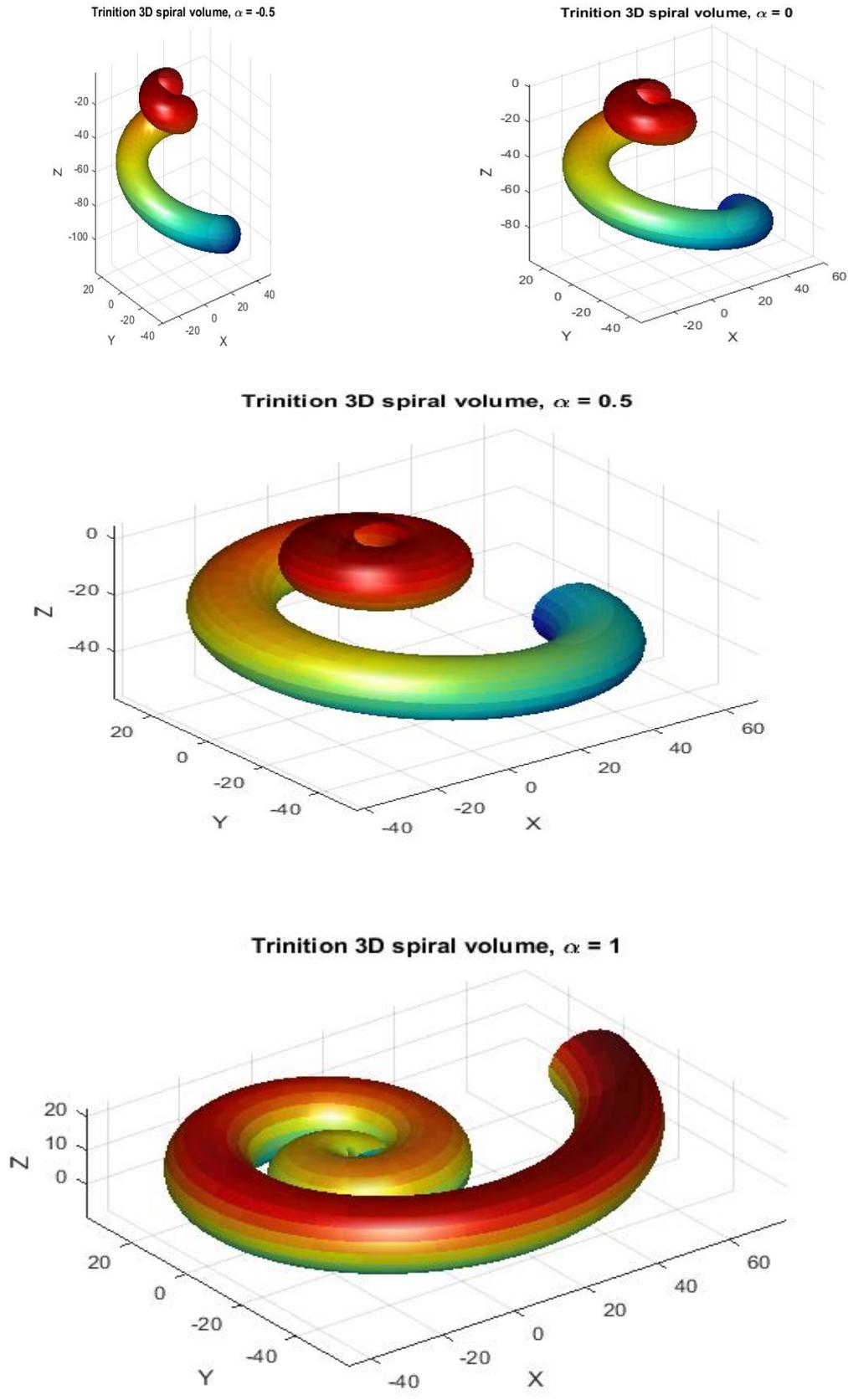

**Figure 23:** Trinition deformable spiral volume for different values of alpha

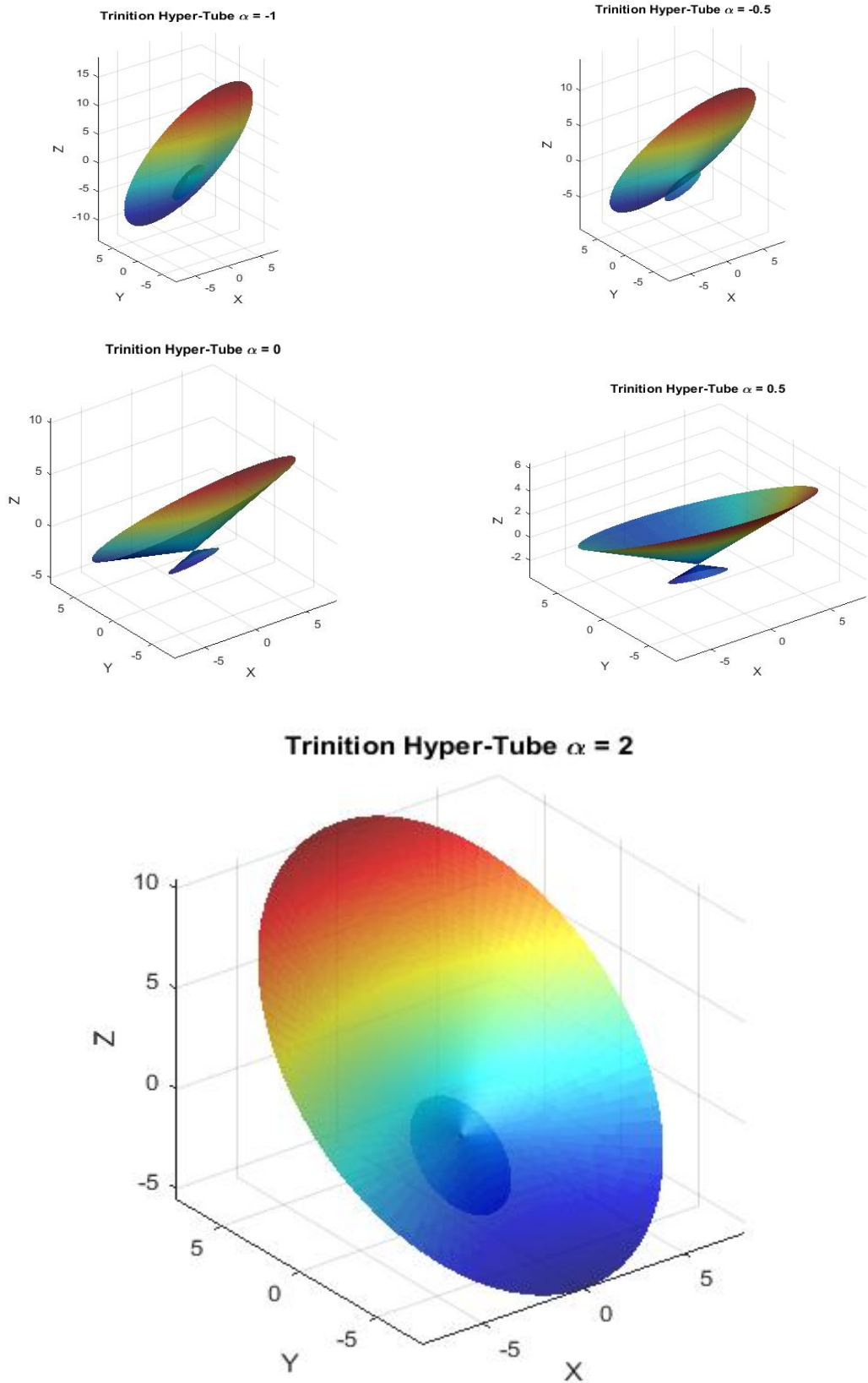

**Figure 24**: Hyper-tube for different values of alpha developed within the framework of Trinition

# Archiometry (corresponding trigonometry with the Trinition)

Here we will consider z with a real x and imaginary part $v = yi + zj$. Naturally, the idea is to view the $x$-direction as a longitudinal or our radical axis consider the plane $(x, z)$ −components as forming our generalized angular plane. Nevertheless, because the product $ij$ is not symmetric naturally, the classical Euclidean quadratic form will be altered to suit the Trinition's algebra structure.

Let us first start by defining a generalized $R$ for $z$. Here the imaginary part will contribute via quadratic form that will incorporate non-commutative corrections. For example, we have that $ij + ji = 2(1 - \alpha)$. Therefore, we will define a modified angular norm by,

$$r_{im} = \sqrt{y^2 + z^2 - 2(1-\alpha)yz} \tag{90}$$

Thus, our candidate for the full modulus is,

$$R_\alpha = \sqrt{x^2 + r_{im}^2} \tag{91}$$

With the above in hand, we will attempt to provide a polar-like decomposition,

$$z = R \cos\theta + u \sin\theta \tag{92}$$

In our case, the unit element ∪ will play the role of a generalized imaginary unit and will be given by,

$$u = \frac{V}{r_{im}} \tag{93}$$

We note that $\theta$ is the generalized angle calculated by the ratio between $x$ and $r_{im}$. This is like the polar angle in spherical coordinates. We recall that the selected definition of $\theta$ accounts for the non-commutative relation between $i$ and $j$ via the parameter $\alpha$. Let us pay attention only on the $y - z$ plane. Let's assume that $x = 0$, then we have a pured imaginary Trinition,

$$Y = yi + zj \tag{94}$$

Here the modulus by,

$$r = \sqrt{y^2 + z^2 - 2(1-\alpha)yz} \tag{95}$$

Here we can introduce an angle $\theta$ so that the unit element in the $y - z$ plane is

$$\cup = \frac{yi + zj}{r} \tag{96}$$

Therefore, the by hypercomplex number $Y$ can be expressed in polar form as;

$$Y = r(\cos\theta + \cup \sin\theta) \qquad (97)$$

This is indeed our direct corresponding Euler formula for complex numbers. However, we should note that the specific form of the trigonometric functions is altered due to the Trinition multiplication. Now, we will introduce the exponential forms in both cases.

The exponential forms in our hypercomplex framework is defined in the classical way through its power series;

$$exp(H) = \sum_{J=0}^{\infty} \frac{H^n}{n!} \qquad (98)$$

This is in general for any hypercomplex number $H$. In particular, for a general Trinition number $z = x + V = x + ix + yj$ one will notice that the real part $x$ commutes with the imaginary part $V$. This indeed helps us obtain the normal splitting,

$$exp(z) = exp(x)exp(V) \qquad (99)$$

We can then follow the footsteps of Euler's formula, that is,

$$exp(V) = \cos\theta + u\sin\theta \qquad (100)$$

Here the function $\cos\theta$ and $\sin\theta$ are defined by their power series. Noting that,

$$V^2 = (yi + zj)^2 = -y^2 - z^2 + yz(ij + ji) \qquad (101)$$
$$= -r_{im}^2 + yz(2(1-\alpha))$$
$$= -r_{im}^2 + 2(1-\alpha)yz$$

In the $y - z$ case, we will obtain a similar case. We have that $Y = yi + zj$, $i^2 = j^2 = -1$, $ij = 1 - \alpha + \alpha k$, $ji = 1 - \alpha - \alpha k$.

The aim is to find the unit $\cup$. To obtain the form;

$$exp(Y) = \cos\theta + \cup \sin\theta \qquad (102)$$

Where $u^2 = -1$, which will correspond to $i^2 = -1$ in the classical case, we compute first $Y^2 = -y^2 - z^2 + 2(1-\alpha)yz$. We rewrite $Y^2$ as,

$$Y^2 = -\theta^2 = -(y^2 + z^2 - 2(1-\alpha)yz) \qquad (103)$$

We define the norm $\cup = \frac{Y}{\theta} = \frac{yi+zj}{\theta}$ then $\cup^2 = -1$.

Now having $\cup^2 = -1$, the power series expansion of $exp(Y)$ provides,

$$exp(Y) = \sum_{n=0}^{\infty} \frac{(\cup \theta)^n}{n!} = \sum_{m=0}^{\infty} \frac{(-1)^m \theta^{2m}}{(2m)!} + \cup \sum_{m=0}^{\infty} \frac{(-1)^m \theta^{2m+1}}{(2m+1)!} \quad (104)$$

We can then recognize the Taylor series of cosine and sine,

$$exp(Y) = \cos \theta + \cup \sin \theta \quad (105)$$

### Defining Trisinus and Triconus

Let us define a Trinition number,

$$T = yi + zj \quad (106)$$

$$Trisinus\ (T) = \sum_{n=0}^{\infty} \frac{(-1)^n}{(2n+1)!} T^{2n+1}$$

$$Triconus\ (T) = \sum_{n=0}^{\infty} \frac{(-1)^n}{(2n)!} T^{2n}$$

Noting that

$$T^2 = -y^2 - z^2 + 2(1-\alpha)yz = -r^2 + +2(1-\alpha)yz$$

$$exp(Y) = Triconus\ (T) + Trisinus\ (T)$$

To capture the full Trinition, we reformulate,

$$z = x + T \quad (107)$$

$$exp(z) = exp(xi) \exp(T) \quad (108)$$

$$Trinus\ (T) = \frac{Trisinus\ (T)}{Triconus\ (T)}$$

Let $T_1, T_2 \in T\alpha$, then,

$$exp(T_1) \exp(T_2) = exp(T_1 + T_2) \quad (109)$$

From the above one will derive for identities.

$$Triconus\ (T_1 + T_2) = Triconus\ (T_1)$$

$$Triconus\ (T_2) - Trisinus\ (T_1)\ Trisinus\ (T_2)$$

$$Trisinus\ (T_1 + T_2) = Trisinus\ (T_1)\ Triconus(T_2) + Triconus(T_1)\ Trisinus\ (T_2)$$

We have that,

$$\frac{d}{d\theta} Triconus(\theta) = -Trisinus\ (\theta)$$

$$\frac{d}{d\theta} Trisinus(\theta) = Triconus\ (\theta)$$

With $z = x + T$.

$$exp(z) = exp(xi)\ exp(T) \tag{110}$$

$$= exp(xi)\ [Triconus(T)\ Trisinus\ (T)]$$

We can see that $exp(xi)$ is the scale magnitude, while the exponential, oscillatory behaviour in the y-plane. We note that if $T = yi + zj$,

$$T^2 = -A^2$$

$T^{2n} = (-A)^n$ and $T^{2n} = (-A)^n\ T$

In $x - y - z$ space,

$$exp(T) = \sum_{n=0}^{\infty} \frac{(T)^n}{n!} = \sum_{n=0}^{\infty} \frac{T^{2n}}{(2n)!} + \sum_{n=0}^{\infty} \frac{T^{2n+1}}{(2n+1)!} = \sum_{n=0}^{\infty} \frac{(-A)^n}{(2n)!} + T \sum_{n=0}^{\infty} \frac{(-A)^n}{(2n+1)!} \tag{111}$$

Therefore,

$$Triconus(T) = \sum_{n=0}^{\infty} \frac{(-A)^n}{(2n)!} \tag{112}$$

$$Trisinus(T) = \frac{T}{\theta} \sin\theta = \frac{T}{\sqrt{A}} \sum_{n=0}^{\infty} \frac{(-A)^n}{(2n+1)!} \tag{113}$$

This shows that indeed if $z = x + T$,

$$exp(z) = exp(xi)\ exp(T) = exp(x)[Triconus(T)\ Trisinus\ (T)] \tag{114}$$

In what follows, in the $y - z$ plane, the squared magnitude is;

$$N_{yz}^2 = y^2 + z^2 - 2(1 - \alpha)yz \tag{115}$$

In the $x - y - z$ space the squared magnitude is;

$$N_{xyz}^2 = x^2 + y^2 + z^2 - 2(1 - \alpha)yz \tag{116}$$

Our goal is to define Triconus, Trisinus and the Trinus.

1) Case 1, $y - z$ plane.

$$y^2 + z^2 - 2(1-\alpha)yz = (y - (1-\alpha)z)^2 + z^2\alpha(2-\alpha) \tag{117}$$

We divide by $y^2 + z^2 - 2(1-\alpha)yz$ on both sides to obtain,

$$\frac{(y-(1-\alpha)z)^2}{N_{yz}^2} + \frac{z^2\alpha(2-\alpha)}{N_{yz}^2} = 1 \tag{118}$$

From the above, we define the Triconus $(\theta)$ as,

$$Triconus\ (\theta) = \frac{y - (1-\alpha)z}{\sqrt{y^2 + z^2 - 2(1-\alpha)yz}} \tag{119}$$

$$Trisinus\ (\theta) = \frac{z\sqrt{\alpha(2-\alpha)}}{\sqrt{y^2 + z^2 - 2(1-\alpha)yz}}. \tag{120}$$

$$\text{Trinus}(\theta) = \frac{y - (1-\alpha)z}{z\sqrt{\alpha(2-\alpha)}} \tag{121}$$

We present here the Trinus, Triconus and the Trisinus as function of $\theta$ for different values of alpha.

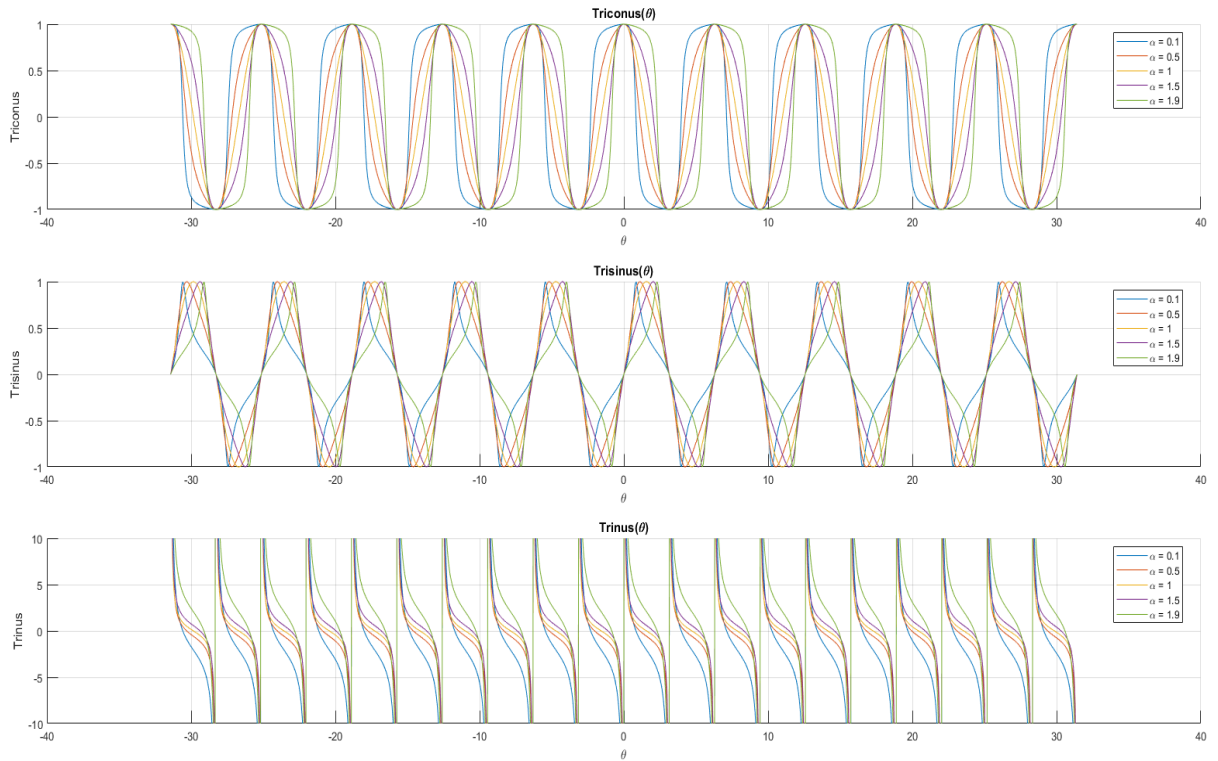

**Figure 25**: Trinus, Triconus and the Trisinus as function of $\theta$ for different values of alpha in y-z subspace .

Now from the above formulas, we can find the archiotwist as:

$$\theta = arctrinus\left(\frac{y-(1-\alpha)z}{z\sqrt{\alpha(2-\alpha)}}\right) \qquad (122)$$

$\theta$ will be called the archiotwist. When $\alpha = 1$, we have, $\theta$ become an angle and Trisinus becomes cosines. We present below the archiotwist for different values of alpha

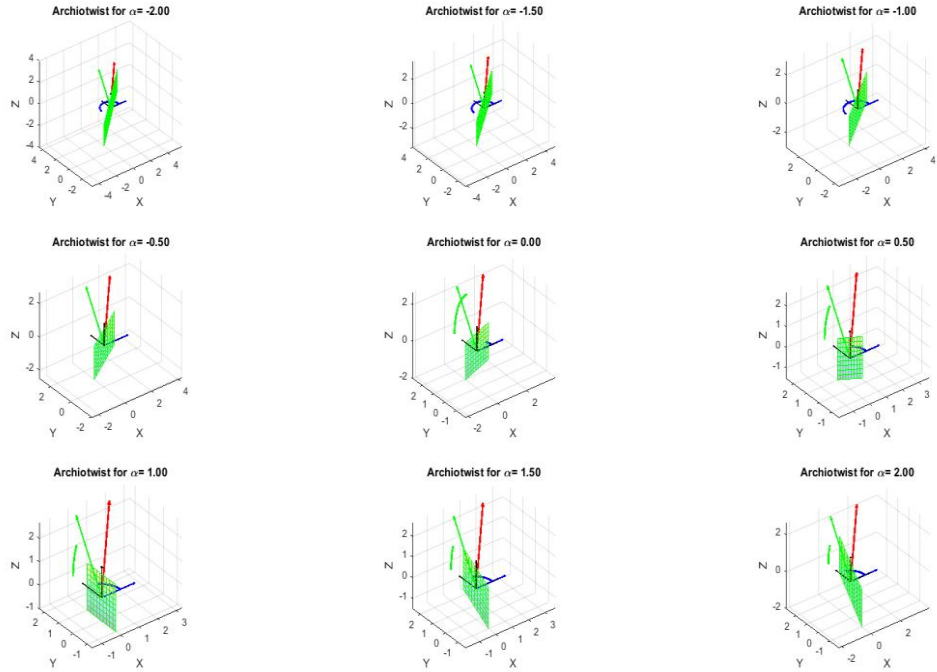

**Figure 26**: An Archiotwists for different values of alpha in y-z subspace

2) $x - y - z$ space.

Here to be careful, we split the coordinate x from the modified y-z subspace. We define a polar Archiotwists such that,

$$Triconus\,(\phi) = \frac{x}{N_{xyz}}, Trisinus\,(\phi) = \frac{N_{yz}}{N_{xyz}},$$

With,

$$N_{xyz} = \sqrt{x^2 + y^2 + z^2 - 2(1-\alpha)yz} \qquad (124)$$

We can verify that;

$$Triconus\ (\phi)^2 + Trisinus\ (\phi)^2 = \frac{x^2 + N_{yz}^2}{N_{xyz}^2} \quad (124)$$

$$= \frac{x^2 + (y^2 + z^2 - 2(1-\alpha)yz)}{x^2 + (y^2 + z^2 - 2(1-\alpha)yz)} = 1$$

Therefore, in $x - y - z$, we will have,

$$Triconus\ (\phi) = \frac{x}{\sqrt{x^2 + (y^2 + z^2 - 2(1-\alpha)yz)^2}} \quad (125)$$

$$Trisinus\ (\phi) = \frac{\sqrt{y^2 + z^2 - 2(1-\alpha)yz}}{\sqrt{x^2 + (y^2 + z^2 - 2(1-\alpha)yz)^2}} \quad (126)$$

$$Trinus\ (\phi) = \frac{\sqrt{y^2 + z^2 - 2(1-\alpha)yz}}{x}$$

We will present between the defined function depending to the variable $\phi$ for different values of alpha.

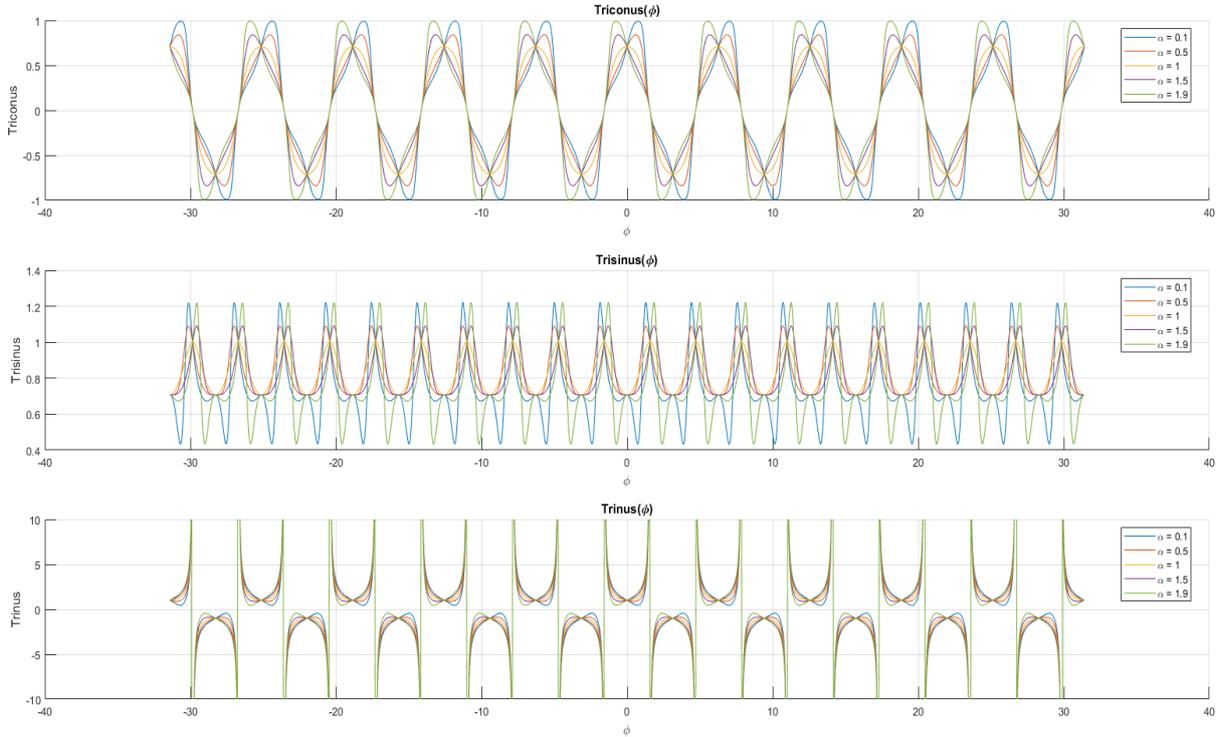

**Figure 27**: Trinus, Triconus and the Trisinus as function of $\theta$ for different values of alpha in x-y-z subspace

From the above formula, we present the archiotwist in x-y-z space.

$$\phi = artTrinus\left(\frac{\sqrt{y^2 + z^2 - 2(1-\alpha)yz}}{x}\right)$$

We present below the angle for different values of alpha.

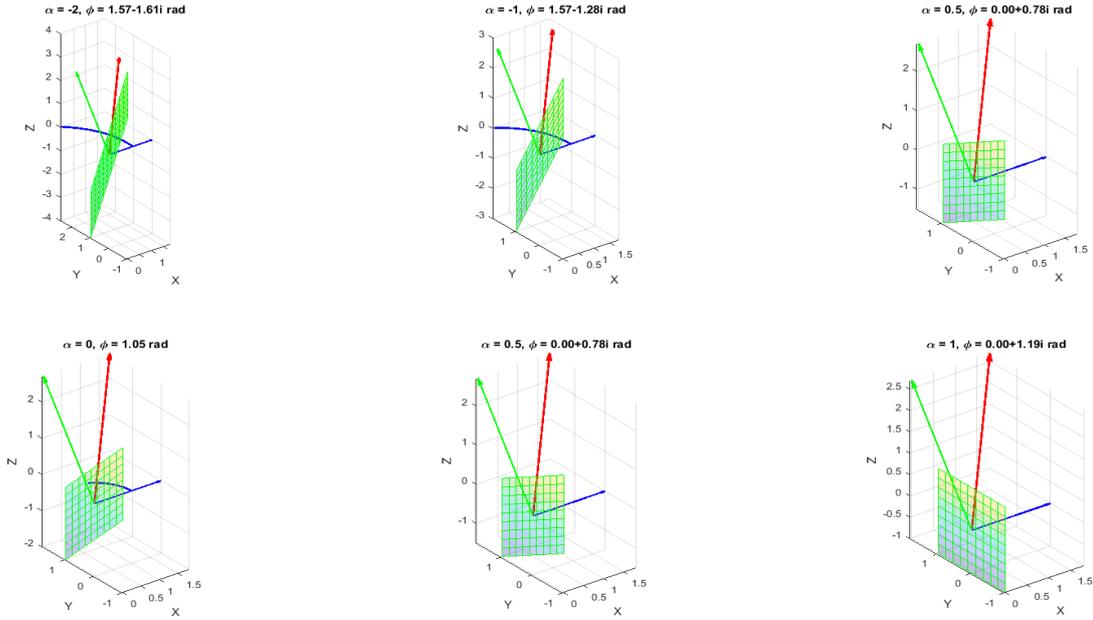

**Figure 28**: Archiotwists for different values of alpha in x-y-z

We can see that in both cases, Triconus$^2$+Trisinus$^2$= 1.

We note that ∪ can also be defined as;

$$\cup = \frac{(y-(1-\alpha)z)i + \sqrt{\alpha(2-\alpha)}zj}{\sqrt{y^2 + z^2 - 2(1-\alpha)yz}} \quad (127)$$

This will be called an effective imaginary unit for the $y - z$ plane. Indeed ∪$^2$= −1. We can now individually obtain a Euler-type formula in $y - z$ subspace,

$$exp(\cup \theta) = Triconus\,(\theta) + Trisinus\,(\theta) \quad (128)$$

But, in $x - y - z$ space, we have,

$$z = R\,exp(i\theta)\,exp(\cup \phi) \quad (129)$$

Here we have that,

- $exp(\cup \phi) = Triconus\ (\phi) +\cup Trisinus\ (\phi)$ this accounts for the rotation in the deformed $y - z$ plane. On the other hand,
- $exp(i\theta) = Triconus\ (\theta) + iTrisinus\ (\theta)$ this will handle the longitudinal part along the $x - axis$. We note that in the representation of x-direction imaginary unit. This is consistent since the overall algebra is non-commutative, so the order maters. Indeed, this double-exponential configuration is wiser since it separates our geometry into two orthogonal parts. One captures the defamation of rotation in the $y - z$ plane and the other the standard rotation in the x-axis.

We can see that,

$$z^n = r^n \exp(in\theta) \exp(n \cup \phi) \tag{130}$$

These two definitions are indeed the generalization of Euler and Moivre's theorem.

Let us find the $nth$ roots within the Archiometric, to do this, we will start with the generalized De Moivre's formula,

$$z^n = r^n \exp(in\theta) \exp(n \cup \phi) \tag{131}$$

We seek solutions to;

$$\lambda^n = z = r \exp(i\theta) \exp(\cup \phi) \tag{132}$$

Let us set;

$$\lambda = \Omega \exp(i\bar{\theta}) \exp(\bar{\phi} \cup) \tag{133}$$

We take the nth power, we get,

$$\left(\Omega \exp(i\bar{\theta}) \exp(\bar{\phi} \cup)\right)^n = \Omega^n \exp(in\bar{\theta}) \exp(n \cup \bar{\phi}) \tag{134}$$

Settling equal to 2 leads to;

A. $\Omega^n = r \rightarrow \Omega = r^{\frac{1}{n}}$

B. $n\bar{\theta} = \theta + 2k\pi \rightarrow \bar{\theta}_R = \frac{\theta + 2k\pi}{n}, k = 0,1 \dots n - 1$

C. $n\bar{\phi} = \phi \rightarrow \bar{\phi} = \frac{\phi}{n}$

Therefore, the nth roots of z are,

$$\lambda = r^{\frac{1}{n}} \exp\left(\frac{i\theta + 2k\pi}{n}\right) \exp\left(\cup \frac{\phi}{n}\right), k = 0,1, \dots n - 1$$

# The Dance of Motion: Velocity, Energy and Acceleration in Trinition

Let $z(t) = x(t) + y(t)i + z(t)j$

We differentiate $z(t)$ to obtain,

$$v(t) = \frac{\partial z(t)}{\partial t} = \dot{x}(t) + \dot{y}(t)i + \dot{z}(t)j \tag{135}$$

The above definition replaces the instantaneous rate of change in each coordinate. We will assume for simplicity that,

$$z(t) = R(t)\sqrt{x^2(t) + \theta^2(t)} \tag{136}$$

Where

$$\theta^2(t) = y^2(t) + z^2(t) - 2(1-\alpha)y(t)z(t) \tag{137}$$

$$\cup(t) = \frac{(y(t) - (1-\alpha)z(t))i + \sqrt{\alpha(2-\alpha)}z(t)j}{\sqrt{y^2(t) + z^2(t) - 2(1-\alpha)y(t)z(t)}} \tag{138}$$

$$A(t) = \frac{\partial u(t)}{\partial t} = \ddot{x}(t) + \ddot{y}(t)i + \ddot{z}(t)j \tag{137}$$

$$\|V(t)\|_\alpha^2 = \dot{x}^2(t) + \dot{y}^2(t) + \dot{z}^2(t) - 2(1-\alpha)\dot{y}(t)\dot{z}(t) \tag{138}$$

Therefore, in the Trinition space, the kinetic energy becomes,

$$E_K^\alpha = \frac{1}{2}m[\dot{x}^2(t) + \dot{y}^2(t) + \dot{z}^2(t) - 2(1-\alpha)\dot{y}(t)\dot{z}(t)] \tag{139}$$

Now, consider a particle of mass M that moves in a potential energy $E_\alpha^P$. One of the common way is to define the Lagrangian. However, for simplicity here, we assume a harmonic oscillator potential. That is, we will let,

$$E_\alpha^P(z) = \frac{1}{2}K\|z(t)\|_\alpha^2 \tag{140}$$

$$E_\alpha^P(z) = \frac{1}{2}K[x^2(t) + y^2(t) + z^2(t) - 2(1-\alpha)y(t)z(t)]$$

The kinetic energy T can be defined using the squared speed computed earlier.

Thus,

$$T = \frac{1}{2}m\|V\|_\alpha^2 = \frac{1}{2}m[\dot{x}^2(t) + \dot{y}^2(t) + \dot{z}^2(t) - 2(1-\alpha)\dot{y}(t)\dot{z}(t)] \quad (141)$$

Here we can obtain our Lagrangian as,

$$L = T - E_\alpha^P = \frac{1}{2}m[\dot{x}^2 + \dot{y}^2 + \dot{z}^2 - 2(1-\alpha)\dot{y}\dot{z}] \quad (142)$$
$$-\frac{1}{2}K[x^2 + y^2 + z^2 - 2(1-\alpha)yz]$$

Now we will derive the Euler-Lagrange equations.

*For the x-coordinate,*

$$\frac{\partial L}{\partial x} = -kx, \frac{\partial L}{\partial \dot{x}} = m\dot{x} \quad (143)$$

We will for the x-direction have,

$$\frac{\partial(m\dot{x})}{\partial t} - (-kx) = m\ddot{x} + kx = 0 \quad (144)$$

Therefore, we will have,

$$m\ddot{x} + kx = 0 \quad (145)$$

*For the y-coordinate,*

$$\frac{\partial L}{\partial y} = -kx + k(1-\alpha)z, \frac{\partial L}{\partial \dot{y}} = m\dot{y} - m(1-\alpha)\dot{z} \quad (146)$$

Then,

$$\frac{\partial}{\partial t} = (m\dot{y} - m(1-\alpha)\dot{z}) = m\ddot{y} - m(1-\alpha)\ddot{z} \quad (147)$$

The Euler-Lagrange equation for y is,

$$m\ddot{y} - m(1-\alpha)\ddot{z} + ky - k(1-\alpha)z = 0 \quad (148)$$

*For the Euler-Lagrange z-coordinate, we will have,*

$$m\ddot{z} - m(1-\alpha)\ddot{y} - k(1-\alpha)y = 0 \quad (149)$$

Let us decouple the y and z equations. Let us set $F = 1 - \alpha$, and $S(t) = y(t) + z(t)$, then we will obtain,

$$m(\ddot{y} + \ddot{z}) - mF(\ddot{y} + \ddot{z}) + k(y + z) - kF(z + y) = 0 \quad (150)$$

$$m(1 - F)(\ddot{y} + \ddot{z}) + k(1 - F)(z + y) = 0 \quad (151)$$

We have that $1 - F = 1 - 1 + \alpha \neq 0$ unless $\alpha = 0$, then,

$$\ddot{\delta} + \frac{k}{m}\delta = 0 \tag{152}$$

$\bar{\delta}(t) = y(t) - z(t)$, we obtain,

$$\ddot{\bar{\delta}} + \frac{k}{m}\bar{\delta} = 0 \tag{153}$$

The coupling of equations y and z leads into two independent harmonic oscillators.

For the canonical momenta and velocity-momentum relations.

$$P_x = m\dot{x}, P_y = m\dot{y} - m(1-\alpha), P_y = m\dot{z} - m(1-\alpha)\dot{y}, F = 1 - \alpha \tag{154}$$

Then the relations for the y and z components becomes,

$$\begin{pmatrix} P_y \\ P_z \end{pmatrix} = m \begin{pmatrix} 1 & -F \\ -F & 1 \end{pmatrix} \begin{pmatrix} \dot{y} \\ \dot{z} \end{pmatrix} \tag{155}$$

$$M = m \begin{bmatrix} 1 & -F \\ -F & 1 \end{bmatrix} \tag{156}$$

$$det = m^2(1 - F^2) \tag{157}$$

Therefore, for $\alpha = 0$ and $\alpha \neq 2$ the inverse exists.

$$\dot{y} = \frac{P_y + FP_3}{m(1 - F^2)}, \dot{z} = \frac{FP_y + P_3}{m(1 - F^2)}, \dot{x} = \frac{P_x}{m} \tag{158}$$

Hamiltonian formulation: The classical Hamiltonian is defined by,

$$H = P_x\dot{x} + P_y\dot{y} + P_z\dot{z} - L \tag{159}$$

$$P_x\dot{x} = \frac{P_x^2}{m} \tag{160}$$

And for y and z parts;

$$P_y\dot{y} + P_z\dot{z} = \frac{1}{m(1 - F^2)}\left(P_y^2 + 2FP_yP_3 + P_z^2\right) \tag{161}$$

Then the kinetic energy in terms of momenta is,

$$T = \frac{1}{2}m[\dot{x}^2 + \dot{y}^2 + \dot{z}^2 - 2F\dot{y}\dot{z}] \tag{162}$$

$$= \frac{P_x^2}{2m} + \frac{P_y^2 + P_z^2 + 2FP_yP_z}{2m(1-F^2)} \tag{163}$$

Now the Hamiltonian is Trinition is given by;

$$H = \frac{P_x^2}{2m} + \frac{P_y^2 + P_z^2 + 2(1-\alpha)P_yP_z}{2m(1-(1-\alpha)^2)} + \frac{1}{2}K[x^2 + y^2 + z^2 - 2(1-\alpha)yz] \tag{164}$$

Now we will show some numerical simulations of these new defined functions, in Trinition Space, which is involved in the motion of a harmonic oscillator due to a deformed geometric structure. Our simulation aims to investigate the effect of introducing a deformation parameter to the effective metric in the model which, in turn, alters familiar quantities, such as position norms, velocity magnitudes, and the conservation of energy. We derive trajectories of the three spatial coordinates evolving in a harmonic form with their deformations included for their interactions by setting up initial conditions for the system. The numeric technique involves calculating velocities and accelerations via differencing as well as reciprocated norms that incorporate the new space geometries. Different plots, including two-dimensional relationships, the position and velocity over time, and the phase-space plots are shown to demonstrate how the system works and the relationship of different variables as the system changes. It also creates an animation which shows how the three-dimensional path of the oscillator changes as a function of time to help visualize the effect of the deformation parameter. This geometric significance should not be underestimated to this end, the value of $\alpha$ imparts a controlled deformation to Trinition Space. This parameter specifically changes how these coordinates interact with one another i.e., how the y and z components relate to each other through a coupling term that now depends on $(1-\alpha)$. For $\alpha = 1$, the standard Euclidian geometry is recovered as expected, and for all other values the modified norms and dynamical properties are far from classical expectations. Reducing α leads to a larger influence of the coupling term, impacting the oscillatory motion, magnitudes of velocities, and the system's energy balance. Such phenomena manifest in the calculated norms and phase-space trajectories, revealing the impact of the deformation parameter on the essential features of motion and energy dispersion within this extended paradigm.

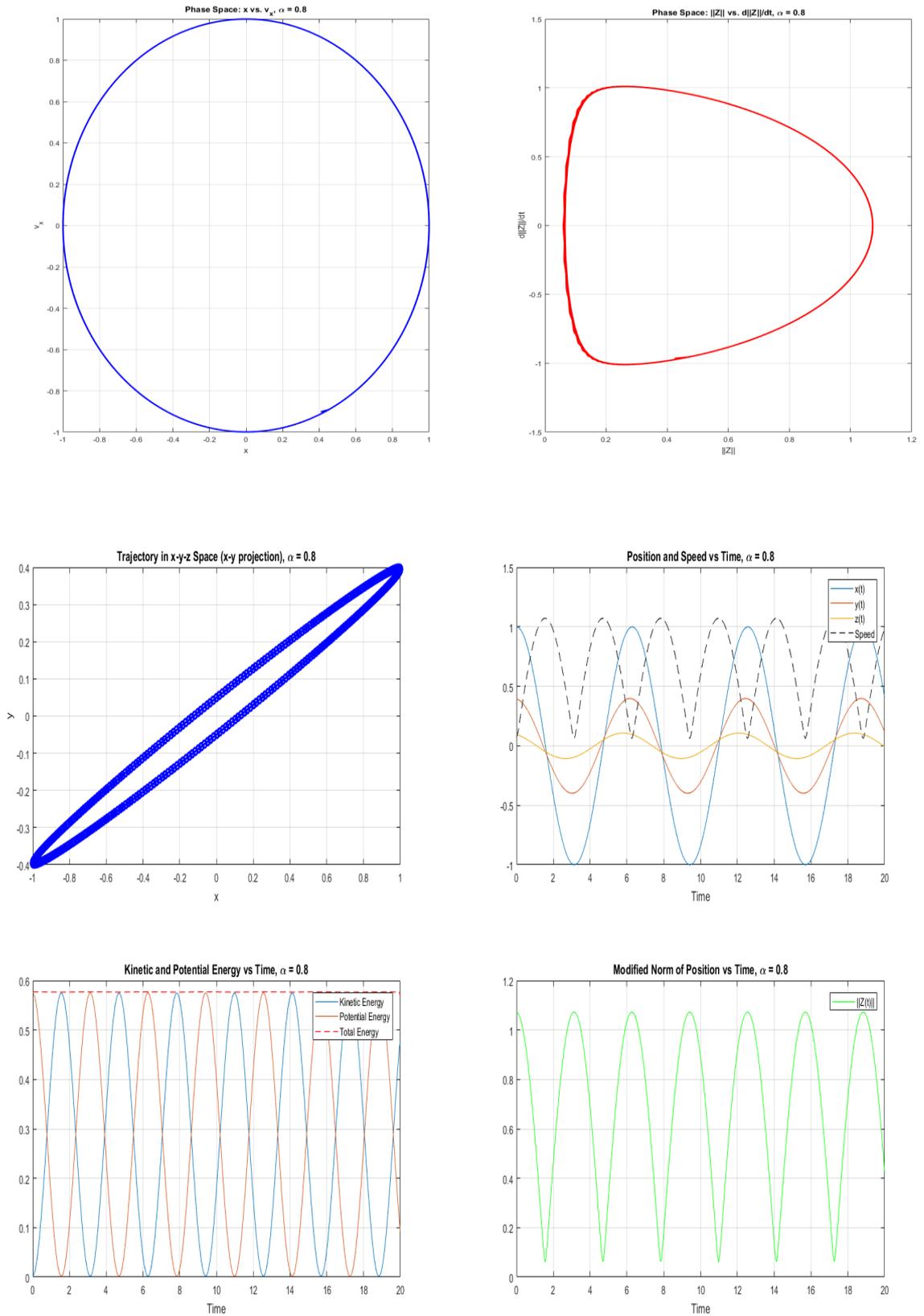

**Figure 29**: Simulation of equations of motion within the Trinition framework for alpha =0.8

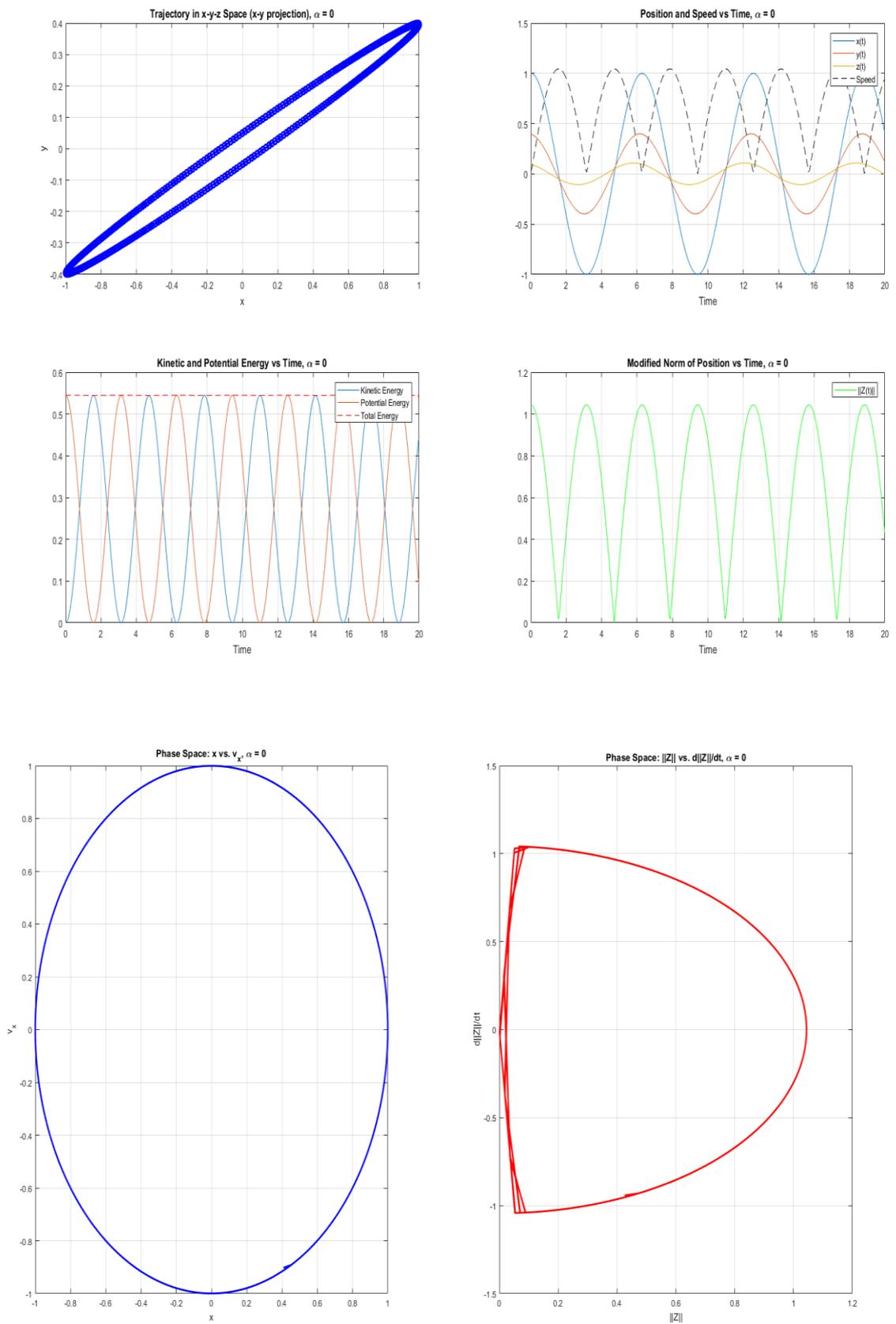

**Figure 30**: Simulation of equations of motion within the Trinition framework for alpha =0.8

## Integral transform

Let us start with Laplace transform, let us consider the Trinition Laplace variable $s$. This $s$ could be a Trinition number or at least include the appropriate deformed imaginary unit, then,

$$\mathcal{L}_\alpha(f(t))(s) = \int_0^\infty f(t)\exp(-st)\, dt \tag{165}$$

$$\int_0^\infty f(t)[Triconus(st) - Trisinus(st)]\, dt \tag{166}$$

$$Triconus(st) = \sum_{n=0}^\infty \frac{(-1)^n (st)^{2n}}{(2n)!} \tag{167}$$

$$Trisinus(st) = \sum_{n=0}^\infty \frac{(-1)^n (st)^{2n+1}}{(2n+1)!} \tag{168}$$

This includes the deformation through the non-communicative product. This means the region convergence may according to $\alpha$.

$$\|S\|_\alpha^2 = S_1^2 + S_2^2 + S_3^2 - 2(1-\alpha)S_2 S_3 \tag{169}$$

The Fourier Transform, let the unit $\cup$ derived from $T = yi + zi$ be considered.

$$\mathcal{F}_\alpha(f(t))(w) = \int_{-\infty}^\infty f(t)\exp(-wt\cup)\, dt \tag{170}$$

Our Fourier Transform may reveal how frequencies interact in space with intrinsic anisotropy. Finally, we present the transfer function,

$$G(S) = \mathcal{L}_\alpha(g(t))(S) \tag{171}$$

# Reimagining spacetime: Einstein's Legacy meets Trinition Geometry

Einstein's formulation of spacetime, where the line element in the Minkowski space was introduced as:

$$ds^2 = -cdt^2 + dx^2 + dy^2 + dz^2 \qquad \qquad \textbf{172}$$

has been a pillar of our understanding of gravity and the makeup of the universe for decades. The elegant form combines time and three-dimensional space into a four-dimensional manifold endowed with a Lorentzian metric that defines gravity. Yet, despite these accolades, Einstein's spacetime can sometimes be problematic when reconciling macroscopic experience with the quantum reality at play in the universe or when we want to describe delicate anisotropies in the universe. And here is where a fresh perspective can step in the possibility of extending/enriching the paradigm of Einstein with the new Trinition framework. Within the framework of Trinition, we define a line element event

$$E = t + Z = t + x + yi + zj \qquad (173)$$

Becomes

$$ds^2 = -cdt^2 + dx^2 + dy^2 + dz^2 - 2(1-\alpha)dydz \qquad (174)$$

The Trinition retake on Einstein's spacetime is not simply an extension; it is a daring generalization. It maintains the known four-dimensional architecture while replacing the metric with a deformable, variable one interpolating between the classical Euclidean and quaternionic cases. Although a mere starting point, the Trinition framework may have the potential to significantly impact the way we scope gravitational interactions and the underlying geometry of our indeed universe under noncommutativity, hidden degrees of freedom, and potentially anisotropic deformations.

This generalization could indeed be a breakthrough in our understanding of spacetime, a breakthrough that honors the legacy of Einstein and opens up entirely new avenues of research in quantum gravity, cosmology, etc. The obtained formula could significantly modify the explicit Einstein's field equations. For instance, we could interpret the metric components

$$g_{00} = -c^2, g_{11} = 1, g_{22} = 1, g_{33} = 1, \qquad g_{23} + g_{32} = -2(1-\alpha)$$

This could add some values to the black hole physics, quasinormal modes, gravitational wave emission, potential experimental signatures, and the connection and curvature. Because the author is not an expert in this field, he will leave expecting to investigate and decide

## Some remarks

We shall provide alternative formulas for the Triconus, Trisinus, and Trinus. Now consider a general Trinition number with $Z = x + iy + zj$ and it deformed norm defined by $r = \sqrt{x^2 + (y^2 + z^2 - 2(1-\alpha)yz)^2}$, here we assume that $\alpha \neq 0 \text{ and } 2$, we introduce the polar-angle $\varphi$ through

$$\cos(\varphi) = \frac{x}{r}, \quad \sin(\varphi) = \frac{\sqrt{y^2 + z^2 - 2(1-\alpha)yz}}{r}$$

In many applications, one finds that when a particular coordinate choice is made, for example,

$$x = r\cos(\varphi), \quad y = x = r\cos(\varphi), z = r\sin(\varphi)$$

So that we can obtain

$$z^2 + y^2 = r^2 \text{ and } 2(1-\alpha)yz = 2(1-\alpha)r^2\cos(\varphi)\sin(\varphi) = (1-\alpha)r^2\sin(2\varphi)$$

$$y^2 + z^2 - 2(1-\alpha)yz = r^2[1 - (1-\alpha)\sin(2\varphi)]$$

Then writing

$$B(\varphi) = 1 - (1-\alpha)\sin(2\varphi)$$

We can now redefine our Archiometry functions as:

$$Trinconus(\varphi) = \frac{\cos(\varphi)}{\sqrt{\cos(\varphi)^2 + (B(\varphi))^2}}$$

$$Trisinus(\varphi) = \frac{B(\varphi)}{\sqrt{\cos(\varphi)^2 + (B(\varphi))^2}}$$

$$Trinus(\varphi) = \frac{B(\varphi)}{\cos(\varphi)}$$

# Generalization and application to partial differential equation

Wave equations have historically been one of the cornerstones of mathematical physics, governing various natural phenomena, including sound propagation, electromagnetism, quantum mechanics, and fluid dynamics [23, 24. 25]. Classical wave equation are heavily based on Euclidean space and inner product structures [26]. Yet, recently to broaden the understanding of wave phenomena, there has been a series of alternative mathematical frameworks developed, yielding insights particularly in distorted and non-Euclidean regimes [27,28]. Trinition is a brand new hypercomplex system that is a natural extension of classical vector spaces, which leads to a deformable inner product structure and thereby modifies the core mathematical operations. However, in Trinition, we have an additional deformation parameter, α, which influences the interactions among basis elements. Such deformations affect not just the algebraic operations (in particular, the algebra of observables) but also the related differential structures entering wave equations, yielding new formulations of such equations. A fundamental problem in extending classical analysis to this new context is their understanding of the functional spaces that solutions to the Trinition wave equation would normally lie in. This requires a generalization of the Sobolev spaces, the spaces that are the natural arenas for studying solutions of PDE, to account for the new bilinear form. Certain classical inequalities like the Sobolev embedding theorem, and the Gagliardo-Nirenberg inequality, in particular, will have to be reformulated in the light of the Trinitarian perspective In this section, we elaborate systematically the Sobolev space theory in Trinition and demonstrate how standard inequalities look like with respect to the deformed inner product. We restate with new versions the Sobolev embedding and Gagliardo-Nirenberg inequalities [29] of the Trinition deformation factor. We then discuss applications to the wave equation: that these deformations preserve existence, uniqueness and regularity in Trinition spaces. We reveal their effects via numerical simulations, providing fresh understanding of wave propagation in non-Euclidean and hypercomplex environments. From Newton's second law, together with the conservation of energy in the Euclidean space, one can derive the classical wave equation. It assumes that there is a standard inner product that makes vector norms and gradients obey the familiar Pythagorean identity. But this assumption is not always appropriate in systems providing the significance to non-Euclidean effects, deformable metrics, or nonlinear inner products. Trinition reframes it as a tunable inner product that locally adjusts the geometric structure of space itself. Therefore, in Trinition, the wave equation will be given differently, due to the deformation of metric, will take the following general form:

$$\frac{\partial^2 w(t,X)}{\partial t^2} = \nabla_{T_\alpha} \cdot \left(D_{T_\alpha} \nabla_{T_\alpha} w(t,X)\right) + H(w) \tag{175}$$

Where $\nabla_{T_\alpha}$ is the gradient under the framework of the deformable Trinition metric, $D_{T_\alpha}$ is the deformation factor that will alters wave propagation characteristics, $H(w)$ is a nonlinear function. There are several advantages associated to this. First, we will incorporates geometric deformation, that is the wave speed and the propagation characteristics will now depends on the parameter $\alpha$, which allow the tunable deformations. This will indeed provide a framework for investigating waves in complex or anisotropic media. Secondly, we will obtain a generalized classical wave behaviour, thus, the classical wave equation emerges as a special case when the parameter $\alpha \to 1$. This indeed helps for a continuous transition from standard Euclidean wave equations to deformed wave equations. Thirdly, since, framework is deformable, indeed, the Sobolev embedding theorem should be revised, this will lead to a new functional space constraint. This will affect the energy estimates and the solution stability, potentially explaining waves phenomena in exotic materials. In addition to this, the new wave equation under the Trinition framework could be used in quantum field theory, nonlinear optics, and fluid mechanics. The deformation parameters could model dispersion effects, viscoelastic, or even nonlocal interaction. These are just a few reasons why, we believe, Trinition is a powerful framework with which one is able to re-examine classical wave mechanics, which in turn will yield an even deeper understanding of wave propagation in general. As we presented before

$$Z \in T_\alpha$$

if: $Z = x + yi + zj$, $i^2 = j^2 = -1$, $ij = 1 - \alpha + \alpha k$, $ji = 1 - \alpha - \alpha k$

Thus within the framework of the above

$$\|Z\|_\alpha = \sqrt{x^2 + y^2 + z^2 - 2(1-\alpha)yz}$$

The above norm suggest a deformation especially the co-relation of the $y - z$ space. This indeed is a step forward toward understanding our nature. We wish to define the Trinition basis and metric. We recall that in Euclidean space, we have the standard basis $e_x, e_y, e_z$ with dot products defining the metric. Now in the Trinition space, we will modify these basis vectors to encode the Trinition algebra. In the standard Euclidean space $R^3$, the basis vectors are:

$e_x = (1,0,0)$, $e_y = (0,1,0)$ and $e_z = (0,01)$. These satisfy the dot product relations as:

$e_x e_x = 1$, $e_y e_y = 1$, $e_z e_z = 1$, and $e_i e_j \neq 0$ $(i \neq j)$

But with the new defined Trinition, can we introduce a modified basis for example $\{E_x, E_y, E_z\}$ that will follow Trinition multiplication rules.

## The deformed Inner Product and Metric

Let us consider a Trinition vector $v = (x, y, z)$ the inner product was given by

$$< v, v >= x^2 + y^2 + z^2 - 2(1 - \alpha)yz. \tag{176}$$

The above quadratic form can be written in a matrix form as follow:

$$< v, v >= (x, y, z) g \begin{pmatrix} x \\ y \\ z \end{pmatrix} \tag{177}$$

Indeed, with the metric tensor

$$g_{ij} = \begin{pmatrix} 1 & 0 & 1 \\ 0 & 1 & -(1-\alpha) \\ 0 & -(1-\alpha) & 1 \end{pmatrix}. \tag{178}$$

By checking we have that

$$< v, v >= (x, y, z) g \begin{pmatrix} x \\ y \\ z \end{pmatrix} = x^2 + y^2 + z^2 - 2(1 - \alpha)yz \tag{179}$$

This lead us to the inverse metric and the Laplace-Beltrami operator [29]. Since the metric is constant because it is independent of the coordinates $x, y, z$, then we can compute the inverse metric $g^{ij}$ by inverting $g$. We have noticed that, $g$ is block-diagonal in $x$ and $(y, z)$ sector. We have then $g_{xx} = 1$ so that $g^{xx} = 1$. Now for the 2x2 remaining block

$$G_{yz} = \begin{bmatrix} 1 & -(1-\alpha) \\ -(1-\alpha) & 1 \end{bmatrix}. \tag{180}$$

The determinant of the above is given by:

$$Det(G_{yz}) = 1 - (1-\alpha)^2 = 1 - (1 - 2\alpha + \alpha^2) = 2\alpha - \alpha^2 = \alpha(2-\alpha). \tag{181}$$

Therefore, the inverse of $G_{yz}$ is:

$$G_{yz}^{-1} = \frac{1}{\alpha(2-\alpha)} \begin{bmatrix} 1 & (1-\alpha) \\ (1-\alpha) & 1 \end{bmatrix}. \tag{182}$$

Therefore, in full, the inverse metric is:

$$g^{ij} = \begin{pmatrix} 1 & 0 & 0 \\ 0 & \dfrac{1}{\alpha(2-\alpha)} & \dfrac{1-\alpha}{\alpha(2-\alpha)} \\ 0 & \dfrac{1-\alpha}{\alpha(2-\alpha)} & \dfrac{1}{\alpha(2-\alpha)} \end{pmatrix}. \tag{183}$$

But because the metric components are constant, we can determinant is :

$$|g| = \alpha(2-\alpha). \tag{184}$$

Because $\sqrt{|g|} = \sqrt{\alpha(2-\alpha)}$ is constant for a fixed alpha. We have

$$\Delta w = g^{ij} \partial_i \partial_j w. \tag{185}$$

By substituting the inverse metric components, we obtain

$$\Delta w = \partial_x^2 w + g^{yy}\partial_y^2 w + 2g^{yz}\partial_y \partial_z w + g^{zz}\partial_z^2 w. \tag{186}$$

With

$$g^{yy} = g^{zz} = \frac{1}{\alpha(2-\alpha)}, \quad g^{yz} = \frac{1-\alpha}{\alpha(2-\alpha)}. \tag{187}$$

So that,

$$\Delta w = \partial_x^2 w + \frac{1}{\alpha(2-\alpha)}\partial_y^2 w + 2\frac{1-\alpha}{\alpha(2-\alpha)}g^{yz}\partial_y\partial_z w + \frac{1}{\alpha(2-\alpha)}\partial_z^2 w. \tag{188}$$

Now we have to verify that $g^{ik}g_{kj} = \delta^i_j$. We will then compute the product $g^{ik}g_{kj}$, which should then result in the identity matrix $\delta^i_j$. Let us start with $x = i$, checking $g^{ik}g_{kj}$: But we have that $g_{xj} = (1,0,0)$ and $g^{ix} = (1,0,0)$, we get:

$$\sum_k g^{xk}g_{kj} = g^{xx}g_{yj} + g^{xy}g_{xj} + g^{xz}g_{zj} = 1.1 + 0.0 + 0.0 = 1. \tag{189}$$

Hence, the first row $g^{ik}g_{kj}$ is correct. $y = i$, checking $g^{ik}g_{kj}$, then computing:

$$\sum_k g^{yk}g_{kj} = g^{xx}g_{yj} + g^{xy}g_{zj}. \tag{190}$$

Now by expanding for each $j$: For $j = y$

$$g^{yy}g_{yy} + g^{yz}g_{zy} = \frac{1-\alpha}{\alpha(2-\alpha)}(1) + \frac{(1-\alpha)}{\alpha(2-\alpha)}(-(1-\alpha)) = \frac{2\alpha - \alpha^2}{\alpha(2-\alpha)} = 1. \tag{191}$$

For $j = z$

$$g^{yy}g_{yz} + g^{yz}g_{zz} = \frac{1-\alpha}{\alpha(2-\alpha)}(1) - \frac{1-\alpha}{\alpha(2-\alpha)} = 0. \tag{192}$$

$z = i$, checking $g^{ik}g_{kj}$, then computing:

$$\sum_k g^{zk}g_{kj} = g^{zy}g_{yj} + g^{zz}g_{zj}. \tag{193}$$

For $j = y$

$$g^{zy}g_{yy} + g^{zz}g_{zy} = \frac{1-\alpha}{\alpha(2-\alpha)}(1) - \frac{1-\alpha}{\alpha(2-\alpha)} = 0. \tag{194}$$

For $j = z$

$$g^{zy}g_{yz} + g^{zz}g_{zz} = \frac{1-\alpha}{\alpha(2-\alpha)}(-(1-\alpha)) + \frac{1}{\alpha(2-\alpha)}(1) = 1. \tag{195}$$

We can then conclude that $g^{ik}g_{kj} = \delta^i_j$, therefore, the given inverse metric components

$$g^{yy} = g^{zz} = \frac{1}{\alpha(2-\alpha)}, \quad g^{yz} = \frac{1-\alpha}{\alpha(2-\alpha)}.$$

Now we can reconstruct all the partial differential equations using the Laplacian, we present here the wave equation within the Trinition framework.

## Wave equation within the Trinition framework

Now by replacing the standard Laplacian in the classical wave equation we have

$$\frac{\partial^2 u}{\partial t^2} = c^2 \nabla^2 u$$

with the Laplace-Beltrami operator we derived above. This provides the Trinition wave equation:

$$\frac{\partial^2 u}{\partial t^2} = c^2 \left\{ \partial_x^2 u + \frac{1}{\alpha(2-\alpha)} \partial_y^2 w + 2 \frac{1-\alpha}{\alpha(2-\alpha)} \partial_y \partial_z u + \frac{1}{\alpha(2-\alpha)} \partial_z^2 u \right\}. \tag{196}$$

The next question that arise naturally is that, does the solution exist and is unique and stable? This should now lead us to rigorously establish existence, uniqueness and the stability, we will need to redefine Sobolev spaces in Trinition framework and also prove that, our deformed inner product still allows a well-posed functional analysis. We present below some numerical simulation of the deformable wave equation for different values of alpha.

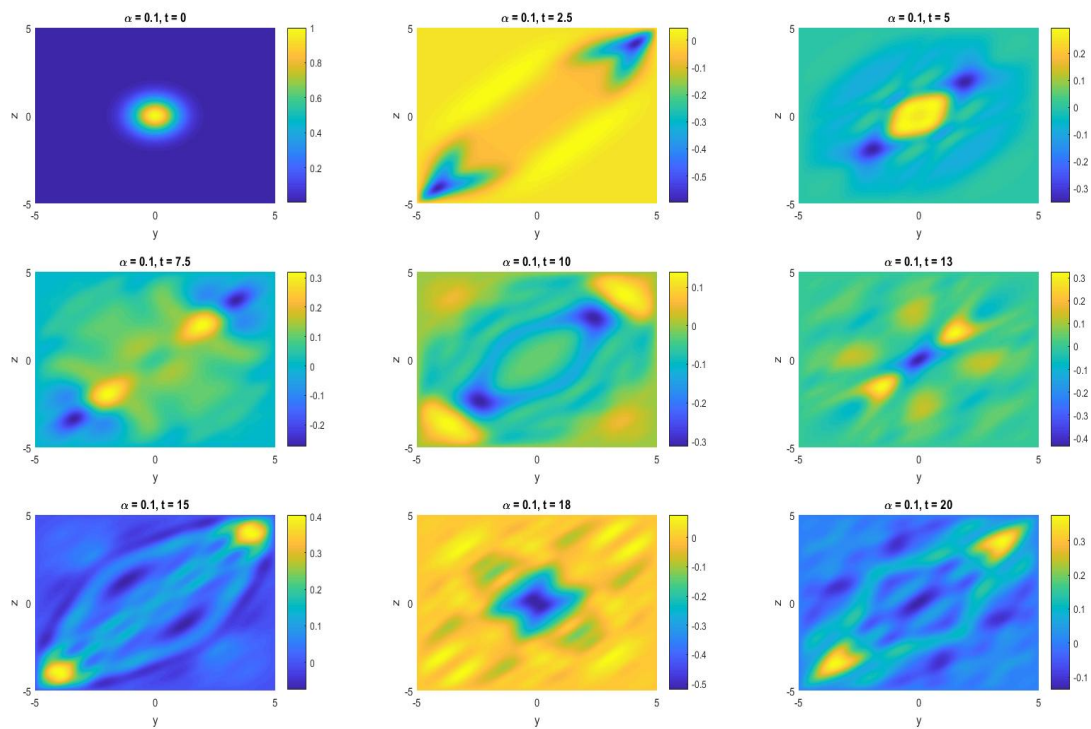

**Figure 31**: Simulation of the two dimensional wave equation for alpha 0.1

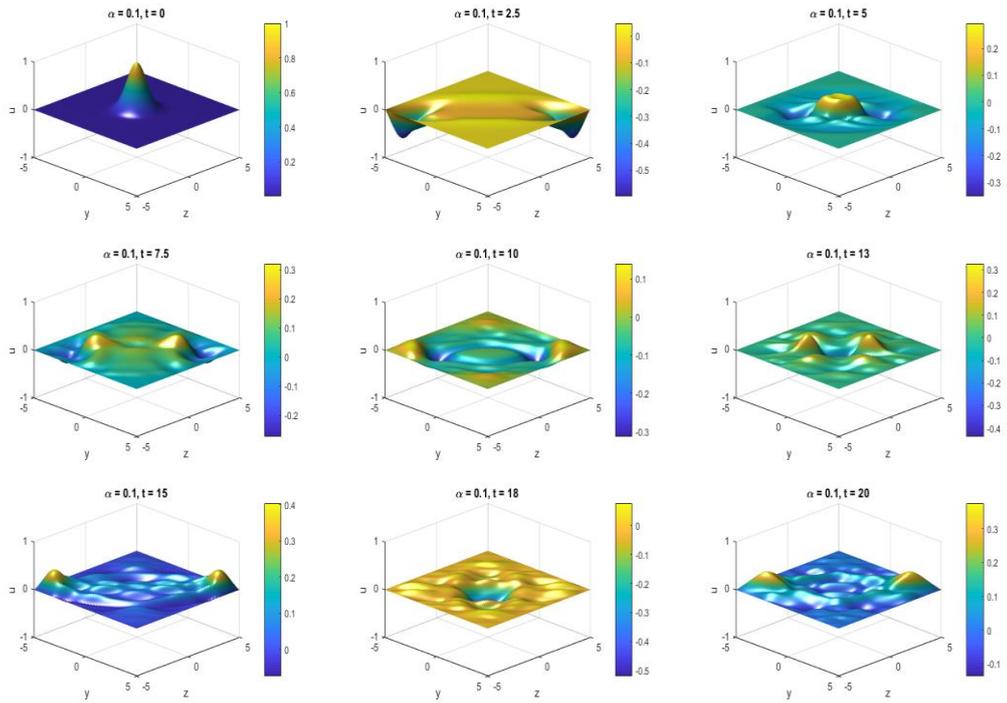

**Figure 32:** Simulation of the three dimensional wave equation for alpha 0.1

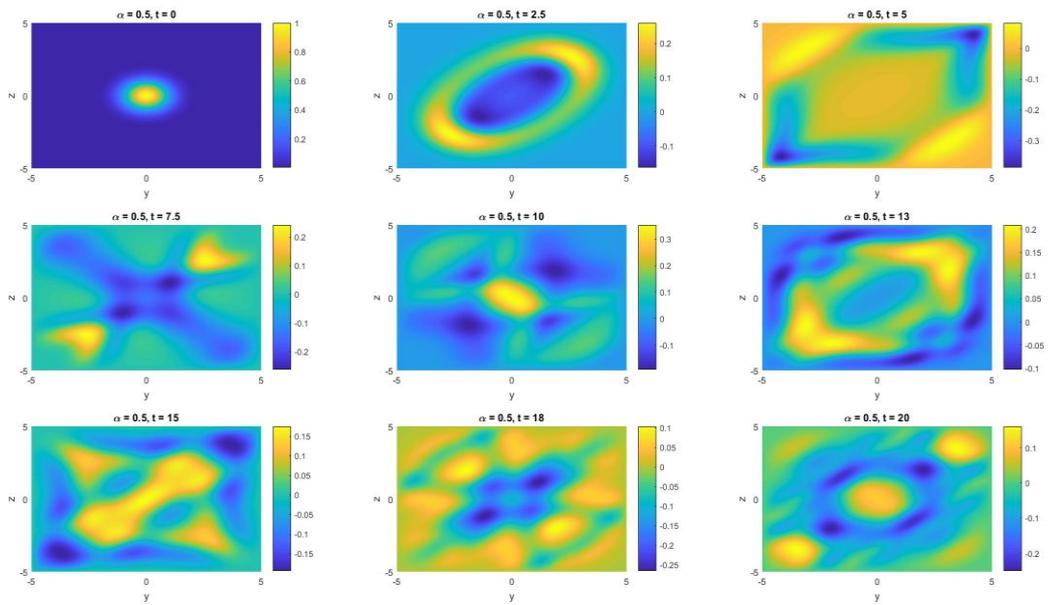

**Figure 33**: Simulation of the two dimensional wave equation for alpha 0.5

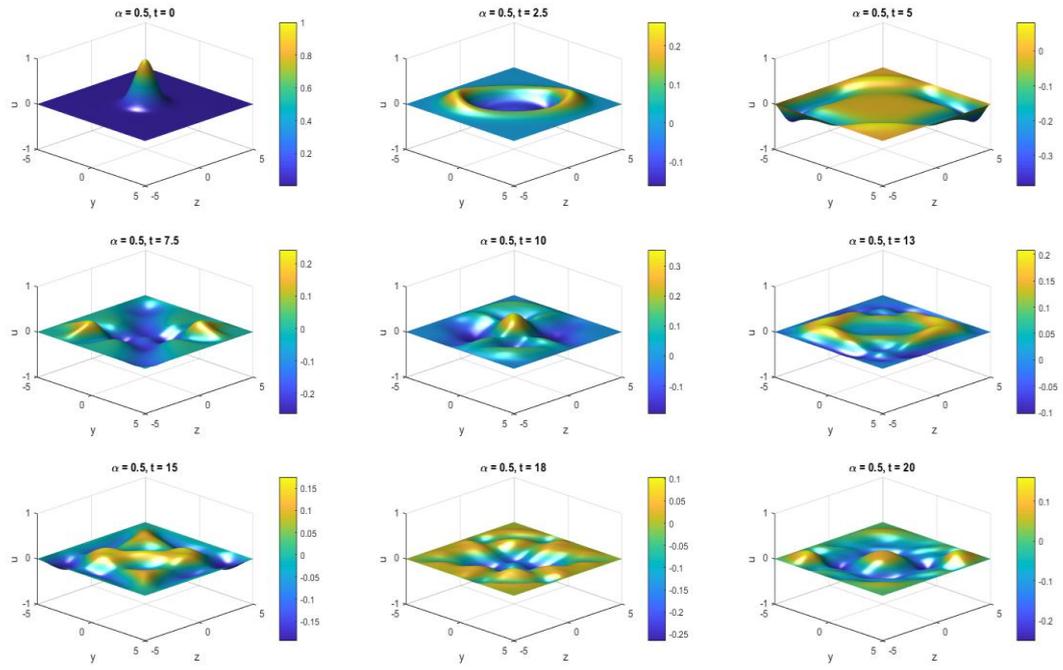

**Figure 34:** Simulation of the three dimensional wave equation for alpha 0.5

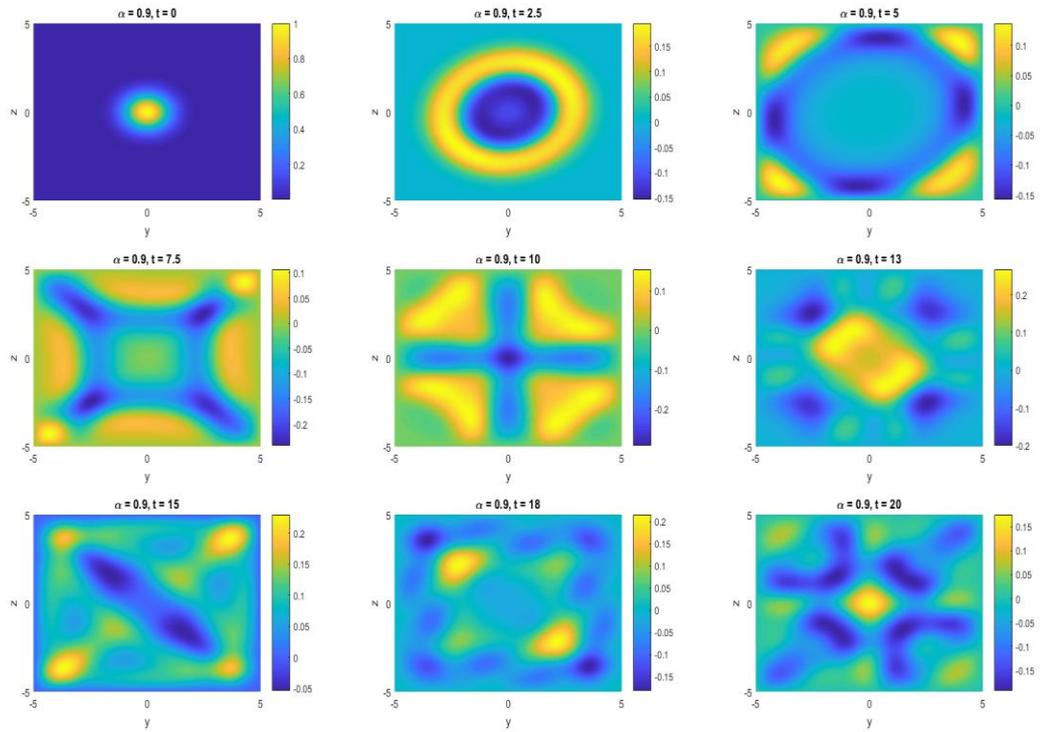

**Figure 35**: Simulation of the two dimensional wave equation for alpha 0.9

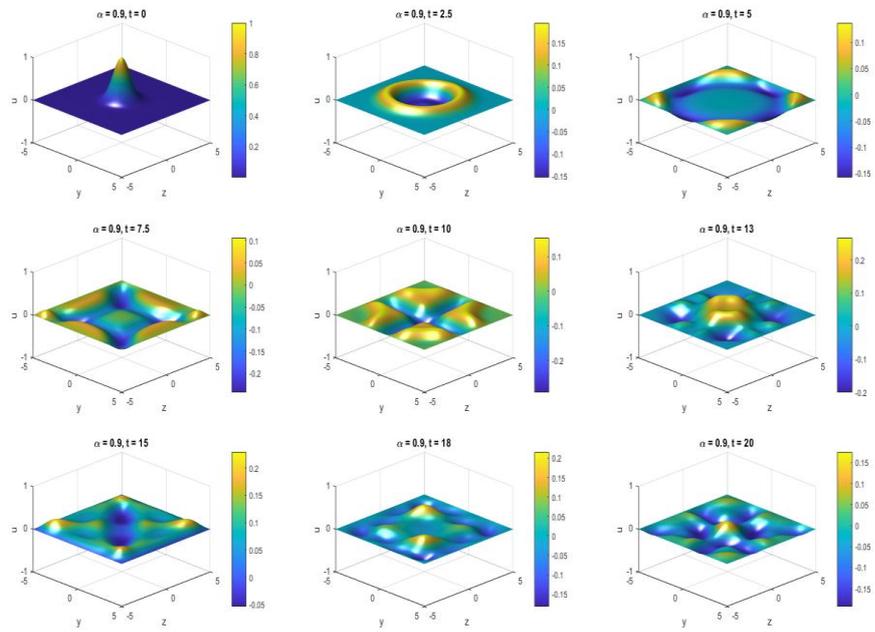

**Figure 36:** Simulation of the three dimensional wave equation for alpha 0.9

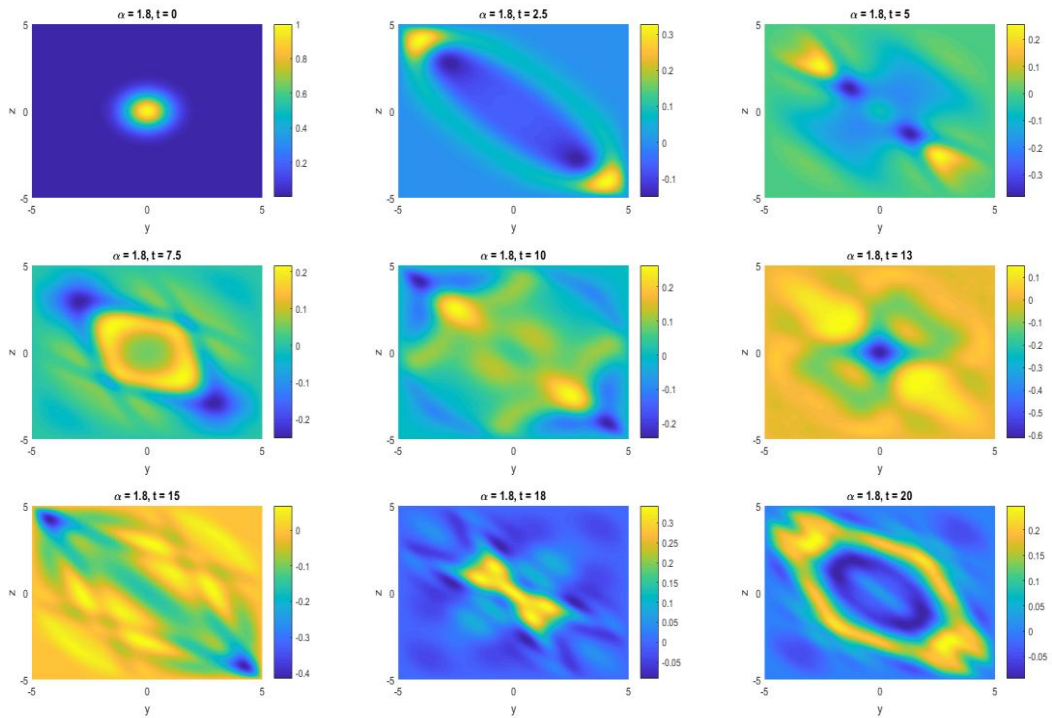

**Figure 37**: Simulation of the two dimensional wave equation for alpha 1.8

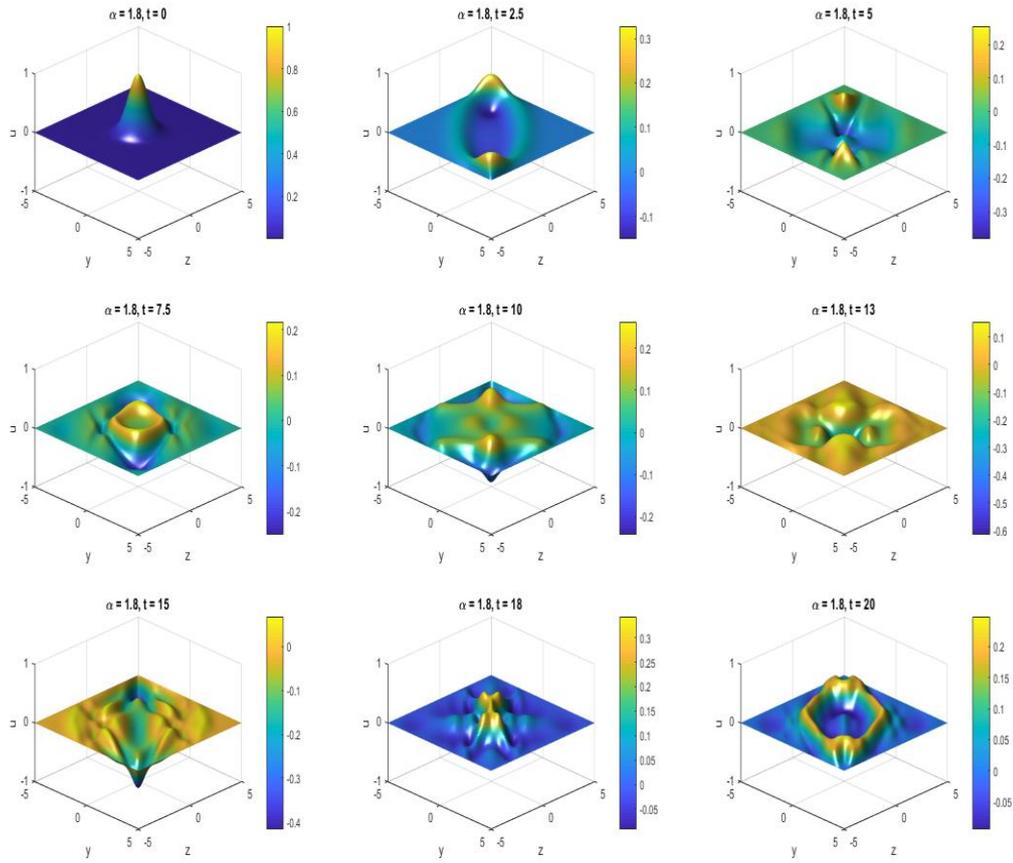

**Figure 38**:Simulation of the three dimensional wave equation for alpha 1.8

# Trinition modified Sobolev Spaces

With what we presented above, we have that the inner product is now deformed by the metric components $g^{yy} = g^{zz} = \frac{1}{\alpha(2-\alpha)}$, $g^{yz} = \frac{1-\alpha}{\alpha(2-\alpha)}$, we need to define a new space $H^1_{Trinition}(\Omega)$. We suggest the following definition:

$$H^1_{Trinition}(\Omega) = \{u \in L^2(\Omega)| \nabla_{Trinition} u \in L^2(\Omega)\} \quad (197)$$

Where in our case the deformed inner product is:

$$<u,v>_\alpha = \int_\Omega \left\{uv + g^{xx}\frac{\partial u}{\partial x}\frac{\partial v}{\partial x} + g^{yy}\frac{\partial u}{\partial y}\frac{\partial v}{\partial y} + g^{zz}\frac{\partial u}{\partial z}\frac{\partial v}{\partial z} + 2g^{yz}\frac{\partial u}{\partial y}\frac{\partial v}{\partial z}\right\} d\Omega. \quad (198)$$

The corresponding norm will then be given as:

$$\|u\|_{H^1_{Trinition}} = \left\{ \|u\|^2_{L^2} + g^{xx}\left\|\frac{\partial u}{\partial x}\right\|^2_{L^2} + g^{yy}\left\|\frac{\partial u}{\partial y}\right\|^2_{L^2} + g^{zz}\left\|\frac{\partial u}{\partial z}\right\|^2_{L^2} \right. \tag{199}$$

$$\left. + 2g^{yz}\int_\Omega \frac{\partial u}{\partial y}\frac{\partial v}{\partial z}d\Omega \right\}^{1/2}.$$

The above space is complete; therefore, it is a Hilbert space.

## Weak formulation and well-posedness under the Trinition framework

To do this, we will consider the Trinition wave –modified equation under the following form:

$$\frac{\partial^2 u}{\partial t^2} - g^{xx}\frac{\partial^2 u}{\partial x^2} - g^{yy}\frac{\partial^2 u}{\partial y^2} - g^{zz}\frac{\partial^2 u}{\partial z^2} - 2g^{yz}\frac{\partial^2 u}{\partial y\partial z} = 0 \tag{200}$$

For the weak formulation, we will multiply by a test function $v \in H^1_{Trinition}(\Omega)$ and then integrating over $\Omega$ we obtain the following:

$$\int_\Omega \frac{\partial^2 u}{\partial t^2}vd\Omega - \int_\Omega g^{xx}\frac{\partial^2 u}{\partial x^2}vd\Omega - \int_\Omega g^{xx}\frac{\partial^2 u}{\partial y^2}g^{yy}vd\Omega - \int_\Omega g^{zz}\frac{\partial^2 u}{\partial x^2}vd\Omega \tag{201}$$

$$- \int_\Omega 2g^{yz}\frac{\partial^2 u}{\partial y\partial z}vd\Omega = 0$$

We will present the integration by part of one component, which will be the same for others.

$$\int_\Omega g^{xx}\frac{\partial^2 u}{\partial x^2}vd\Omega = \int_\Omega g^{xx}\frac{\partial u}{\partial x}\frac{\partial v}{\partial x}d\Omega - \int_{\partial\Omega} g^{xx}\frac{\partial u}{\partial x}vdS \tag{202}$$

If we do not have Neumann boundary conditions, therefore

$$\int_\Omega g^{xx}\frac{\partial^2 u}{\partial x^2}vd\Omega = \int_\Omega g^{xx}\frac{\partial u}{\partial x}\frac{\partial v}{\partial x}d\Omega \tag{203}$$

Therefore after integration by part, we obtain:

$$\int_\Omega \frac{\partial^2 u}{\partial t^2} v d\Omega + \int_\Omega g^{xx} \frac{\partial u}{\partial x}\frac{\partial v}{\partial x} d\Omega + \int_\Omega g^{yy} \frac{\partial u}{\partial y}\frac{\partial v}{\partial y} d\Omega + \int_\Omega g^{zz} \frac{\partial u}{\partial z}\frac{\partial v}{\partial z} d\Omega \qquad (204)$$

$$+ 2\int_\Omega g^{yz}\left\{\frac{\partial u}{\partial y}\frac{\partial v}{\partial z} + \frac{\partial u}{\partial z}\frac{\partial v}{\partial y}\right\} d\Omega$$

$$= \int_{\partial\Omega} \left\{g^{xx}\frac{\partial u}{\partial x} + g^{yy}\frac{\partial u}{\partial y} + g^{zz}\frac{\partial u}{\partial z} + 2g^{yz}\frac{\partial u}{\partial z}\right\} v dS$$

If we do not use the Neumann, the above becomes

$$\int_\Omega \frac{\partial^2 u}{\partial t^2} v d\Omega + \int_\Omega g^{xx} \frac{\partial u}{\partial x}\frac{\partial v}{\partial x} d\Omega + \int_\Omega g^{yy} \frac{\partial u}{\partial y}\frac{\partial v}{\partial y} d\Omega + \int_\Omega g^{zz} \frac{\partial u}{\partial z}\frac{\partial v}{\partial z} d\Omega \qquad (205)$$

$$+ 2\int_\Omega g^{yz}\left\{\frac{\partial u}{\partial y}\frac{\partial v}{\partial z} + \frac{\partial u}{\partial z}\frac{\partial v}{\partial y}\right\} d\Omega = 0$$

This defines the bilinear form within the Trinition as:

$$B(u,v) = \int_\Omega g^{xx} \frac{\partial u}{\partial x}\frac{\partial v}{\partial x} d\Omega + \int_\Omega g^{yy} \frac{\partial u}{\partial y}\frac{\partial v}{\partial y} d\Omega + \int_\Omega g^{zz} \frac{\partial u}{\partial z}\frac{\partial v}{\partial z} d\Omega \qquad (206)$$

$$+ 2\int_\Omega g^{yz}\left\{\frac{\partial u}{\partial y}\frac{\partial v}{\partial z} + \frac{\partial u}{\partial z}\frac{\partial v}{\partial y}\right\} d\Omega$$

Therefore, for existence, we check the coercivity:

$$B(u,u) = \int_\Omega g^{xx} \left|\frac{\partial u}{\partial x}\right|^2 d\Omega + \int_\Omega g^{yy} \left|\frac{\partial u}{\partial y}\right|^2 d\Omega + \int_\Omega g^{zz} \left|\frac{\partial u}{\partial z}\right|^2 d\Omega \qquad (207)$$

$$+ 2\int_\Omega g^{yz}\left\{\frac{\partial u}{\partial y}\frac{\partial u}{\partial z}\right\} d\Omega$$

If $g^{xx}, g^{yy} > 0$, then $B(u,u) \geq C\|u\|^2_{H^1_{Trinition}}$, this ensures existence through the Lax-Milgram theorem. We now define the energy function as:

$$E(t) = \frac{1}{2}\int_\Omega \left\{\left|\frac{\partial^2 u}{\partial t^2}\right|^2 + g^{xx}\left|\frac{\partial u}{\partial x}\right|^2 + g^{yy}\left|\frac{\partial u}{\partial y}\right|^2 + g^{zz}\left|\frac{\partial u}{\partial z}\right|^2 + 2g^{zy}\frac{\partial u}{\partial y}\frac{\partial u}{\partial z}\right\}d\Omega \tag{208}$$

Taking the time derivative:

$$\frac{dE(t)}{dt} = \frac{1}{2}\int_\Omega \frac{\partial}{\partial t}\left\{\left|\frac{\partial^2 u}{\partial t^2}\right|^2 + g^{xx}\left|\frac{\partial u}{\partial x}\right|^2 + g^{yy}\left|\frac{\partial u}{\partial y}\right|^2 + g^{zz}\left|\frac{\partial u}{\partial z}\right|^2 + 2g^{zy}\frac{\partial u}{\partial y}\frac{\partial u}{\partial z}\right\}d\Omega \tag{209}$$
$$= 0$$

So $E(t) = E(0)$ implying stability. We will now generalize some well-known inequalities.

## Generalizing the Sobolev and Garliardo-Nirenberg inequalities in the Trinition framework

When we move to a deformed setting like our Trinition framework we refined the gradient norm using the deformed metric. For example, if the inner product is given by the metric tensor $g^{ij}$, then the deformed gradient norm will be:

$$\|\nabla_T u\|_{L^2}^2 = \int_\Omega \sum_{i,j} g^{i,j}\partial_i u \partial_j u \, dx \tag{210}$$

We note that this is an integral of a quadratic form in the derivatives. It is not of course equivalent to summing the product of individual $L^2$ norms of the derivatives. In order words, in our case the correct way to incorporate the deformation is to defined the Sobolev norms as:

$$\|u\|_{H_T^1(\Omega)}^2 = \int_\Omega \left(u^2 + \sum_{i,j} g^{i,j}\partial_i u \partial_j u\right) dx \tag{211}$$

In this case the deformed Sobolev embedding would be

$$\|u\|_{L^q(\Omega)} \leq C\|u\|_{H_T^1(\Omega)} \tag{212}$$

We recall that the expression

$$\left(\sum_{i,j} g^{ij}\|\partial_i u\|_{L^2}\|\partial_j u\|_{L^2}\right)^{1/2}$$

is not the same as the $L^2-$norm of the gradient computed with the metric $g^{ij}$. The proper formulation in the Trinition setting to work with the deformed Sobolev norm dined by inner product

$$\|u\|_{H^1_T(\Omega)} = \int_\Omega \left(uv + \sum_{i,j} g^{i,j}\partial_i u\partial_j u\right)dx \tag{213}$$

And then, the corresponding embedding inequality becomes

$$\|u\|_{L^q(\Omega)} \leq C\left\{\int_\Omega \left(u^2 + \sum_{i,j} g^{i,j}\partial_i u\partial_j u\right)dx\right\}^{1/2} \tag{214}$$

This is how the deformation properly enters within the Sobolev space formulation, and also guarantees that the deformation factors (which depends on the metric) enter through the integral rather than as factor products of norms.

Gagliardo–Nirenberg inequality

$$\|\nabla_T u\|_{L^r} = \left(\int_\Omega \left\{\sum_{i,j} g^{i,j}(x)\partial_i u(x)\partial_j u(x)\right\}^{\frac{r}{2}} dx\right)^{\frac{1}{r}} \tag{215}$$

And then the inequality is defined as

$$\|u\|_{L^q(\Omega)} \leq C_T \|\nabla_T u\|_{L^r(\Omega)}^\theta \|u\|_{L^s(\Omega)}^{1-\theta} \tag{216}$$

# Generalization with 3 deformable parameters

In this section, the deformation is introduced also in $x-y$ and $x-z$. Consider a Trinition number $Z = x + yi + zj$. We defined the following norm

$$\|Z\|_{\alpha,\gamma,\beta} = \sqrt{x^2 + y^2 + z^2 - 2(1-\alpha)xy - 2(1-\gamma)xz - 2(1-\beta)yz} \tag{217}$$

We want to express it as:

$$\|Z\|_{\alpha,\gamma,\beta}^2 = (x,y,z)G\begin{pmatrix}x\\y\\z\end{pmatrix} \tag{218}$$

With the aim that, the metric tensor $G = (g_{ij})$ is determined by matching the respective coefficients. From the norm, we have that $g_{11} = g_{22} = g_{33} = 1$. The coefficient of $xy$ is $-2(1-\alpha)xy$ in quadratic form, the cross term $xy$ do appear as $2g^{12}xy$. Therefore,

$$-2(1-\alpha) = 2g_{12} \to g_{12} = g_{21} = -(1-\alpha)$$

Similarly, we will have

$$-2(1-\gamma) = 2g_{13} \to g_{13} = g_{31} = -(1-\gamma)$$

$$-2(1-\beta) = 2g_{23} \to g_{23} = g_{32} = -(1-\beta)$$

Therefore, the metric tensor in matrix form is given from the above as:

$$G = \begin{pmatrix} 1 & -(1-\alpha) & -(1-\gamma) \\ -(1-\alpha) & 1 & -(1-\beta) \\ -(1-\gamma) & -(1-\beta) & 1 \end{pmatrix}. \tag{219}$$

We will now compute the inverse Metric $g^{23}$. Thus, we wish to find the inverse matrix $G^{-1}$ that will verify that $G^{-1}G = I$. For simplicity, we will let

$$a = 1-\alpha, \quad b = 1-\gamma, \quad c = 1-\beta$$

Such that the above metric becomes

$$G = \begin{pmatrix} 1 & -a & -c \\ -a & 1 & -b \\ -c & -b & 1 \end{pmatrix}. \tag{220}$$

After calculation, we found that the determinant of the deformed metric is given by:

$$\det(G) = 1 - 2abc - (a^2 + b^2 + c^2). \tag{221}$$

We know that the inverse is obtained with the following formula

$$G^{-1} = \frac{adj(G)}{det(G)}$$

Where the $adj(G)$ is the adjugate matrix, in our case it is given by:

$$adj(G) = \begin{pmatrix} 1 - b^2 & a + bc & ab + c \\ a + bc & 1 - c^2 & b + ac \\ ab + c & b + ac & 1 - a^2 \end{pmatrix}. \tag{222}$$

Therefore, we have the inverse as:

$$G^{-1} = \frac{1}{(1 - 2abc - (a^2 + b^2 + c^2))} \begin{pmatrix} 1 - b^2 & a + bc & ab + c \\ a + bc & 1 - c^2 & b + ac \\ ab + c & b + ac & 1 - a^2 \end{pmatrix}. \tag{223}$$

Now we shall provide some interpretation and also verify that the obtained coefficients really meet the requirements.

$$\sum_l g^{il} g_{lj} = \delta^i_j. \tag{224}$$

Noting that, the above expression for $g^{ij}$ was derived to hold the condition, and one can verify it by just multiplying $g^{ij}$ and $g_{jl}$ explicitly. For instance for the (224) component, we show that

$$\sum_l g^{1l} g_{l1} = g^{11} g_{11} + g^{12} g_{21} + g^{13} g_{31}. \tag{225}$$

Let us make substitution:

$$g^{11} = \frac{1 - b^2}{\det(G)}, \quad g^{12} = \frac{a + bc}{\det(G)}, \quad g^{13} = \frac{ab + c}{\det(G)}, g_{11} = 1, g_{21} = -a, g_{31} = -c$$

$$g^{11} g_{11} = \frac{1 - b^2}{\det(G)}, \quad g^{12} g_{21} = \frac{a + bc}{\det(G)}(-a), \quad g^{13} g_{31} = \frac{ab + c}{\det(G)}(-c)$$

$$\frac{1-b^2}{\det(G)} + \frac{a+bc}{\det(G)}(-a) + \frac{ab+c}{\det(G)}(-c) = \frac{1-2abc-(a^2+b^2+c^2)}{\det(G)} = 1. \qquad (226)$$

## Geometrical figures with three deformable parameters

Here is a complete description of how the standard geometric figures are re-formulated when distances are taken using the deformed Trinition norm. Within this framework, the squared "distance" (or norm) of a point $Z = (x, y, z)$

$$\|Z\|_{\alpha,\gamma,\beta}^2 = x^2 + y^2 + z^2 - 2(1-\alpha)xy - 2(1-\gamma)xz - 2(1-\beta)yz$$

(227)

This way, the replaced distance gives rise to the modified distance, which will result in the equation of the geometric figure being changed for other figures defined in terms of distances. In the sequel we derive the formulas of some classical figures in this new framework. The distance between any two points A and B is given by:

$$\begin{aligned} d(A-B) &= \|A-B\|_{\alpha,\gamma,\beta}^2 \\ &= (x-a)^2 + (y-b)^2 + (z-c)^2 - 2(1-\alpha)(x-a)(y-b) \\ &\quad - 2(1-\gamma)(x-a)(z-c) - 2(1-\beta)(y-b)(z-c) \end{aligned} \qquad (228)$$

Therefore the generalized sphere in 3D with center B and radius R is defined as the set of points $P(x, y, z)$ that satisfy:

$$\begin{aligned} d(P-B) &= \|A-B\|_{\alpha,\gamma,\beta}^2 \\ &= (x-a)^2 + (y-b)^2 + (z-c)^2 - 2(1-\alpha)(x-a)(y-b) \\ &\quad - 2(1-\gamma)(x-a)(z-c) - 2(1-\beta)(y-b)(z-c) = R^2 \end{aligned} \qquad (229)$$

Note from the above formula that, Due to the cross-terms, the figure does not have to be a sphere in the classical sense it can be an ellipsoid or a more general quadratic surface depending on the deformation parameters. Noting that in the above the radius $R$ is fixed but in general, it can be also deformed using the same norm, this will give birth to a very strange object especially as the deformed parameters change.

Let us look at a circle in a 2D subspace, let us consider the $2D$ $xy$-plane or any plane defined by two coordinates with a deformed norm. Let us consider for example $xy$-plane with $z = 0$, the deformed norm will become

$$\|Z\|_\alpha^2 = x^2 + y^2 - 2(1 - \alpha)xy \tag{230}$$

Therefore, a circle centred at the point $A(a, b)$ with a non-deformed radius R will be given by

$$R^2 = (x - a)^2 + (y - b)^2 - 2(1 - \alpha)(x - a)(y - b) \tag{231}$$

But while this has a familiar appearance of a classical circle equation, the terms in cross means that the points satisfying this equation will, in general, be rotated or sheared with respect to the standard circle actually defining an ellipse (or a tilted ellipse) in the classical Euclidean sense. Ellipse and general quadratic surfaces, in a broader context, we can express any quadratic surface in the deformed space in matrix form. Relative to some center $P_0$, define the coordinate vector as:

$$X = \begin{pmatrix} x - a \\ y - b \\ z - c \end{pmatrix} \tag{232}$$

Then, using the defined deformed norm, we have:

$$d(Z - A) = \|Z - A\|_{\alpha,\gamma,\beta}^2 = X^T G X \tag{233}$$

Where the metric tensor G associated with the deformed inner product is given by:

$$G = \begin{pmatrix} 1 & -(1 - \alpha) & -(1 - \gamma) \\ -(1 - \alpha) & 1 & -(1 - \beta) \\ -(1 - \gamma) & -(1 - \beta) & 1 \end{pmatrix}. \tag{234}$$

Thus, a more general quadratic surface that could be a sphere, ellipsoid, hyperboloid and other geometrical figure that can be obtained within the framework of the Trinition deformed geometry is given by:

$$R^2 = X^T G X \tag{235}$$

Noting that with the precise nature of the surface determined by the eigenvalues and eigenvector of the matric tensor $G$. Such shows that all classical geometry figures using the distance concept altered by application of the Trinition norm are simply converted to deformed versions of themselves. Their shapes, sizes, orientations are now functions of the deformation parameters, introducing a new rich structure in geometry. We present below some geometrical structure within the extended Trinition in Figure 39 and 40

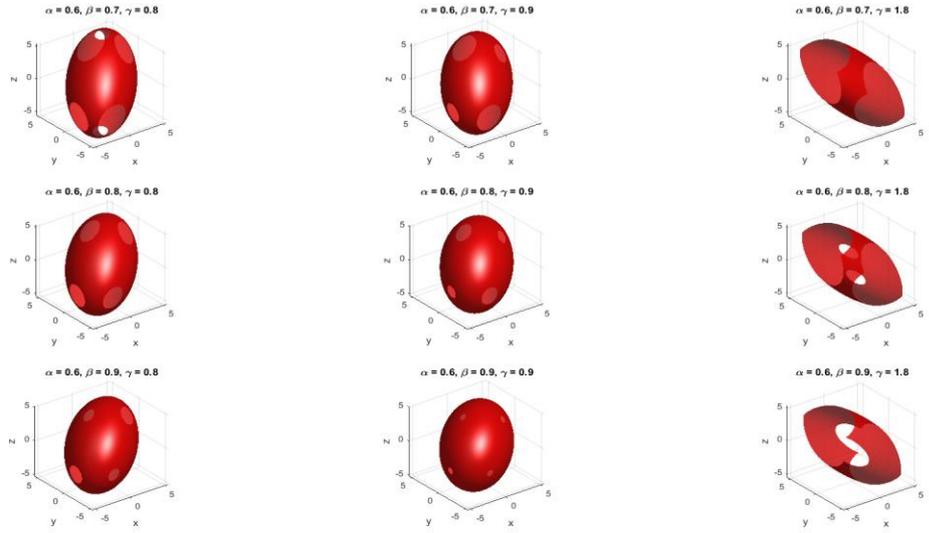

**Figure 39**: A Trinition sphere

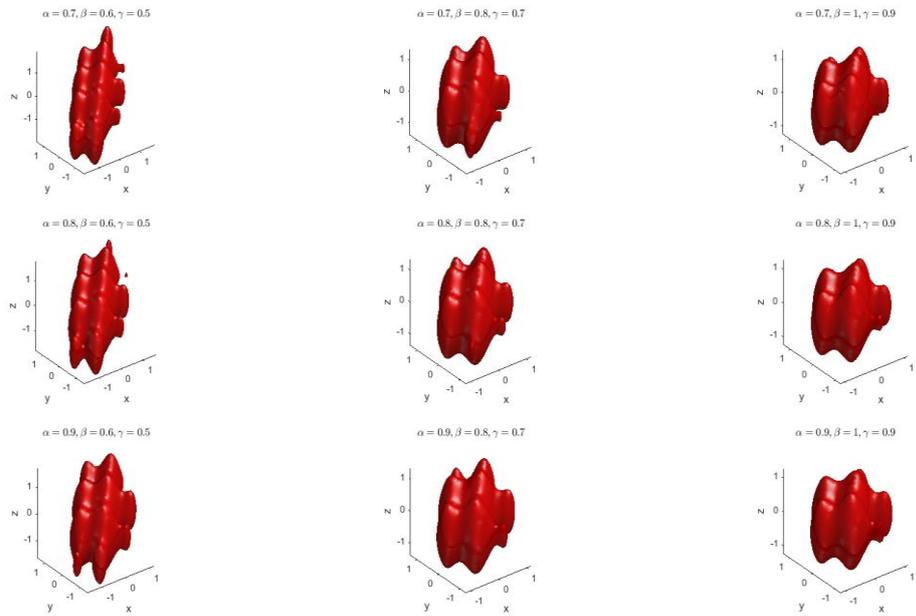

**Figure 40**: A Triniboid

# A generalized wave equation with 3 deformable parameters

We shall present a wave equation under this framework. Following the derivation presented earlier, in a curved or deformed spatial geometry, the Laplace-Beltrami operator acting on a scalar function $u(x,y,z,t)$ is

$$\Delta_G u(x,y,z,t) = \sum_{i,j} \frac{1}{\sqrt{|g|}} \partial_i(\sqrt{g} g^{ij} \partial_j u). \tag{236}$$

We recall that the determinant of the metric tensor is:

$$|g| = 1 - (1-\alpha)^2 - (1-\beta)^2 - (1-\gamma)^2 - 2(1-\alpha)(1-\beta)(1-\gamma) \tag{237}$$

And the inverse of the metric tensor $g^{ij}$ that satisfies $g^{ik} g_{kj} = \delta^i_j$ and was computed as:

$$g^{ij} = \frac{1}{(1 - 2abc - (a^2 + b^2 + c^2))} \begin{pmatrix} 1-b^2 & a+bc & ab+c \\ a+bc & 1-c^2 & b+ac \\ ab+c & b+ac & 1-a^2 \end{pmatrix}. \tag{238}$$

Here the indices $i,j$ run over the spatial coordinates $x, y\ z$ then the Trinition-modified wave becomes:

$$\frac{\partial^2 u}{\partial t^2} = c^2 \sum_{i,j} \frac{1}{\sqrt{|g|}} \partial_i(\sqrt{g} g^{ij} \partial_j u) \tag{239}$$

The numerical simulation are presented in Figures 41, 42, 43 and 44.

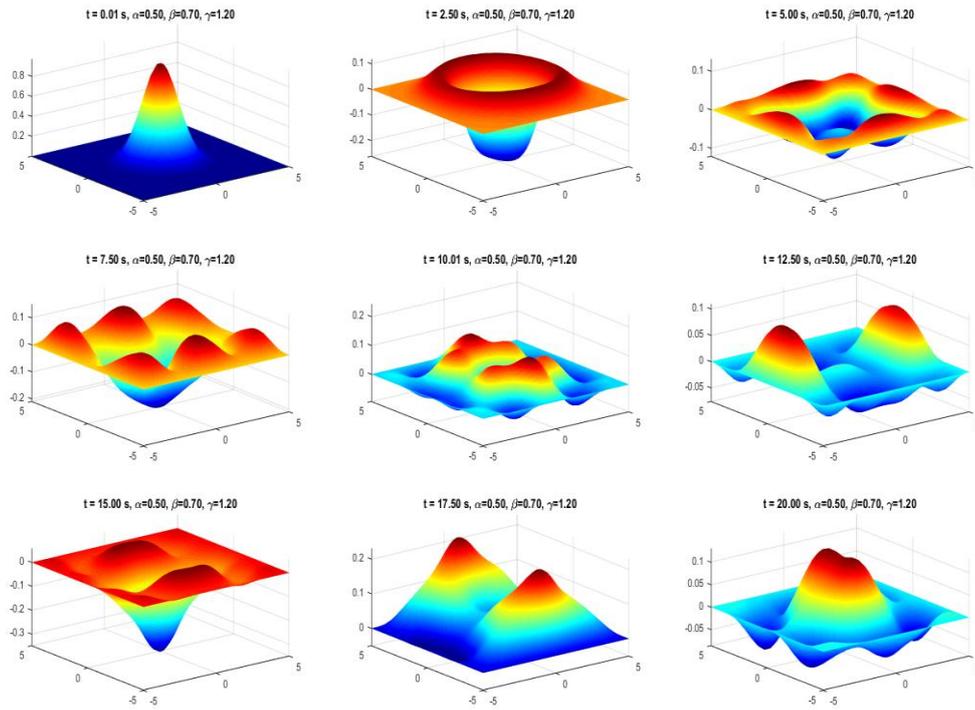

**Figure 41**: Numerical simulation of a Trinition modified wave equation for (0.50, 0.70, 1.20)

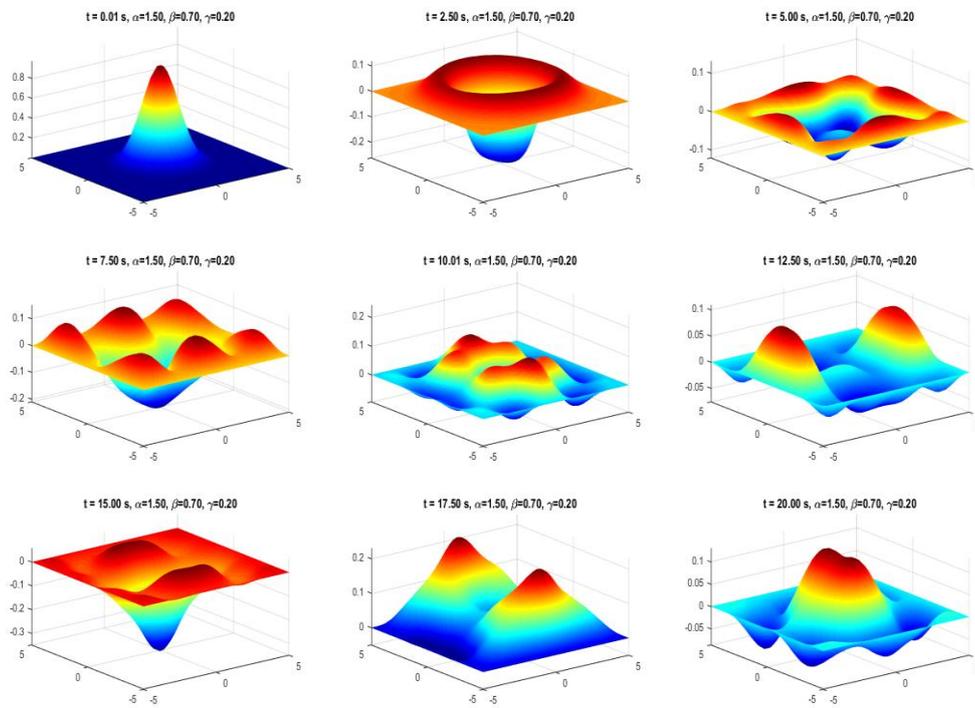

**Figure 42**: Numerical simulation of a Trinition modified wave equation for (1.50, 0.70, 0.20)

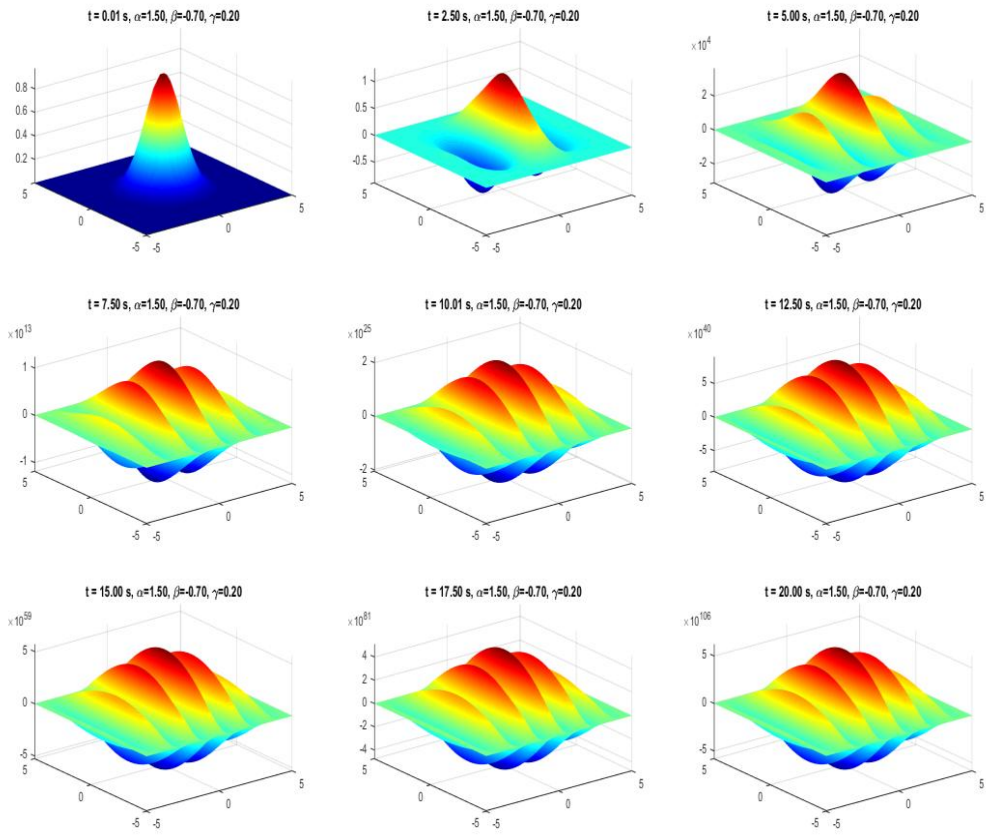

**Figure 43**: Numerical simulation of a Trinition modified wave equation for (1.50, 0.70, 0.20)

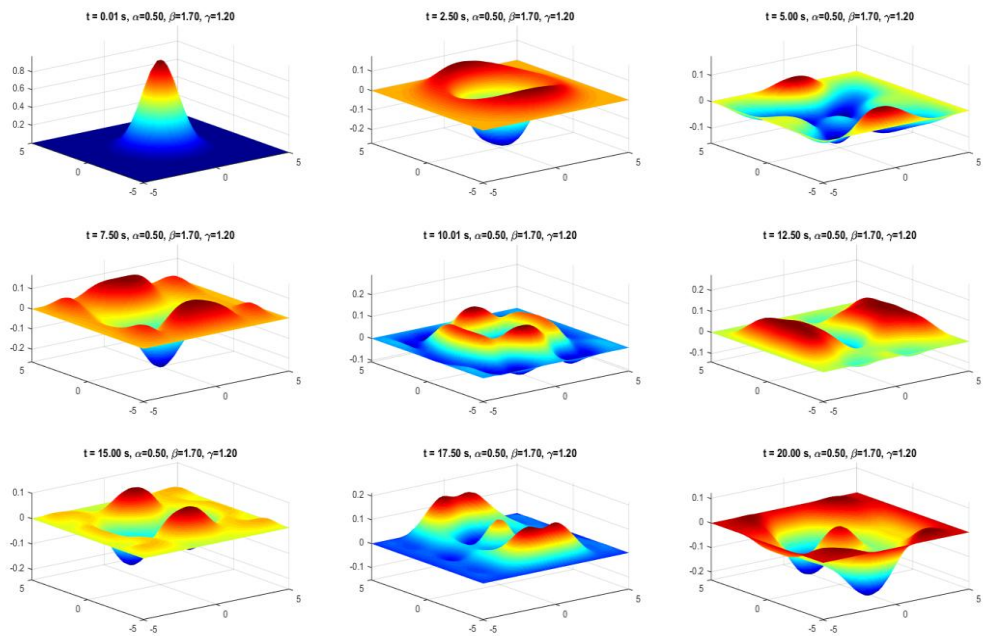

**Figure 44:** Numerical simulation of a Trinition modified wave equation for (0.50, 1.70, 1.20)

# The Minkowski spacetime with 3 deformable and implication to Einstein equations

We recall that in special relativity, the Minskowski metric in coordinates $(t, x, y, z)$ is given by:

$$ds^2_{Mink} = -c^2 dt^2 + dx^2 + dy^2 + dz^2 \tag{240}$$

Using the Trinition norm with one deformable parameter, the above equation was modified as:

$$ds^2_T = -c^2 dt^2 + dx^2 + dy^2 + dz^2 - 2(1-\alpha) dy dz \tag{241}$$

In this section, we now combine the time component with deformed spatial metric introduced in this paper to obtain the full spacetime metric. In the block form, the metric will be given as:

$$g_{\mu\nu} = \begin{pmatrix} -c^2 & 0 \\ 0 & G \end{pmatrix}, \tag{242}$$

We can reformulate out the above expression explicitly as:

$$ds^2_T = -c^2 dt^2 + dx^2 + dy^2 + dz^2 - 2(1-\alpha) dx dy - 2(1-\gamma) dx dz - 2(1-\beta) dy dz \tag{243}$$

This is the suggested new spacetime equation within the Trinition framework with three deformable blocks. We note that the parameters $\alpha, \beta$ and $\gamma$ control the coupling between $x$ and $y$, $x$ and $z$ and $y$ and $z$ directions respectively. There are physical and mathematical implications for this change. The Anisotropic geometry, The off-diagonal terms suggest that the x, y and z directions are no longer mutually orthogonal in the traditional sense. Rather, the principal directions of the space will thus be given by the eigenvalues and eigenvectors of G, and in this case, diagonalizing G will give three eigenvalues that can differ from 1, leading to anisotropy. The modified Einstein equations, In general relativity, the curvature tensors such as the Christoffel symbols and the Ricci tensor would all be affected if the Minkowski metric is replaced by the Trinition-modified one. The field equations are :

$$G_{\mu\nu} = \frac{8\pi G}{c^4} T_{\mu\nu} \tag{244}$$

which would then generate extra factors of those terms or `effective stress-energy contributions' in vacuum, that may be interpreted as encoding anisotropic or exotic gravitational effects. So that has three deformable blocks that correspond to the $xy$, $xz$, and $yz$ couplings. This defines a generalized anisotropic version of the classical Minkowski metric, with deep implications for the geometry of spacetime as well as the dynamics of propagating fields(e.g. waves, gravitational fields, etc.) within it.

## Consequences on curvatures

Now if $\alpha, \beta$ and $\gamma$ are constant, the coefficients of the metric are constants. Then, it can be demonstrated (with the aid of some standard expression) that the entire Riemann metric is flat, i.e., all of its partial derivatives are zero. So the Christoffel spot vanish and the Riemann curvature tensor is zero. That is to say, while the metric looks deformed it has off-diagonal terms, it is flat, this is because any constant metric is locally isometric to Minkowski space in terms of a coordinate transformation that diagonalizes it. However, in the case that any of these deformation parameters are spatially and/or temporally variable (eg:$\alpha = \alpha(x, y, z, t)$), the metric moves away from constancy. This means that the Christoffel symbols are nonzero, and the Riemann curvature tensor (and thus the Ricci and Einstein tensors) will typically also be nonzero in that region. The deformation leads to nontrivial curvature, which can be used to model anisotropic or exotic gravitational effects.

This also have effects on the Christoffel symbols and the Riemann tensor. We recall that the Christoffel symbols are computed from

$$\Gamma^\lambda_{\mu\nu} = \frac{1}{2} g^{\lambda\sigma} \{\partial_\mu g_{\sigma\nu} + \partial_\nu g_{\mu\nu} - \partial_\sigma g_{\mu\nu}\} \tag{245}$$

Here, when the off-diagonal terms like $-2(1-\alpha)dxdy, -2(1-\gamma)dxdz,$ and $-2(1-\beta)dydz$ are present and change with position, the derivatives of $g_{\mu\nu}$ contribute additional terms compared to the standard Minkowski case. Then the Riemann curvature tensor

$$R^\rho_{\sigma\mu\nu} = \partial_\mu \Gamma^\rho_{\nu\sigma} - \partial_\nu \Gamma^\rho_{\mu\sigma} + \Gamma^\rho_{\mu\lambda}\Gamma^\rho_{\nu\sigma} - \Gamma^\rho_{\nu\lambda}\Gamma^\rho_{\mu\sigma} \tag{246}$$

Then acquires nonzero components that reflect the anisotropy introduced by the deformation parameters. The Ricci tensor $R_{\mu\nu}$ and the Einstein tensor $G_{\mu\nu} = R_{\mu\nu} - \frac{1}{2} R g_{\mu\nu}$ will differ from zero even in the regions that would be vacuum in the classical Minkowski case. The following is observed under the three blocks deformation in the Trinition framework in the spatial metric.

***When deformation parameters are constant:*** This metric is a constant deformation of Minkowski space. The fact that it seems non-diagonal turns out not to matter once one realizes that the curvature tensors are vanishing, leaving plenty of room for the connection itself to be flat. The Einstein equations in vacuum are modified in this case, but any ensuing anisotropy is simply a coordinate effect.

***For variable deformation parameters:*** The metric thereby becomes genuinely dynamic, with its nonzero derivatives giving rise to nonvanishing curvature. The Einstein tensor $G_{\mu\nu}$ takes contributions from the spatial deformation(s) it has, which may include model of effective energy-momentum by the geometry itself. This paves the way for new solutions of Einstein's equations in which the anisotropy is woven into the very fabric of spacetime, allowing for potential new gravitational phenomena. These effects are of particular interest for modeling scenarios where spacetime itself might be intrinsically anisotropic or subjected to additional "deformation fields." They could find applications in cosmology like (anisotropic) early universe models, astrophysics like gravity fields around (non-)spherical matter distributions or alternative theories of gravity.

## Variable order deformation factors

An extended explanation of how one can find elegant expressions for the Christoffel symbols and the Einstein tensor whenever the spatial metric is deformed through three independent parameters is below. For our example, we take deformation parameters dependent on a single coordinate (say, $x$) so the metric is function of $x$ only. This simplifies the computation, but nevertheless provides insight into the new features afforded by the Trinitron framework. We will consider a modified spacetime metric of the form

$$ds_T^2 = -c^2 dt^2 + dx^2 + dy^2 + dz^2 - 2(1 - \alpha(x, y, z, t))dxdy - 2(1 - \gamma(t, x, y, z))dxdz - 2(1 - \beta(t, x, y, z))dydz \quad (247)$$

But for simplication purpose, we will work on

$$ds_T^2 = -c^2 dt^2 + dx^2 + dy^2 + dz^2 - 2(1 - \alpha(x))dxdy - 2(1 - \gamma(x))dxdz - 2(1 - \beta(x))dydz \quad (248)$$

Thus, in the coordinate order $(t, x, y, z)$ the metric tensor $g_{\mu\nu}$ is given by:

$$g_{\mu\nu} = \begin{pmatrix} -c^2 & 0 & 0 & 0 \\ 0 & 1 & -[1-\alpha(x)] & -[1-\gamma(x)] \\ 0 & -[1-\alpha(x)] & 1 & -[1-\beta(x)] \\ 0 & -[1-\gamma(x)] & -[1-\beta(x)] & 1 \end{pmatrix}. \tag{249}$$

For the analysis, we will suppose that the deformation parameters $\alpha(x), \beta(x)$ and $\gamma(x)$ depend only on $x$, then, their partial derivative with respect to $y$ and $z$ are zero. The computation of the Christoffel symbols are then presented, we recall that the Christoffel symbols are defined by

$$\Gamma^{\lambda}_{\mu\nu} = \frac{1}{2} g^{\lambda\sigma} \{\partial_{\mu} g_{\sigma\nu} + \partial_{\nu} g_{\mu\nu} - \partial_{\sigma} g_{\mu\nu}\}$$

where $g^{\lambda\sigma}$ is the inverse of $g_{\lambda\sigma}$. In our case, since the deformation parameters depends on $x$, the only nonzero derivatives will be with respect to $x$. For example, if we denote $x = x, y = x^2, z = x^3$, then, we will ignore the time as the metric components are time-independent.

Here for example $g_{11} = 1, g_{12} = -(1-\alpha(x))$, $\partial_1 g_{12} = \alpha'(x)$, but $\partial_1 g_{11} = 0$, thus, for $\Gamma^2_{11}$:

$$\Gamma^2_{11} = g^{22} \alpha'(x)$$

Then, the Riemann curvature tensor is given by

$$R^{\rho}_{\sigma\mu\nu} = \partial_{\mu} \Gamma^{\rho}_{\nu\sigma} - \partial_{\nu} \Gamma^{\rho}_{\mu\sigma} + \Gamma^{\rho}_{\mu\lambda} \Gamma^{\lambda}_{\nu\sigma} - \Gamma^{\rho}_{\nu\lambda} \Gamma^{\lambda}_{\mu\sigma}$$

and the Ricci tensor is obtained by

$$R^{\rho}_{\sigma\mu\nu} = R_{\sigma\nu}$$

Then, the Einstein tensor is the given by:

$$G_{\mu\nu} = R_{\mu\nu} - \frac{1}{2} R g_{\mu\nu}$$

With scalar curvature

$$R = g^{\mu\nu} R_{\mu\nu}$$

This will then also have consequences on the Einstein field equations. We will stop here since we not expect in this field and will let experts to make more investigations.

# Conclusion

In the Transition hypercomplex number concept, we have left behind the cozy realms of classical mathematics and classical physics, opening a new sublime view where geometry, algebra, analysis, and even dynamical systems all twist and turn under the pervasive influence of α. Trinition straddles the line in between establishing itself as a bridge, a continuous spectrum of Algebraic constructs uniting 3D known kingdoms with an emergent 4D behaviour, sliding through this phenomenon based by a lone fractional parameter. It introduces one of the most exciting discoveries of the field: $\alpha$-geometry, a completely generalized metric that distorts Euclidean distance into a pliable norm. In this geometry, our familiar shapes spheres, cylinders, ellipsoids loosen up. Instead, they are free to degeneracy or bloom into new topological forms, at times collapsing along an axis or sprouting forth as something altogether unexpected. By way of this process, this is governed by the commutative versus noncommutative influences of the Trinition product, and it has already explained how one and the same mathematical relationship partial cross-term coupling can give rise to whole families of "cousin" shapes that none of classical geometry could ever have contemplated.

Alongside that, we introduced the *Archiometry*, a new kind of trigonometry designed for this new geometry. In the layers of our classical understanding, where circles and angles govern the behaviour of cos, sin and tan, arises instead the realm of **Archiometry** with Triconus, Trisinus and Trinus functions that contain, still, the spirit of trigonometric identities, like de Moivre's theorem and Pythagorean relations but football more the footprint of partial noncommutativity. It's as if the familiar has been raised to a more flexible dimension, replete with *Archiotwists* (similar to angles) that upend our usual conception of "rotation," and of "phase." Forking into the trigonometric analogs flow and distort, with respect to $\alpha$ smoothly interpolated between near-classical forms at $\alpha = 1$ and exotic cross-coupled behavior for smaller $\alpha$. One of the most striking physically is the Atangana resonance dancers, fractals that emerge generically in Trinition space. These fractals "dance" because they can flicker between known fractal forms and fabulously intricate shapes we'd never seen previously, all resonating under the partial commutation. Adjusting $\alpha$ causes these fractals to continuously evolve over time, a demonstration that even fractal geometry can be changed into an infinitude of successively subtle states, thus revealing continuities that classical fractal families always obscure. In doing so, Trinition asks us to think of fractals, not as static sets, but as dynamic or resonant structures, living structures that "dance" in the presence of cross terms as function of alpha. Finally, on the physical side, the redefinition of equations of motion, Hamiltonians, Lagrangian and energy clearly states that Trinition is no mere algebraic curio: it has great potential in terms of modeling wave propagation, rotational symmetries or quantum-like systems in a partially noncommutative environment. The aspect of introducing Trinition transforms that make Laplace and Fourier invariants as transformations of the $x, y$ space extended to handily configure and transform not just the square norm but the cross term in it, heads us in provided/deployed, to handle a more numeric piece of cross geometry as well analytically be it pure PDE, wave property, or wave-type undulating through the classics and new realms at ease. Trinition, run of the mill, is much more than a tweak of hypercomplex arithmetic. It is a very broad generalization of our understanding of all of geometry, analysis, fractals, physics and so on. We can therefore open and close the tap of $\alpha$, smoothly transitioning between quasi-classical three-dimensional bounds and extravagant four-dimensional noncommutative bouquets, revealing a kaleidoscope of shapes, fractals, and dynamics. In this, Trinition does not simply unify existing structures; it supersedes them, providing a sandbox in which Euclid, Pythagoras and de Moivre can be endeavoured into untraveled territory. It is a lighthouse reminding us of the infinite ability of mathematics and physics to reformulate themselves to unlock the doors to things we never even dream that we will fully comprehend through our classical eyes. This paper significantly contributed to further development of the Trinition framework, as it moves Trinition from theoretical framework into real-world applications, and generalize the model norm. Our contributions are twofold: 1) By injecting deformation parameters into the inner product we have defined a valid differential structure which offers a dynamical variation of the mathematical structures of classical differential operators. 2) We subsequently prove their utility for generating generalizations of the classical wave and heat equations that exhibit anisotropic and non-Euclidean features. These tools from analysis make it possible to formulate the deformed Sobolev spaces and derive generalized Sobolev and Gagliardo–Nirenberg inequalities, leading to a solid functional-analytic framework that guarantees the well-posedness and regularity of solutions in a new setting. Furthermore, as our generalization of the Laplace transform with respect to the newly defined Trinition metric suggests, it provides an efficient way of generating probability distribution functions with tunable, deformable shapes. Full-scale numerical simulations confirm our theoretical predictions, providing a vivid analysis of the effects of the deformation parameters on the dynamics of the wave propagation and the heat diffusion process. In conclusion, these findings highlight the approach's potential for new horizons in the solution of complex problems in mathematical physics and engineering, especially for anisotropic and non-Euclidean models that extend the domain of applicability of the classical frameworks.